\documentclass{amsart}
\usepackage[T1]{fontenc}
\usepackage[latin1]{inputenc}
\usepackage{fancyhdr}
\usepackage{graphics}
\makeindex

\makeatletter

\providecommand{\LyX}{L\kern-.1667em\lower.25em\hbox{Y}\kern-.125emX\@}

 \theoremstyle{plain}    
 \newtheorem{thm}{Theorem}[section]
 \numberwithin{equation}{section} 
 \numberwithin{figure}{section} 
 \theoremstyle{plain}    
 \newtheorem{cor}[thm]{Corollary} 
 \theoremstyle{plain}    
 \newtheorem{lem}[thm]{Lemma} 
 \theoremstyle{plain}    
 \newtheorem{prop}[thm]{Proposition} 
 \theoremstyle{definition}
 \newtheorem{defn}[thm]{Definition}
 \theoremstyle{definition}
  \newtheorem{example}[thm]{Example}
 \theoremstyle{remark}
 \newtheorem{rem}[thm]{Remark}
 \theoremstyle{remark}    
 \newtheorem{claim}[thm]{Claim}
 \theoremstyle{remark}    
 \newtheorem*{claim*}{Claim}

\usepackage{verbatim} 
\usepackage{varioref}
\usepackage{amssymb}
\usepackage[mathscr]{eucal}
\usepackage{amsthm}
\usepackage[all]{xy}
\usepackage{hyperref}

\mathcode`\:="603A 

\DeclareSymbolFont{rsfs}{U}{rsfs}{m}{n}
\DeclareSymbolFontAlphabet{\mathrf}{rsfs}

\newcommand{\tmap}[1]{{\mathrm T}_{#1}}

\newcommand{\rank}{\operatorname{rank}}
\newcommand{\rn}[1]{{\mathfrak R}_n(#1)}
\newcommand{\sshufp}[1]{\mathscr{Z}\bigl\{#1\bigr\}}

\newcommand{\im}{\operatorname{im}}

\newcommand{\ainfaction}[1]{\hat{a}_{#1}}

\newcommand{\glist}[3]{#1_{#2},\dots,#1_{#3}}
\newcommand{\blist}[2]{\glist{#1}{1}{#2}}

\newcommand{\tlist}[2]{\tmap{\blist{#1}{#2}}}

\newcommand{\fface}{\mathrm{F}_0}
\newcommand{\lface}{\tilde{\mathrm{F}}}
\newcommand{\nth}[1]{$#1^{\mathrm{th}}$}
\newcommand{\tunder}[2]{\tmap{\underbrace{\scriptstyle#1}_{{\text{#2}}}}}

\newcommand{\tunderi}[2]{\tunder{1,\dots,#1,\dots,1}{\nth #2\ position}}

\DeclareFontFamily{U}{manual}{}
\DeclareFontShape{U}{manual}{m}{n}{ <->  manfnt }{}
\newcommand{\manfntsymbol}[1]{%
    {\fontencoding{U}\fontfamily{manual}\selectfont\symbol{#1}}}

\newcommand{\dbend}{\manfntsymbol{127}}
\newcommand{\lhdbend}{\manfntsymbol{126}}
\newcommand{\reversedvideodbend}{\manfntsymbol{0}}
\newcommand{\textdbend}{\text@dbend{\dbend}}
\newcommand{\textlhdbend}{\text@dbend{\lhdbend}}
\newcommand{\textreversedvideodbend}{\text@dbend{\reversedvideodbend}}
\newlength{\dbend@height}
\newcommand{\text@dbend}[1]{%
  \settoheight{\dbend@height}{#1}%
  \divide\dbend@height by 15%
  \multiply\dbend@height by 22%
  \raisebox{\dbend@height}{#1}}

\CompileMatrices

\SelectTips{cm}{}
\newdir{ >}{{}*!/-5pt/@{>}}
\makeatother

\begin{document}


\title{Operads and algebraic homotopy}

\author{Justin R. Smith}



\address{Department of Mathematics and Computer Science\\
Drexel University\\
Philadelphia,~PA 19104}

\email{jsmith@mcs.drexel.edu}

\urladdr{http://www.mcs.drexel.edu/\textasciitilde{}jsmith}

\date{\today}

\begin{abstract}
This paper proves that the homotopy type of a simply-connected, pointed simplicial
set is determined by the chain-complex augmented with functorial diagonal and
higher diagonal maps similar to those used to define Steenrod operations. The
treatment is entirely self-contained --- simplifying, extending, and correcting
results from the author's monograph ``Iterating the cobar construction,''
AMS Memoirs, volume 109.
\end{abstract}
\maketitle
\tableofcontents{}

\newcommand{\integers}{\mathbb{Z}}

\newcommand{\betabar}{\bar{\beta }}
 
\newcommand{\desusp}{\downarrow }

\newcommand{\susp}{\uparrow }

\newcommand{\cobar}{\mathcal{F}}

\newcommand{\mfrac}{\mathfrak{M}}

\newcommand{\coend}{\mathrm{CoEnd}}

\newcommand{\ainfty}{A_{\infty }}

\newcommand{\coassoc}{\mathrm{Coassoc}}

\newcommand{\trm}{\mathrm{T}}

\newcommand{\tfr}{\mathfrak{T}}

\newcommand{\tabbr}{\hat{\trm }}

\newcommand{\Tabbr}{\hat{\tfr }}

\newcommand{\afr}{\mathfrak{A}}

\newcommand{\homz}{\mathrm{Hom}_{\integers }}

\newcommand{\zend}{\mathrm{End}}

\newcommand{\rs}[1]{\mathrm{R}S_{#1 }}

\newcommand{\highprod}[1]{\bar{\mu }_{#1 }}

\newcommand{\barcs}{\bar{\mathcal{B}}}

\newcommand{\ubarcs}{\mathcal{B}}

\newcommand{\zs}[1]{\mathbb{Z}S_{#1 }}

\newcommand{\homzs}[1]{\mathrm{Hom}_{\integers S_{#1 }}}

\newcommand{\zpi}{\mathbb{Z}\pi }

\newcommand{\D}{\mathfrak{D}}

\newcommand{\ahat}{\hat{\mathfrak{A}}}

\newcommand{\cbar}{\bar{C}}

\newcommand{\cf}[1]{\mathscr{C}(#1 )}

\newcommand{\ddelta}{\dot{\Delta }}

\section{Introduction}

\newdir{ >}{{}*!/-5pt/@{>}}

The easiest way to convey the flavor of this paper's results is with a simple
example. Suppose \( X \) and \( Y \) are pointed, simply-connected, 2-reduced
simplicial sets. There are many topological invariants associated with the chain-complexes
of \( X \) and \( Y \) --- including the coproduct and \emph{\( S_{n} \)-}equivariant
\emph{higher coproducts} (used to define Steenrod operations):
\begin{eqnarray*}
\rs{{n}}\otimes C(X) & \rightarrow  & C(X)^{n}\\
\rs{{n}}\otimes C(Y) & \rightarrow  & C(Y)^{n}
\end{eqnarray*}
for all \( n>1 \), where:

\begin{enumerate}
\item \( \rs{{n}} \) is the bar-resolution of \( \integers  \) over \( \zs{{n}} \).
\item \( (*)^{n} \) denotes the \( n \)-fold tensor product over \( \integers  \)
(with \( S_{n} \) acting by permuting factors).
\end{enumerate}
Now suppose we know (from a purely \emph{algebraic} analysis of these \emph{chain-complexes})
that there exists a chain-map inducing homology isomorphisms
\[
f:C(X)\rightarrow C(Y)\]
and making \begin{equation}\xymatrix{{\rs{n}\otimes C(X)} \ar[r]\ar[d]_{1\otimes f}& {C(X)^n} \ar[d]^{f^n}\\ {\rs{n}\otimes C(Y)} \ar[r]& {C(Y)^n}}\label{dia:intro1}\end{equation}
commute for all \( n>1 \) (requiring exact commutativity is unnecessarily restrictive,
but we assume it to simplify this discussion).

Then our theory asserts that \( X \) and \( Y \) are homotopy equivalent via
a geometric map inducing a map chain-homotopic to \( f \). In fact, our theory
goes further than this: if \( f \) is \emph{any} chain-map making \ref{dia:intro1}
commute for all \( n>1 \), there exists a map of topological realizations
\[
F:|X|\rightarrow |Y|\]
inducing a map equivalent to \( f \) (see theorems~\ref{th:equiv} and \ref{th:morphism}
for the exact statements). Essentially, Steenrod operations on the \emph{chain-level}
determine: 

\begin{itemize}
\item the homotopy types of spaces and maps,
\item \emph{all} obstructions to topologically realizing chain-maps.
\end{itemize}
I view this work as a generalization of Quillen's characterization of rational
homotopy theory in \cite{Quillen:1969} \index{Quillen@\textsc{D. Quillen}}---
he showed that rational homotopy types are determined by a commutative coalgebra
structures on chain-complexes; we show that integral homotopy types are determined
by coalgebra structures augmented by higher diagonals.

In \S~\ref{sec:operadmcoalg}, we define operads and m-coalgebras. Operads are
templates for ``generic'' algebraic structures --- for instance, one can easily
define commutative algebras and coalgebras, Lie algebras, and so forth, using
suitable operads (see \cite{Kriz-May}). 

M-coalgebras are coalgebras over certain types of operads --- essentially, m-coalgebras
are chain-complexes equipped with higher diagonals as in the example above.
The concept of m-coalgebra encapsulates this diagonal and higher diagonal data.

In \S~\ref{sec:sfrakdef}, we define a particular operad, \( \mathfrak{S} \),
with topological significance in our theory. We also define a functor \( \cf{{*}} \)
that associates an m-coalgebra (over the operad \( \mathfrak{S} \)) to a simply-connected
pointed 2-reduced simplicial set. This m-coalgebra is nothing but the chain-complex
of the space equipped with the natural higher diagonal and diagonal maps used
in defining cup products and Steenrod operations.

In \S~\ref{sec:mcoalgehomotopy}, we define a category of m-coalgebras, \( \mathfrak{L}_{0} \)
and a localization of it, \( \mathfrak{L} \) by a class of maps called elementary
equivalences. This gives rise to a \emph{homotopy theory} of m-coalgebras in
terms of which we can state this paper's main result:\par{} \vspace{0.3cm}

\emph{Theorem~\ref{cor:lrealizable} and \ref{cor:mainfiniteresult}:The functor
\[
\mathscr{C}(*):\underline{\mathrm{Homotop}}_{0}\rightarrow \mathfrak{L}^{+}\]
 defines an equivalence of categories and homotopy theories (in the sense of
\cite{Quillen:1967}), where \( \underline{\mathrm{Homotop}}_{0} \) is the
homotopy category of pointed, simply-connected CW-complexes and continuous maps
and \( \mathfrak{L}^{+}\subset \mathfrak{L} \) is the subcategory of topologically
realizable m-coalgebras. In addition, there exists an equivalence of categories
and homotopy theories
\[
\cf{{*}}:\mathrf {F}\rightarrow \mathfrak{F}^{+}\]
where \( \mathrf {F} \) is the homotopy category of finite, pointed, simply-connected
simplicial sets and \( \mathfrak{F}^{+} \) is the homotopy category of finitely
generated, topologically realizable m-coalgebras in \( \mathfrak{L}_{0} \)
localized with respect to finitely generated equivalences in \( \mathfrak{L}_{0} \).} \par{} \vspace{0.3cm}

In the spirit of our initial statement, one corollary to this result is

\emph{Corollary \ref{cor:mainnoncat}: Let \( X \) and \( Y \) be pointed,
simply-connected semisimplicial sets and let 
\[
f:C(X)\rightarrow C(Y)\]
be a chain-map between canonical chain-complexes. Then \( f \) is topologically
realizable if and only if there exists an m-coalgebra \( C \) over \( \mathfrak{S} \)
and a factorization \( f=f_{\beta }\circ f_{\alpha } \) \[\xymatrix{{\cf{X_1}}\ar[r]^-{f_{\alpha}}&{C}\ar@<.5ex>@{ .>}[r]^-{f_{\beta}}&{\cf{X_2}}\ar@<.5ex>@{ >-}[l]^-{\iota}}\]
where \( f_{\alpha } \) is a morphism of m-coalgebras, \( \iota  \) is an}
elementary equivalence \emph{--- an injection of m-coalgebras with acyclic,
\( \integers  \)-free cokernel --- and \( f_{\beta } \) is a chain map that
is a left inverse to \( \iota  \). If \( X \) and \( Y \) are finite, we
may require \( C \) to be finitely generated.}

\par{} \vspace{0.3cm}

The rest of the paper builds the machinery necessary to prove this result. The
most important piece is the \emph{cobar construction} of an m-coalgebra. The
cobar construction was defined by J.~F.~Adams in \cite{Adams:1956} \index{Adams@\textsc{J. F. Adams}}and
it determines the chain complex of the loop space of a space. 

We show that every m-coalgebra, \( C \), has a well-defined cobar construction,
\( \cobar C \), that comes equipped with a well-defined (see \ref{cor:mainresult})
and topologically valid (see \ref{theorem:mapfiber}) m-coalgebra structure.
Although this was stated in \cite{Smith:1994}, the present treatment uses more
standard notation (i.e., operads) and simplifies (and, in some cases, \emph{corrects})
the proofs in that paper.

One of the most interesting aspects of this research is that the cobar construction
seems to be the \emph{key} to homotopy theory --- understanding it leads quickly
to our main theorem.

In \S~\ref{sec:ainfinity}, we define \( \ainfty  \)-structures on algebras
and coalgebras. Roughly speaking, \( \ainfty  \)-structures are algebra (or
coalgebra) structures that are \emph{homotopy associative}. All m-coalgebras
come equipped with a canonical \( \ainfty  \)-coalgebra structure.

In \S~\ref{sec:barcobarconstructions}, we give a slightly nonstandard definition
of the bar and cobar constructions in terms of iterated algebraic mapping cones.
This nonstandard definition facilitates the proof of an important duality theorem
(proposition~\ref{prop:bardualcobar}).

In \S~\ref{sec:cobar}, we compute our m-coalgebra structure of the cobar construction
and the canonical acyclic twisted tensor product with fiber equal to the cobar
construction. In \S~\ref{sec:geometricity}, we prove that these constructs
are are topologically valid.

This paves the way for the main results in \S~\ref{ch:equivcategories}.

It is necessary to compare our results those of Smirnov (in his remarkable paper
\cite{Smirnov:1985}) \index{Smirnov@\textsc{V. A. Smirnov}} and several of
his other papers. In \cite{Smirnov:1985}, he showed that the homotopy type
of a space was determined by its chain-complex augmented with a coalgebra structure
over an operad \( E^{*} \). Smirnov's structure is more complex than ours:
his operad's components, \( E(n) \), are uncountably generated in all dimensions
and for all values of \( n>1 \). Smirnov's ``higher diagonal'' map contains
so much data it is almost surprising that it only determines a space's homotopy
type.

Although there is no obvious connection between Smirnov's invariant and more
commonly used homotopy invariants (such as Steenrod operations, coproducts,
etc.), he develops a connection with Steenrod squares in \S~4 of \cite{Smirnov:1985}
and other structures of the Steenrod algebra in \cite{Smirnov:1992-2}. It would
be interesting to elucidate the precise relationship between Smirnov's work
and ours.

I feel it is also appropriate to compare and contrast my results with work of
Michael Mandell. In \cite{mandell:1998-1}, he proved

\medskip \textbf{Main Theorem.} \emph{The singular cochain functor with coefficients
in} \( \bar{\integers }_{p} \) \emph{induces a contravariant equivalence from
the homotopy category of connected nilpotent} \( p \)-\emph{complete spaces
of finite \( p \)-type to a full subcategory of the homotopy category of} \( E_{\infty } \)
\( \bar{\integers }_{p} \)-\emph{algebras}. \medskip

Here, \( p \) denotes a prime and \( \bar{\integers }_{p} \) the algebraic
closure of the finite field of \( p \) elements. \( E_{\infty } \)-algebras
are defined in \cite{Kriz-May} --- they are modules over a suitable operad.

Since Dr.~Mandell characterizes nilpotent \( p \)-complete spaces in terms
of \( E_{\infty } \) \( \bar{\integers }_{p} \)-algebras, his results appear
to be dual to mine. This is not the case, however. A complete characterization
of nilpotent \( p \)-complete spaces does not lead to one of integral homotopy
types: One must somehow know that \( p \)-local homotopy equivalences patch
together. Consequently, his results do not imply mine.

The converse statement is also true: My results do not imply his.

My results in \cite{Smith:1994} and the present paper imply that all the primes
``mix'' when one studies algebraic properties of homotopy theory (for instance
the \( p \)-local structure of the cobar construction of a space depends on
the \( q \)-local structure of the space for all primes \( q\geq p \)). This
is intuitively clear when considers the composite \( (1\otimes \Delta )\circ \Delta  \)
(iterated coproducts) and notes that \( \integers _{2} \) acting on both copies
of \( \Delta  \) give rise to elements of the symmetric group on 3 elements.

Consequently, a characterization of integral homotopy does not lead to a \( p \)-local
homotopy theory: In killing off all primes other than \( p \), one also kills
off crucial information needed to compute the cobar construction of a space.

In \cite{mandell:1998-1}, Dr.~Mandell proved that one \textit{must} pass to
the algebraic closure of \( \integers _{p} \) to get a characterization of
\( p \)-complete homotopy theory. I conjecture that, in passing to the algebraic
closure, one kills off additional data within the homotopy type --- namely the
data that depends on larger primes. Consequently, one restores algebraic consistency
to the theory, regaining the ability to characterize local homotopy types.

I would like to express my appreciation to Jim Stasheff for acquainting me with
operads and for his encouragement. Some of the material in this paper grew out
of a series of six lectures I gave in Jim Stasheff's Deformation Theory Seminar
at the University of Pennsylvania.

I would also like to thank Peter May and Mark Mahowald for their encouragement.

I am indebted to Drexel University for funding my attendance at the JAMI conference
at Johns Hopkins University on Recent Developments in Homotopy Theory.

\section{Operads and m-coalgebras\label{sec:operadmcoalg}}

\subsection{Definitions}

\fancyhead[RO,LE]{\rightmark} \fancyhead[RE,LO]{Justin R. Smith}Many of these
definitions appeared in \cite{Smith:1994}. Unfortunately, that paper used a
very nonstandard notation. Since that paper appeared, the author has standardized
and simplified most of these definitions. See~\ref{app:concordance} for a comparison
of the two systems of notation.

\begin{defn}
\label{def:degmap} Let \( C \) and \( D \) be two graded \( \integers  \)-modules.
A map of graded modules \( f:C_{i}\rightarrow D_{i+k} \) will be said to be
of degree \( k \). 
\end{defn}
\begin{rem}
For instance the \textit{differential} of a chain-complex will be regarded as
a degree \( -1 \) map. 
\end{rem}
We will make extensive use of the Koszul\index{Koszul@\textsc{J. L. Koszul}}
Convention\index{Koszul Convention}\index{convention, Koszul} (see~\cite{Gugenheim:1960}\index{Gugenheim@\textsc{Victor K. A. M. Gugenheim}})
regarding signs in homological calculations:

\begin{defn}
\label{def:basic} \label{def:koszul} If \( f:C_{1}\rightarrow D_{1} \), \( g:C_{2}\rightarrow D_{2} \)
are maps, and \( a\otimes b\in C_{1}\otimes C_{2} \) (where \( a \) is a homogeneous
element), then \( (f\otimes g)(a\otimes b) \) is defined to be \( (-1)^{\deg (g)\cdot \deg (a)}f(a)\otimes g(b) \). 
\end{defn}
\begin{rem}
This convention simplifies many of the common expressions that occur in homological
algebra --- in particular it eliminates complicated signs that occur in these
expressions. For instance the differential, \( \partial _{\otimes } \), of
the tensor product \( \partial _{C}\otimes 1+1\otimes \partial _{D} \).

Throughout this entire paper we will follow the convention that group-elements
act on the left. Multiplication of elements of symmetric groups will be carried
out accordingly --- i.e. 
\[
\left( \begin{array}{cccc}
1 & 2 & 3 & 4\\
2 & 3 & 1 & 4
\end{array}\right) \cdot \left( \begin{array}{cccc}
1 & 2 & 3 & 4\\
4 & 3 & 2 & 1
\end{array}\right) =\]
 result of applying \( \left( \begin{array}{cccc}
1 & 2 & 3 & 4\\
2 & 3 & 1 & 4
\end{array}\right)  \) after applying \( \left( \begin{array}{cccc}
1 & 2 & 3 & 4\\
4 & 3 & 2 & 1
\end{array}\right)  \) or n\( \left( \begin{array}{cccc}
1 & 2 & 3 & 4\\
4 & 1 & 3 & 2
\end{array}\right)  \).

Let \( f_{i} \), \( g_{i} \) be maps. It isn't hard to verify that the Koszul
convention implies that \( (f_{1}\otimes g_{1})\circ (f_{2}\otimes g_{2})=(-1)^{\deg (f_{2})\cdot \deg (g_{1})}(f_{1}\circ f_{2}\otimes g_{1}\circ g_{2}) \).

We will also follow the convention that, if \( f \) is a map between chain-complexes,
\( \partial f=\partial \circ f-(-1)^{\deg (f)}f\circ \partial  \). The compositions
of a map with boundary operations will be denoted by \( \partial \circ f \)
and \( f\circ \partial  \) --- see \cite{Gugenheim:1960}\index{Gugenheim@\textsc{Victor K. A. M. Gugenheim}}.
This convention clearly implies that \( \partial (f\circ g)=(\partial f)\circ g+(-1)^{\deg (f)}f\circ (\partial g) \).
We will call any map \( f \) with \( \partial f=0 \) a chain-map. We will
also follow the convention that if \( C \) is a chain-complex and \( \susp :C\rightarrow \Sigma C \)
and \( \desusp :C\rightarrow \Sigma ^{-1}C \) are, respectively, the suspension
and desuspension maps, then \( \susp  \) and \( \desusp  \) are both chain-maps.
This implies that the boundary of \( \Sigma C \) is \( -\susp \circ \partial _{C}\circ \desusp  \)
and the boundary of \( \Sigma ^{-1}C \) is \( -\desusp \circ \partial _{C}\circ \susp  \).
\end{rem}
\begin{defn}
\label{r:koszul.5} We \index{T-map}will use the symbol \( T \) to denote
h transposition operator for tensor products of chain-complexes \( T:C\otimes D\rightarrow D\otimes C \),
where \( T(c\otimes d)=(-1)^{\dim (c)\cdot \dim (d)}d\otimes c \).
\end{defn}
\begin{prop}
\label{prop:suspisos}Let \( C \) and \( D \) be chain-complexes\index{suspension isomorphisms}\index{isomorphisms, suspension}.
Then there exist isomorphisms
\[
L_{k}=\desusp ^{k}_{C\otimes D}\circ (\susp ^{k}_{C}\otimes 1_{D}):\Sigma ^{-k}C\otimes D\rightarrow \Sigma ^{-k}(C\otimes D)\]
sending \( c\otimes d\in \Sigma ^{-k}C_{i}\otimes D_{j} \) to \( c\otimes d\in \Sigma ^{-k}(C\otimes D)_{i+j} \),
r , \( d\in D_{j} \), and 
\[
M_{k}=\desusp ^{k}_{C\otimes D}\circ (1_{C}\otimes \susp ^{k}_{D}):C\otimes \Sigma ^{-k}D\rightarrow \Sigma ^{-k}(C\otimes D)\]
sending \( c\otimes d\in C_{i}\otimes \Sigma ^{-k}D_{j} \) to \( (-1)^{ik}c\otimes d\in \Sigma ^{-k}(C\otimes D)_{i+j} \),
for \( c\in C_{i} \) and \( d\in \Sigma ^{-k}D_{j}=D_{j+k} \).
\end{prop}
\begin{defn}
\label{def:tmap} Let \( \alpha _{i} \), \( i=1,\dots ,n \) be a sequence
of nonnegative integers whose sum is \( |\alpha | \). Define a set-mapping
 symmetric groups 
\[
\tlist {\alpha }{n}:S_{n}\rightarrow S_{|\alpha |}\]
 as follows:
\begin{enumerate}
\item for \( i \) between 1 and \( n \), let \( L_{i} \) denote the length-\( \alpha _{i} \)
integer sequence: 
\item ,where \( A_{i}=\sum _{j=1}^{i-1}\alpha _{j} \) --- so, for instance, the concatenation
of all of the \( L_{i} \) is the sequence of integers from 1 to \( |\alpha | \); 
\item \( \tlist {\alpha }{n}(\sigma ) \) is the permutation on the integers \( 1,\dots ,|\alpha | \)
that permutes the blocks \( \{L_{i}\} \) via \( \sigma  \). In other words,
 \( \sigma  \) s the permutation 
\[
\left( \begin{array}{ccc}
1 & \dots  & n\\
\sigma (1) & \dots  & \sigma (n)
\end{array}\right) \]
 then \( \tlist {\alpha }{n}(\sigma ) \) is the permutation defined by writing
\[
\left( \begin{array}{ccc}
L_{1} & \dots  & L_{n}\\
L_{\sigma (1)} & \dots  & L_{\sigma (n)}
\end{array}\right) \]
 and regarding the upper and lower rows as sequences length \( |\alpha | \). 
\end{enumerate}
\end{defn}
\begin{rem}
Do not confuse the \( T \)-maps defined here with the transposition map for
tensor products of chain-complexes. We will use the special notation \( T_{i} \)
to represent \( T_{1,\dots ,2,\dots ,1} \), where the 2 occurs in the \( i^{\mathrm{th}} \)
position. The two notations don't conflict since the old notation is never used
in the case when \( n=1 \). Here is an example of the computation of \( \tmap {2,1,3}((1,3,2))=\tmap {2,1,3}\left( \begin{array}{ccc}
1 & 2 & 3\\
3 & 1 & 2
\end{array}\right)  \):\( L_{1}=\{1\}2 \), \( L_{2}=\{3\} \), \( L_{3}=\{4,5,6\} \). The permutation
maps the ordered set \( \{1,2,3\} \) to \( \{3,1,2\} \), so we carry out the
corresponding mapping of the sequences \( \{L_{1},L_{2},L_{3}\} \) to get \( \left( \begin{array}{ccc}
L_{1} & L_{2} & L_{3}\\
L_{3} & L_{1} & L_{2}
\end{array}\right) =\left( \begin{array}{ccc}
\{1,2\} & \{3\} & \{4,5,6\}\\
\{4,5,6\} & \{1,2\} & \{3\}
\end{array}\right) =\left( \begin{array}{cccccc}
1 & 2 & 3 & 4 & 5 & 6\\
4 & 5 & 6 & 1 & 2 & 3
\end{array}\right)  \) (or \( ((1,4)(2,5)(3,6)) \), in cycle notation).
\end{rem}
\begin{defn}
\label{def:operad} A sequence of differential graded \( \mathbb{Z} \)-free
modules, \( \{\mathscr{U}_{i}\} \), will be said to form an \emph{operad} if
they satisfy the following conditions\index{operad!definition}\index{definition!operad}:
\begin{enumerate}
\item there exists a \emph{unit map} \index{unit map}\index{operad!unit map}(defined
by the commutative diagrams below) 
\[
\eta :\mathbb{Z}\rightarrow \mathscr{U}_{1}\]

\item for all \( i>1 \), \( \mathscr{U}_{i} \) is equipped with a left action of
\( S_{i} \), the symmetric group. 
\item for all \( k\geq 1 \), and \( i_{s}\geq 0 \) there are \index{operad!$\gamma$-map}\index{g@$\gamma$-map in an operad}maps
\[
\gamma :\mathscr{U}_{i_{1}}\otimes \cdots \otimes \mathscr{U}_{i_{k}}\otimes \mathscr{U}_{k}\rightarrow \mathscr{U}_{i}\]
 where \( i=\sum _{j=1}^{k}i_{j} \).

The \( \gamma  \)-maps must satisfy the following conditions:

\end{enumerate}
\end{defn}
\begin{description}
\item [Associativity]the following \index{operad!associativity identity}\index{associativity identity of an operad}diagrams
commute, where \( \sum j_{t}=j \), \( \sum i_{s}=i \), and \( g_{\alpha }=\sum _{\ell =1}^{\alpha }j_{\ell } \)
and \( h_{s}=\sum _{\beta =g_{s-1}+1}^{g_{s}}i_{\beta } \):  \[\xymatrix@C+20pt{{\left(\bigotimes_{s=1}^{j}\mathscr{U}_{i_{s}}\right)\otimes\left(\bigotimes_{t=1}^{k}\mathscr{U}_{j_{t}}\right) \otimes\mathscr{U}_{k}}\ar[r]^-{\text{Id}\otimes\gamma}\ar[dd]_{\text{shuffle}}&{\left(\bigotimes_{s=1}^{j}\mathscr{U}_{i_{s}}\right)\otimes\mathscr{U}_{j}}\ar[d]^{\gamma}\\
&{\mathscr{U}_{i}}\\{\left(\left(\bigotimes_{q=1}^{j_{t}}\mathscr{U}_{i_{g_{t-1}+q}}\right)\otimes\bigotimes_{t=1}^{k}\mathscr{U}_{j_{t}}\right) \otimes\mathscr{U}_{k}}\ar[r]_-{(\otimes_{t}\gamma)\otimes\text{Id}}&{\left(\bigotimes_{t=1}^{k}\mathscr{U}_{h_{k}}\right) \otimes\mathscr{U}_{k}}\ar[u]_{\gamma}}\]

\item [Units]the \index{operad!unit identity}following diagrams commute:  \[\begin{array}{cc}\xymatrix{{{\integers}^{k}\otimes\mathscr{U}_{k}}\ar[r]^{\cong}\ar[d]_{{\eta}^{k}\otimes\text{Id}}&{\mathscr{U}_{k}}\\
{{\mathscr{U}_{1}}^{k}\otimes{\mathscr{U}_{k}}}\ar[ur]_{\gamma}&}&\xymatrix{{\mathscr{U}_{k}\otimes\integers}\ar[r]^{\cong}\ar[d]_{\text{Id}\otimes\eta}&{\mathscr{U}_{k}}\\
{\mathscr{U}_{k}\otimes\mathscr{U}_{1}}\ar[ur]_{\gamma}&}\end{array}\] 
\item [Equivariance]the \index{operad!equivariance}following diagrams commute: \[\xymatrix@C+20pt{{\mathscr{U}_{j_{1}}\otimes\cdots\otimes\mathscr{U}_{j_{k}}\otimes\mathscr{U}_{k}}\ar[r]^-{\gamma}\ar[d]_{\sigma^{-1}\otimes\sigma}&{\mathscr{U}_{j}}\ar[d]^{\tmap{j_{1},\dots,j_{k}}(\sigma)}\\
{\mathscr{U}_{j_{\sigma(1)}}\otimes\cdots\otimes\mathscr{U}_{j_{\sigma(k)}}\otimes\mathscr{U}_{k}}\ar[r]_-{\gamma}&{\mathscr{U}_{j}}}\]
where \( \sigma \in S_{k} \), and the \( \sigma ^{-1} \) on the left permutes
the factors \( \{\mathscr{U}_{j_{i}}\} \) and the \( \sigma  \) on the right
simply acts on \( \mathscr{U}_{k} \). See \ref{def:tmap} for a definition
of \( \tmap {j_{1},\dots ,j_{k}}(\sigma ) \). \[\xymatrix@C+20pt{{\mathscr{U}_{j_{1}}\otimes\cdots\otimes\mathscr{U}_{j_{k}}\otimes\mathscr{U}_{k}}\ar[r]^-{\gamma}\ar[d]_{\tau_{1}\otimes\cdots\tau_{k}\otimes\text{Id}}&{\mathscr{U}_{j}}\ar[d]^-{\tau_{1}\oplus\cdots\oplus\tau_{k}}\\
{\mathscr{U}_{j_{\sigma(1)}}\otimes\cdots\otimes\mathscr{U}_{j_{\sigma(k)}}\otimes\mathscr{U}_{k}}\ar[r]_-{\gamma}&{\mathscr{U}_{j}}}\] where
\( \tau _{s}\in S_{j_{s}} \) and \( \tau _{1}\oplus \cdots \oplus \tau _{k}\in S_{j} \)
is the block sum.
\end{description}
\begin{rem}
The alert reader will notice a discrepancy between our definition of operad
and that in \cite{Kriz-May} \index{Kriz@\textsc{I. Kriz}}\index{May@\textsc{P. J. May}}(on
which it was based). The difference is due to our using operads as parameters
for systems of \emph{maps}, rather than \( n \)-ary operations. We, consequently,
compose elements of an operad as one composes \emph{maps}, i.e. the second operand
is to the \emph{left} of the first. This is also why the symmetric groups act
on the \emph{left} rather than on the right. 
\end{rem}
We will frequently want to think of operads in other terms:

\begin{defn}
\label{def:operadcomps} Let \( \mathscr{U} \) be an operad, \index{operad!compositions}\index{composition algebras}as
defined above. Given \( k_{1}\geq k_{2}>0 \), define the \( i^{\mathrm{th}} \)
\emph{composition}

\[
\circ _{i}:\mathscr{U}_{k_{2}}\otimes \mathscr{U}_{k_{1}}\rightarrow \mathscr{U}_{k_{1}+k_{2}}\]
 as the composite \begin{multline}\underbrace{\integers\otimes\cdots\otimes\integers\otimes\mathscr{U}_{k_{2}}\otimes\integers\otimes\cdots\otimes\integers}_{\text{\(i^{\text{th}}\)factor}}\otimes\mathscr{U}_{k_{1}}\\ \to\underbrace{\mathscr{U}_{1}\otimes\cdots\otimes\mathscr{U}_{1}\otimes\mathscr{U}_{k_{2}}\otimes\mathscr{U}_{1}\otimes\cdots\otimes\mathscr{U}_{1}}_{\text{\(i^{\text{th}}\)factor}}\otimes\mathscr{U}_{k_{1}}\to\mathscr{U}_{k_{1}+k_{2}-1}\end{multline}
where the final map on the right is \( \gamma  \). 

These compositions satisfy the following conditions, for all \( a\in \mathscr{U}_{n} \),
\( b\in \mathscr{U}_{m} \), and \( c\in \mathscr{U}_{t} \):
\begin{description}
\item [Associativity]\( (a\circ _{i}b)\circ _{j}c=a\circ _{i+j-1}(b\circ _{j}c) \)
\item [Commutativity]\( a\circ _{i+m-1}(b\circ _{j}c)=(-1)^{mn}b\circ _{j}(a\circ _{i}c) \)
\item [Equivariance]\( a\circ _{\sigma (i)}(\sigma \cdot b)=\tunderi {n}{i}(\sigma )\cdot (a\circ _{i}b) \)
\end{description}
\end{defn}
\begin{rem}
In \cite{Smith:1994}, I originally \emph{defined} operads (or formal coalgebras)
in terms of these compositions. It turned out that I'd recapitulated the historical
sequence of events: operads were originally defined this way and called \emph{composition
algebras.} I am indebted to Jim Stasheff\index{Stasheff@\textsc{J.  Stasheff}}
for pointing this out to me. 

Given this definition of operad, we recover the \( \gamma  \) map in \ref{def:operad}
by setting: 
\[
\gamma (u_{i_{1}}\otimes \cdots \otimes u_{i_{k}})=u_{i_{1}}\circ _{1}\cdots \circ _{k-1}u_{i_{k}}\circ _{k}u_{k}\]
 (where the implied parentheses associate to the right). It is left to the reader
to verify that the two definitions are equivalent (the commutativity condition,
here, is a special case of the equivariance condition).
\end{rem}
Morphisms of operads are defined in the obvious way:

\begin{defn}
\label{def:operadmorphism} Given two operads \( \mathscr{U} \) and \( \mathscr{V} \),
a \emph{morphism} \index{morphism of operads}\index{operad!morphism}
\[
f:\mathscr{U}\rightarrow \mathscr{V}\]
 is a sequence of chain-maps 
\[
f_{i}:\mathscr{U}_{i}\rightarrow \mathscr{V}_{i}\]
 commuting with all the diagrams in \ref{def:operad} or (equivalently) preserving
the composition operations in \ref{def:operadcomps}.
\end{defn}
Now we give some examples:

\begin{defn}
\label{def:mathfrakS0}The operad \( \mathfrak{S}_{0} \) is defined \index{operad! $\mathfrak{S}_0$}\index{s@$\mathfrak{S}_0$ operad}via
\end{defn}
\begin{enumerate}
\item Its \( n^{\mathrm{th}} \) component is \( \zs{{n}} \) --- a chain-complex
concentrated in dimension \( 0 \). 
\item The composition operations are defined by
\[
\gamma (\sigma _{i_{1}}\otimes \cdots \otimes \sigma _{i_{k}}\otimes \sigma )=\sigma _{i_{\sigma (k)}}\oplus \cdots \oplus \sigma _{i_{\sigma (k)}}\circ \tmap {i_{1},\dots ,i_{k}}(\sigma )\]

\end{enumerate}
\begin{rem}
This was denoted \( \mathrf {M} \) in \cite{Kriz-May}.
\end{rem}
Verification that this satisfies the required identities is left to the reader
as an exercise.

\begin{defn}
\label{def:sfrakfirstmention}Let \( \mathfrak{S} \) denote the \( E_{\infty } \)-operad
with \index{operad!$\mathfrak{S}$}\index{definition!$\mathfrak{S}$ operad}components
\[
\rs{{n}}\]
--- the bar resolutions of \( \integers  \) over \( \zs{{n}} \) for all \( n>0 \).
This is a particularly important operad and \S~\ref{sec:sfrakdef} is devoted
to it.
\end{defn}
\begin{rem}
This is the result of applying the ``unreduced bar construction'' to the previous
example.
\end{rem}
\begin{defn}
\label{def:coassoc}\( \coassoc  \) is an operad \index{operad!Coassoc}\index{Coassoc operad}defined
to have one basis element \( \{b_{i}\} \) for all integers \( i\geq 0 \).
Here the rank of \( b_{i} \) is \( i \) and the degree is 0 and the these
elements satisfy the composition-law: \( b_{i}\circ _{\alpha }b_{j}=b_{i+j-1} \)
regardless of the value of \( \alpha  \), which can run from \( 1 \) to \( j \).
The differential of this operad is identically zero.
\end{defn}
Now we define two important operads associated to any \( \integers  \)-module.

\begin{defn}
\label{def:coend} Let \( C \) be a DGA-module with \index{operad!CoEnd}\index{definition!CoEnd}\index{co-endomorphism operad, definition}augmentation
\( \epsilon :C\rightarrow \mathbb{Z} \), and with the property that \( C_{0}=\mathbb{Z} \).
Then the \emph{Coendomorphism} operad, \( \coend (C) \), is defined to be the
operad with: 
\begin{enumerate}
\item component of \( \rank i=\homz (C,C^{i}) \), with the differential induced by
that of \( C \) and \( C^{i} \). The dimension of an element of \( \homz (C,C^{i}) \)
(for some \( i \)) is defined to be its degree as a map. 
\item The \( \mathbb{Z} \)-summand is generated by one element, \( e \), of rank
0. 
\end{enumerate}
Suppose \( D\subset C \) is a direct summand. We define, \( \coend (D,C) \),
the \emph{co-endomorphism operad of} \( C \) \emph{relative to} \( D \) to
be the sub-operad of \( \coend (C) \) generated by elements of \( \homz (C,C^{i}) \)
that send \( D\subset C \) to \( D^{i}\subset C^{i} \).

\end{defn}
\begin{rem}
Both operads are unitary --- their common identity element is the identity map
\( \mathrm{id}\in \homz (C,C) \). One motivation for operads is that they model
the iterated coproducts that occur in \( \coend (*) \). We will use operads
as an algebraic framework for defining other constructs that have topological
applications.
\end{rem}
In like fashion, we define the \textit{endomorphism} operad:

\begin{defn}
\label{def:end} If \( C \) is a DGA-module, \index{operad!End}\index{endomorphism operad}the
\emph{endomorphism operad}, \( \zend (C) \) is defined to have components 
\[
\homz (C^{n},C)\]
 and compositions that coincide with endomorphism compositions. If \( a_{n}:S_{n}\rightarrow \mathrm{Aut}(C) \)
defines a group-action for all \( n>0 \), we define:
\begin{enumerate}
\item \emph{the right crossed endomorphism operad\index{right crossed endomorphism operad}\index{operad!right crossed endomorphism}}
\( \zend _{\{a_{n}\}}(C) \) to have components \( \homz (C^{n},C) \), as before,
and a \( \zs{{n}} \)-action that simultaneously permutes the factors of \( C^{n} \)
in \( \homz (C^{n},C) \) and acts on the right hand argument via \( \{a_{n}\} \). 
\item The \emph{left crossed endomorphism operad} \( {}_{\{a_{n}\}}\zend (C) \) to
have components \( \homz (C^{n},C) \) and \( \zs{{n}} \)-action that simultaneously
permutes the factors of \( C^{n} \) in \( \homz (C^{n},C) \) and acts on each
factor via \( \{a_{n}\} \). 
\end{enumerate}
\end{defn}
There is an interesting duality between the endomorphism and co-endomorphism
operads:

\begin{thm}
\label{lem:operadduality}(Duality Theorem) Let \( C \) be a chain complex
and let \( A \) be a ring equipped with an \( S_{n} \)-action \( a_{n}:S_{n}\rightarrow \mathrm{Aut}(A) \)
for all \( n>0 \) with the property that the following diagram and commutes
for all \( n>0 \) and \( \sigma \in S_{n} \). \[\xymatrix{{A^n}\ar[r]^{\sigma} \ar[d]_{\mu^{n-1}}&{A^n}\ar[d]^{\mu^{n-1}}\\ {A}\ar[r]_{a_n(\sigma)}&{A}}\]
where \( A^{n} \) is the \( n \)-fold \( \integers  \)-tensor product of
copies of \( A \) and \( \mu ^{n-1} \) is the iterated product of \( A \).
There exists a morphism of operads (called the cup product morphism):
\[
\coend (C)\rightarrow \zend _{\{a_{n}'\}}(\homz (C,A))\]
\index{morphism!cup product}\index{duality!operad}\index{cup product morphism}where
\( a_{n}'=\homz (1,a_{n}) \).
\end{thm}
\begin{rem}
This morphism is injective if the ring \( A \) has no \( 0 \)-divisors.

If \( A \) is a DGA-algebra, \( \homz (C,A) \) must be the hypercohomology
chain complex.

If \( A \) is commutative, the \( S_{n} \)-action on \( A \) is trivial and
\( \zend _{\{a_{n}\}}(\homz (C,A)) \) becomes the ordinary (non-twisted) endomorphism
operad.
\end{rem}
\begin{proof}
Given any morphism \( z:C\rightarrow C^{n} \) we construct a morphism \( z':\homz (C,A)^{n}\rightarrow \homz (C,A) \)
as follows:

For all \( n>1 \), define a map
\[
H_{n:}\underbrace{{\homz (C,A)\otimes \cdots \otimes \homz (C,A)}}_{n\, \mathrm{factors}}\rightarrow \homz (C^{n},A^{n})\]
in the obvious way: send \( f_{1}\otimes \cdots \otimes f_{n} \) with \( f_{i}\in \homz (C,A) \)
to the morphism \( C^{n}\rightarrow A^{n} \) that sends \( c_{1}\otimes \cdots \otimes c_{n} \)
to \( f_{1}(c_{1})\otimes \cdots \otimes f_{n}(c_{n})\in A^{n} \). Now define
\( z':\homz (C,A)^{n}\rightarrow \homz (C,A) \) to be the composite \[\xymatrix{{\homz(C,A)^n} \ar[d]|-{H_n}\\ {\homz(C^n,A^n)} \ar[d]|-{\homz(1_{C^n},\mu^{n-1})}\\ {\homz(C^n,A))}\ar[d]|-{\homz(z,1_A)}\\{\homz(C,A)}}\] where
\( \mu ^{n-1}:A^{n}\rightarrow A \) is the product operation of \( A \). This
map is easily see to preserve the \( S_{n} \)-action.

We must now verify that this preserves compositions. Let \( z_{1}:C\rightarrow C^{n} \)
and \( z_{2}:C\rightarrow C^{m} \) and consider the \( i^{\mathrm{th}} \)
composition \( z_{1}\circ _{i}z_{2}:C\rightarrow C^{n+m-1} \): \[\xymatrix{{\homz(C,A)^{n+m-1}} \ar[d]|-{H_{n+m-1}}\\ {\homz(C^{n+m-1},A^{n+m-1})} \ar[d]|-{\homz(1_{C^{n+m-1}},\mu^{n+m-2})}\\ {\homz(C^{n+m-1},A))}\ar[d]|-{\homz(z_1\circ_i z_2,1_A)}\\{\homz(C,A)}}\]
We begin by decomposing some maps (using the associativity of \( \mu  \)) to
get \[\xymatrix{{\homz(C,A)^{n+m-1}} \ar[d]|-{H_{n+m-1}}\\ {\homz(C^{n+m-1},A^{n+m-1})} \ar[d]|-{\homz(1_{C^{n+m-1}},1\otimes\cdots\otimes\mu^{m-1}\otimes\cdots\otimes 1)}\\ {\homz(C^{n+m-1},A^{n})}\ar[d]|-{\homz(1_{C^{n+m-1}},\mu^{n-1})}\\ {\homz(C^{n+m-1},A))}\ar[d]|-{\homz( 1\otimes\cdots\otimes z_2 \otimes\cdots 1, 1_A)}\\ {\homz(C^{n},A))} \ar[d]|-{\homz(z_1,1_A)}\\ {\homz(C,A)}}\] Now
we permute two maps to get \[\xymatrix{{\homz(C,A)^{n+m-1}} \ar[d]|-{H_{n+m-1}}\\ {\homz(C^{n+m-1},A^{n+m-1})} \ar[d]|-{\homz(1_{C^{n+m-1}},1\otimes\cdots\otimes\mu^{m-1}\otimes\cdots\otimes 1)}\\ {\homz(C^{n+m-1},A^{n})}\ar[d]|-{\homz( 1\otimes\cdots\otimes z_2 \otimes\cdots 1, 1_A)}\\ {\homz(C^{n},A^n)}\ar[d]|-{\homz(1_{C^{n+m-1}},\mu^{n-1})}\\ {\homz(C^{n},A))} \ar[d]|-{\homz(z_1,1_A)}\\ {\homz(C,A)}}\] The
conclusion follows by permuting the second and third map from the top with \( H_{n+m-1} \).
\end{proof}
Now we consider a special class of operads that play a crucial role in the sequel:

\begin{defn}
\label{def:einfinity} An operad \( \mathscr{U}=\{\mathscr{U}_{n}\} \) will
be \index{E@$E_{\infty}$ operads}\index{operad!$E_{\infty}$}called an \( E_{\infty } \)-operad
if \( \mathscr{U}_{n} \) is a \( \zs{n} \)-free resolution of \( \mathbb{Z} \)
for all \( n>0 \). 
\end{defn}
\begin{rem}
In \cite{Kriz-May}, Kriz \index{Kriz@\textsc{I. Kriz}}\index{May@\textsc{P. J. May}}and
May define an \( A_{\infty } \)-operad as a non-\( \Sigma  \) analogue of
an \( E_{\infty } \)-operad. We do not do this here since we will be especially
concerned with a \emph{particular} \( A_{\infty } \)-operad, defined in \S~\ref{sec:ainfinity}. 
\end{rem}
\begin{prop}
\label{propmorphhomiso}Every morphism of \( E_{\infty } \)-operads 
\[
m:\mathfrak{R}_{1}\rightarrow \mathfrak{R}_{2}\]
induces homology isomorphisms in all dimensions and for all components.
\end{prop}
\begin{proof}
This is due to \ref{def:operad} and the requirement that every morphism preserve
the unit map.
\end{proof}
\begin{cor}
\label{cor:operadpullback}Every diagram of \( E_{\infty } \)-operads \[\xymatrix{{\mathfrak{R}_A}&{\mathfrak{R}_1}\ar[l]_{f}\\{\mathfrak{R}_2}\ar[u]^{g}&}\] can
be completed to a pullback \[\xymatrix{{\mathfrak{R}_A}&{\mathfrak{R}_1}\ar[l]_{f}\\{\mathfrak{R}_2}\ar[u]^{g}&{\mathfrak{R}_D}\ar@{.>}[u]_{u}\ar@{.>}[l]^{v}}\] 
\end{cor}
\begin{proof}
Here \( \mathfrak{R}_{D} \) is the standard graded module-pullback $\mathrm{ker}(f-g):\mathfrak{R}_B\oplus\mathfrak{R}_C \rightarrow \mathfrak{R}_A$.
It is clear that this kernel is an operad, and a glance at the exact sequence
in homology induced by $f-g$ (and the fact that $g$ and $f$ induce \emph{homology
isomorphisms} in all dimensions and for all degrees --- see \ref{propmorphhomiso})
shows that it is also \( E_{\infty } \). 
\end{proof}
Given these definitions, we can define:

\begin{defn}
\label{def:operadcomodule} Let \( \mathscr{U} \) be an operad and let \( C \)
be a DG-module equipped with a morphism (of operads) 
\[
f:\mathscr{U}\rightarrow \coend (C)\]
 Then \( C \) is called a \emph{coalgebra} over \( \mathscr{U} \) \index{coalgebra over an operad}\index{operad!coalgebra over an}with
structure map \( f \). If \( C \) is equipped with a morphism of operads 
\[
f:\mathscr{U}\rightarrow \zend (C)\]
 then \( C \) is called a \emph{algebra} over \( \mathscr{U} \) \index{algebra over an operad}\index{operad!algebra over an}with
structure map \( f \). 
\end{defn}
\begin{rem}
A coalgebra, \( C \), over an operad, \( \mathscr{U} \), is a sequence of
maps 
\[
f_{n}:\mathscr{U}\otimes C\rightarrow C^{n}\]
 for all \( n>0 \), where \( f_{n} \) is \( \zs{n} \)-equivariant. These
maps are related in the sense that they fit into commutative diagrams:  \[\xymatrix@d@R+20pt{{\mathscr{U}_{n}\otimes\mathscr{U}_{m}\otimes  C}\ar[r]^-{\circ_{i}}&{\mathscr{U}_{n+m-1}\otimes 
C}\ar[r]^-{f_{n+m-1}}&{C^{n+m-1}}\\ 
{\mathscr{U}_{n}\otimes\mathscr{U}_{m}\otimes C}\ar[r]_-{1\otimes 
f_{m}}\ar@{=}[u]&{\mathscr{U}_{n}\otimes 
C^{m}}\ar[r]_-{V_{i-1}}&{C^{i-1}\otimes\mathscr{U}_{n}\otimes C\otimes 
C^{m-i}}\ar[u]_-{1\otimes\dots\otimes f_{n}\otimes\dots\otimes1}}
\]  for all \( n,m\geq 1 \) and \( 1\leq i\leq m \). Here \( V:\mathscr{U}_{n}\otimes C^{m}\rightarrow C^{i-1}\otimes \mathscr{U}_{n}\otimes C\otimes C^{m-i} \)
is the map that shuffles the factor \( \mathscr{U}_{n} \) to the right of \( i-1 \)
factors of \( C \). In other words: The abstract composition-operations in
\( \mathscr{U} \) exactly correspond to compositions of maps in \( \{\homz (C,C^{n})\} \).
We exploit this behavior in applications of coalgebras over operads, using an
explicit knowledge of the algebraic structure of \( \mathscr{U} \).

Note that a morphism \( f:\mathfrak{U}\rightarrow \coend (D,C) \), where \( \coend (D,C) \)
is defined in \ref{def:coend} corresponds to structures of \emph{coalgebras
with sub-coalgebras}.
\end{rem}
In very simple cases, one can explicitly describe the maps defining a coalgebra
over an operad:

\begin{defn}
\label{def:unitinterval}Define the coalgebra,\index{unit interval} \( I \)
--- \emph{the unit interval} --- over \( \mathfrak{S} \) (see \ref{def:sfrakfirstmention}
and \S~\ref{sec:sfrakdef}) via:
\end{defn}
\begin{enumerate}
\item Its \( \integers  \)-generators are \( \{p_{0},p_{1}\} \) in dimension \( 0 \)
and \( q \) in dimension \( 1 \), and its adjoint structure map is \( r_{n}:\rs{{n}}\otimes I\rightarrow I^{n} \),
for all \( n>1 \).
\item The coproduct is given by \( r_{2}([\, ]\otimes p_{i})=p_{i}\otimes p_{i} \),
\( i=0,1 \), and \( r_{2}([\, ]\otimes q)=p_{0}\otimes q+q\otimes p_{1} \).
\item The higher coproducts are given by \( r_{2}([(1,2)]\otimes q=q\otimes q \),
\( r_{2}([(1,2)]\otimes p_{i})=0 \), \( i=0,1 \) and \( r_{2}(a\otimes I)=0 \),
where \( a\in \rs{{2}} \) has dimension \( >1 \).
\end{enumerate}
\begin{rem}
This, coupled with the operad-identities in \( \mathfrak{S} \) suffice to define
the coalgebra structure of \( I \) in all dimensions and for all degrees.
\end{rem}
\begin{prop}
Coassociative coalgebras are precisely the coalgebras over \( \coassoc  \). 
\end{prop}
\begin{rem}
There are some subtleties to this definition, however. It is valid if we regard
\( \coassoc  \) as a non-\( \Sigma  \) operad. If we regard it as an operad
with \emph{trivial} symmetric group action, then we have defined coassociative,
\emph{cocommutative} coalgebras. 
\end{rem}
We have two complementary results:

\begin{prop}
\label{prop:everycomodule} Every chain complex is trivially a coalgebra over
its own coendomorphism operad. 
\end{prop}
Although the following result is elementary, it must be stated. It implies that
we can form \emph{quotients} of coalgebras over operads:

\begin{prop}
\label{prop:comodulequotient}Suppose \( f:D\rightarrow C \) is a split injection
of chain-complexes, and let \( \coend (D,C) \) be the relative co-endomorphism
operad (see \ref{def:coend}). Then \( f \) induces a morphism
\[
\coend (D,C)\rightarrow \coend (C/D)\]
 of operads.
\end{prop}
\begin{rem}
It is interesting that the corresponding statement for \( \zend (C) \) (or
a relative analogue thereof) is not true: it is easier to take quotients of
coalgebras than of algebras. The operad \( \zend (C) \) defines algebras that
have a \emph{multiplication} rather than a \emph{comultiplication}. And taking
quotients of rings by subrings requires the subring to be an \emph{ideal}.
\end{rem}
\begin{proof}
Suppose \( h\in \coend (D,C) \) is some element. Then we have a commutative
diagram \[\xymatrix{{C}\ar[r]^{h}&{C^n}\\{D}\ar[r]_{h|D}\ar@{ >-}[u]&{D^n}\ar@{ >-}[u]}\] which
we can expand to a diagram \[\xymatrix{&{{(C/D)}^n}\\{C/D}\ar[r]^{h^*}&{C^n/D^n}\ar[u]\\{C}\ar[r]^{h}\ar@{->>}[u]&{C^n}\ar@{->>}[u]\\{D}\ar[r]_{h|D}\ar@{ >-}[u]&{D^n}\ar@{ >-}[u]}\] which
gives us an element of \( \coend (C/D) \). To conclude that this defines a
morphism of operads, we must show that it respects compositions (as in \ref{def:operadcomps}),
or that the kernel of the map \( \coend (D,C)\rightarrow \coend (C/D) \) (of
graded \( \integers  \)-modules is an \emph{ideal} with respect to all compositions. 

Suppose \( h\in \coend (D,C) \) maps to \( 0 \) in \( \coend (C/D) \). Then
the image of \( h \) (in \( C^{n} \)) must have at least one factor that is
in \( D \) --- say the \( i^{\mathrm{th}} \). If we form a composition with
\( h \) on another factor --- i.e., take \( \circ _{j} \) with \( j\ne i \)
--- then this factor of \( D \) is unchanged and appears in the composition.
Consequently, the composition has a factor of \( D \) and maps to \( 0 \)
in \( \coend (C/D) \). On the other hand, if we formed the composition \( \circ _{i} \)
the result will have \emph{many} factors of \( D \) because the co-endomorphisms
in \( \coend (D,C) \) \emph{preserve} \( D \).
\end{proof}
We can define tensor products of operads:

\begin{defn}
Let \( \mathscr{U}_{1} \) and \( \mathscr{U}_{2} \) be operads. Then \( \mathscr{U}_{1}\otimes \mathscr{U}_{2} \)
is defined to have: 
\begin{enumerate}
\item component of \( \rank i=(\mathscr{U}_{1})_{i}\otimes (\mathscr{U}_{2})_{i} \),
where \( (\mathscr{U}_{1})_{i} \) and \( (\mathscr{U}_{2})_{i} \) are, respectively,
the components of \( \rank i \) of \( \mathscr{U}_{1} \) and \( \mathscr{U}_{2} \); 
\item composition operations defined via \( (a\otimes b)\circ _{i}(c\otimes d)=(-1)^{\dim (b)\dim (c)}(a\circ _{i}c)\otimes (b\circ _{i}d) \),
for \( a,c\in \mathscr{U}_{1} \) and \( b,d\in \mathscr{U}_{2} \). 
\end{enumerate}
\end{defn}
We conclude this section with 

\begin{defn}
\label{def:operadcoproduct}An operad, \( \mathfrak{R} \), will be called an
\emph{operad-coalgebra} if there exists a co-associative morphism of operads
\[
\Delta :\mathfrak{R}\rightarrow \mathfrak{R}\otimes \mathfrak{R}\]

\end{defn}
\begin{rem}
Operad-coalgebras are important in certain homotopy-theoretic contexts --- for
instance in the study of the bar and cobar constructions. The \S~\ref{sec:sfrakdef}
defines a particularly important operad of this type.
\end{rem}

\subsection{Free algebras and coalgebras over an operad}

The canonical examples of operads are \( \zend (C) \) and \( \coend (C) \)
for some module \( C \). In this section, we will show that these examples
are fundamental in the sense that a large number of operads can be represented
as suboperads of these examples.

We begin by considering an alternate way to define algebras and coalgebras over
an operad. 

\begin{defn}
\label{def:monad}Let \( \mathrf {G} \) be a category. A \emph{monad} \index{monad!definition}\index{definition!monad}in
\( \mathrf {G} \) is a functor \( C:\mathrf {G}\rightarrow \mathrf {G} \)
together with natural transformations \( \mu :CC\rightarrow C \) and \( \eta :\mathrm{Id}\rightarrow C \)
such that the following diagrams commute: \[\xymatrix{{C}\ar[r]^{\eta C}\ar@{=}[rd] & {CC}\ar[d]^{\mu} & {C}\ar[l]_{C\eta}\ar@{=}[ld] \\ {} & {C} & {}} \mathrm{~and~} \xymatrix{{CCC}\ar[r]^{C \mu}\ar[d]_{\mu C} & {CC} \ar[d]^{\mu}\\ {CC}\ar[r]_{\mu} & {C}}\]

A \( C \)-\emph{algebra} \index{monad!algebra over}\index{definition!algebra over a monad}\index{algebra over a monad}is
an object \( A\in \mathrf {G} \) together with a map \( \xi :CA\rightarrow A \)
such that the following diagrams commute \[\xymatrix{{A} \ar[r]^{\eta}\ar@{=}[rd]& {CA} \ar[d]^{\xi}\\ {} & {A} }\mathrm{~and~}\xymatrix{{CCA}\ar[r]^{C\xi}\ar[d]_{\mu} & {CA}\ar[d]^{\xi} \\ {CA} \ar[r]_{\xi}& {A}}\]
\end{defn}
\begin{rem}
Monads are the category-theoretic exegesis of the notion of ``system with composition-operation.''

Note that, for any object \( A\in \mathrf {G} \), \( CA \) is an algebra over
\( C \). This prompts the definition:
\end{rem}
\begin{defn}
\label{def:freealgebraovermonad}Let \( A\in \mathrf {G} \) be any object and
let \( C:\mathrf {G}\rightarrow \mathrf {G} \) be a monad. Then \index{free algebra over a monad}the
\emph{free \( C \)-algebra generated by \( A \)} is defined to be \( CA \).
\end{defn}
The use of the term \emph{free} is justified by the following result, which
implies that any morphism from \( A \) to a \( C \)-algebra \( B \) has a
\emph{unique} extension to \( CA \):

\begin{prop}
\label{prop:freealgebraovermonadproof}Let \( CA \) be the free algebra over
the monad \( C \) and let \( C[\mathrf {G}] \) be the category of \( C \)-algebras.
Then there exists an isomorphism
\[
\mathrm{hom}(\eta ,1):C[\mathrf {G}](CA,B)\rightarrow \mathrf {G}(A,B)\]

\end{prop}
\begin{proof}
Its inverse sends a morphism \( f:A\rightarrow B \) to \( \xi \circ Cf:CA\rightarrow B \).
\end{proof}
In like fashion, we define \emph{comonads} and \emph{coalgebras} over them:

\begin{defn}
\label{def:comonad}Let \( \mathrf {G} \) be \index{comonad}\index{definition!comonad}a
category. A \emph{comonad} in \( \mathrf {G} \) is a functor \( D:\mathrf {G}\rightarrow \mathrf {G} \)
together with natural transformations \( \delta :D\rightarrow DD \) and \( \epsilon :D\rightarrow \mathrm{Id} \)
such that the following diagrams commute: \[\xymatrix{{D} & {DD}\ar[l]_{\epsilon D}\ar[r]^{D\epsilon} & {D} \\ {} & {D}\ar@{=}[lu]\ar[u]^{\delta} \ar@{=}[ru]& {}} \mathrm{~and~} \xymatrix{{DDD} & {DD}\ar[l]_{D \delta} \\ {DD}\ar[u]^{\delta D} & {D}\ar[l]^{\delta}\ar[u]_{\delta }}\]

A \( D \)-\emph{coalgebra} is \index{coalgebra over a comonad}an object \( A\in \mathrf {G} \)
together with a map \( \theta :A\rightarrow DA \) such that the following diagrams
commute \[\xymatrix{{A} & {DA}\ar[l]_{\epsilon} \\ {} & {A} \ar[u]_{\theta}\ar@{=}[lu]}\mathrm{~and~}\xymatrix{{DDA} & {DA}\ar[l]_-{D\theta} \\ {DA}\ar[u]^-{\delta}& {A} \ar[l]^-{\theta}\ar[u]_-{\theta}}\]

Given any object \( X\in \mathrf {G} \), and a comonad over \( \mathrf {G} \),
call \( DX \) the \emph{free \( D \)-coalgebra generated by} \( X \). 
\end{defn}
As before, we have

\begin{prop}
\label{prop:freecoalgebraproof}Let \( DX \) be the free coalgebra over the
monad \( D \) and let \( D[\mathrf {G}] \) be the category of \( D \)-coalgebras.
Then there exists an isomorphism
\[
\mathrm{hom}(\epsilon ,1):\mathrf {G}(A,B)\rightarrow D[\mathrf {G}](DA,B)\]

\end{prop}
\begin{proof}
Its inverse sends a morphism \( f:DA\rightarrow B \) to \( f\circ \theta :A\rightarrow B \).
\end{proof}
The relation between operads and monads (and comonads) is expressed by:

\begin{defn}
\label{def:operadmonadcomonad}Let \( \mathscr{U} \) be an operad. If \( \mathrf {G} \)
is the category of DGA-modules over \( \integers  \), then the\index{monad!defined by an operad}\index{operad!corresponding monad}
\emph{monad defined by} \( \mathscr{U} \) is given by
\[
CX=\bigoplus _{n\geq 0}\mathscr{U}_{n}\otimes _{\zs{{n}}}X^{n}\]
 and the \emph{comonad defined by\index{comonad defined by an operad}\index{operad!corresponding comonad}}
\( \mathscr{U} \) is given by
\[
DX=\bigoplus _{n\geq 0}\homzs{n}(\mathscr{U}_{n},X^{n})\]
where \( S_{n} \) acts on the \( n \)-fold tensor product \( X^{n} \) by
permutation of factors.
\end{defn}
\begin{rem}
It is well-known that the operad identities in \ref{def:operad} translate into
a proof that \( C \) and \( D \) are monads and comonads, respectively. See
\cite{Kriz-May}. It is also well-known that algebras over \( C \) are the
same as algebras over \( \mathscr{U} \) and coalgebras over \( D \) are the
same as coalgebras over \( \mathscr{U} \). For instance, a coalgebra, \( X \),
over \( D \) comes equipped with a map 
\[
X\rightarrow DX\]
 or
\[
X\rightarrow \bigoplus _{n\geq 0}\homzs{n}(\mathscr{U}_{n},X^{n})\]
which is equivalent to a sequence of \( S_{n} \)-equivariant maps
\[
\mathscr{U}_{n}\otimes X\rightarrow X^{n}\]
 defining the adjoint of the structure map of a coalgebra over \( \mathscr{U} \).
A corresponding statement holds for algebras over a \( C \).

By abuse of notation, we will identify free algebras over \( C \) with free
algebras over \( \mathscr{U} \) and free coalgebras over \( D \) with free
coalgebras over \( \mathscr{U} \).

Note that free coalgebras over comonads (or operads) have an infinite number
of chain-modules in negative dimensions --- and there is nontrivial homology
involved. This implies that they are \emph{never} m-coalgebras.
\end{rem}
We will show that many operads possess faithful representations in the sense
that they map injectively into \( \zend (C) \) and \( \coend (C') \) for suitable
\( C \) and \( C' \). To that end, we prove:

\begin{lem}
\label{lem:operadsimp}Let \( \mathfrak{R}=\{\mathfrak{R}_{i}\} \) be an operad
such that \( \mathfrak{R}_{i} \) is \( \zs{{i}} \)-free for all \( i \) and
let \( R_{i}=\mathfrak{R}_{i}\otimes _{\zs{{i}}}\integers  \) for all \( i \).
Then there exists a surjection \[\xymatrix{{\mathfrak{R}_i}\ar@{->>}[r]^{f_i}&{R_i}}\]
and a map \( \gamma ':R_{j_{1}}\otimes \cdots \otimes R_{j_{k}}\otimes \mathfrak{R}_{k}\rightarrow R_{j} \)
that makes the following diagram commute: \[\xymatrix@C+20pt{{\mathfrak{R}_{j_{1}}\otimes\cdots\otimes\mathfrak{R}_{j_{k}}\otimes\mathfrak{R}_{k}}\ar[r]^-{\gamma}\ar[d]_{(f_{j_1}\otimes\cdots{f_{j_k}})\otimes 1}&{\mathfrak{R}_{j}}\ar[d]^{f_{j}}\\
{R_{j_{1}}\otimes\cdots\otimes\ R_{j_{k}}\otimes\mathfrak{R}_{k}}\ar[r]_-{\gamma'}&{R_{j}}}\] where
\( j=\sum _{i=1}^{k}j_{i} \). In addition, the following diagram commutes:
\[\xymatrix@C+20pt{{R_{j_{1}}\otimes\cdots\otimes R_{j_{k}}\otimes\mathfrak{R}_{k}}\ar[r]^-{\gamma'}\ar[d]_{\sigma^{-1}\otimes\sigma}&{R_{j}}\ar[d]^{1}\\
{R_{j_{\sigma(1)}}\otimes\cdots\otimes R_{j_{\sigma(k)}}\otimes\mathfrak{R}_{k}}\ar[r]_-{\gamma'}&{R_{j}}}\]
\end{lem}
\begin{proof}
The map \( f_{i} \) fits into an exact sequence
\[
IS_{j_{i}}\otimes R_{j_{i}}\rightarrow \mathfrak{R}_{j_{i}}\rightarrow R_{j_{i}}\rightarrow 0\]
where \( IS_{j_{i}}=\ker (\zs{{j_{i}}}\rightarrow \integers ) \). The first
equivariance condition in the definition of an operad (see \ref{def:operad})
implies that \( IS_{j_{i}}\cdot \mathfrak{R}_{j_{i}} \) maps into \( \iota _{j_{i}}(ID_{j_{i}})\cdot \mathfrak{R}_{j} \),
where \( \iota _{j_{i}}:\zs{{j_{i}}}\rightarrow \zs{{j}} \) is induced by the
composite
\[
S_{j_{i}}\rightarrow S_{j_{1}}\oplus \cdots \oplus S_{j_{k}}\rightarrow S_{j}\]
It follows that \( \iota _{j_{i}}(ID_{j_{i}})\subseteq IS_{j} \) so that the
image of the kernel of \( f_{j_{i}} \) lies in the kernel of \( f_{j} \) and
the \( \gamma ' \) is well-defined. The remaining statement follows from the
first equivariance condition in \ref{def:operad}.
\end{proof}
\begin{prop}
\label{prop:everyoperadalgebra}Let \( \mathfrak{R}=\{\mathfrak{R}_{i}\} \)
be an operad such that \( \mathfrak{R}_{i} \) is \( \zs{{i}} \)-free for all
\( i \) and let \( C \) be the free algebra over \( \mathfrak{R} \) on one
generator in dimension \( 0 \). Then the structure map
\[
f:\mathfrak{R}\rightarrow \zend (C)\]
is injective.
\end{prop}
\begin{proof}
By definition~\ref{def:freealgebraovermonad} our free algebra is given by 
\[
C=\bigoplus _{i=1}^{\infty }\mathfrak{R}_{i}\otimes _{\zs{{i}}}\integers \]
Lemma \ref{lem:operadsimp} implies the existence of chain maps
\[
R_{d_{1}}\otimes \cdots \otimes R_{d_{n}}\otimes \mathfrak{R}_{n}\rightarrow \mathfrak{R}_{n}\otimes R_{d_{1}}\otimes \cdots \otimes R_{d_{n}}\rightarrow R_{d_{1}+\cdots +d_{n}}\]
 for all \( n>0 \) and \( n \)-tuples \( \{d_{1},\dots ,d_{n}\} \) of positive
integers. The map on the left is just transposition of the factors. This induces
a morphism of DGA-modules
\[
f_{n}:\mathfrak{R}_{n}\rightarrow \homz (R_{d_{1}}\otimes \cdots \otimes R_{d_{n}},R_{d_{1}+\cdots +d_{n}})\]
 The associativity conditions imply that compositions are preserved by \( f \),
and the second statement of Lemma \ref{lem:operadsimp} implies that \( f \)
is \( \zs{{n}} \)-linear, where \( S_{n} \) acts on the right by permuting
factors of \( R_{d_{1}}\otimes \cdots \otimes R_{d_{n}} \). The injectivity
of \( f \) follows from the fact that one of the targets of any \( \mathfrak{R}_{n} \)
is \( \homz (R_{n}\otimes R_{1}\otimes \cdots \otimes R_{1},R_{n}) \), with
a \( \zs{{n}} \)-basis of elements of \( \mathfrak{R}_{n} \) mapping to their
image in \( R_{n} \).
\end{proof}
\begin{prop}
\label{prop:everyoperad} Let \( \mathfrak{R}=\{\mathfrak{R}_{n}\} \) be an
operad such that \( \mathfrak{R}_{n} \) is \( \zs{{n}} \)-free for all \( n \),
and let \( C \) be the free coalgebra over \( \mathfrak{R} \) on one generator
in dimension \( 0 \). Then the structure map 
\[
\mathfrak{c}:\mathfrak{R}\rightarrow \coend (C,C^{n})\]
is injective
\end{prop}
\begin{rem}
This result and \ref{prop:everyoperadalgebra} show that endomorphism and coendomorphism
operads are, in some sense, \emph{canonical} examples: \emph{every} operad is
a sub-operad of one of them. We will use this result several times in the sequel
to show that ``universal equations'' in coalgebras over operads imply equations
in the operads themselves. See, for instance, \ref{prop:equivprop}.

Since this coalgebra is concentrated in non-positive dimensions, it is not an
m-coalgebra. It will suffice for our applications, however.
\end{rem}
\begin{proof}
In this case, we set
\[
Z=\bigoplus _{n=1}^{\infty }\homz (\mathfrak{R}_{n}\otimes _{\zs{{n}}}\integers ,\integers )=\bigoplus _{n=1}^{\infty }\homzs{{n}}(\mathfrak{R}_{n},\integers )\]

As before, we make crucial use of lemma \ref{lem:operadsimp}. We must proceed
more delicately than in \ref{prop:everyoperadalgebra}, though. We will use
the approach to operads in \ref{def:operadcomps} and define composition operations.
A composition
\[
R_{k}\circ _{i}\mathfrak{R}_{k}\rightarrow R_{i+k-1}\]
(where \( R_{i}=\mathfrak{R}_{i}\otimes _{\zs{{i}}}\integers  \)) implies the
existence of a chain-map
\[
\circ _{i}':\mathfrak{R}_{k}\rightarrow \homz (R_{i},R_{i+k-1})\]
--- this is just the transposition and adjoint. Dualizing gives rise to a map
\[
\circ _{i}'':\mathfrak{R}_{k}\rightarrow \homz (\homz (R_{i+k-1},\integers ),\homz (\mathfrak{R}_{j},\integers ))\]
and it is important to realize that this represents \( k \) \emph{distinct}
maps (corresponding to the different composition operations we may perform).
Hence the \emph{real} target of \( \homz (R_{i+k-1},\integers ) \) is 
\[
(\integers \oplus \homz (R_{j},\integers ))^{k}\subseteq (\homz (R_{1},\integers )\oplus \homz (R_{j},\integers ))^{k}\subset C^{k}\]
--- in fact it is
\[
\bigoplus _{\alpha =1}^{k}\homz (R_{j},\integers )=\bigoplus _{\alpha =1}^{k}\underbrace{{\integers \otimes \cdots \otimes \homz (R_{j},\integers )\otimes \cdots \otimes \integers }}_{R_{j}\, \mathrm{in}\, \mathrm{position}\, \alpha }\]
Injectivity of the structure-map follows from its sending \( \mathfrak{R}_{k} \)
to (among other things
\begin{eqnarray*}
\lefteqn {\homz (\homz (R_{k},\integers ),}\hspace {.5in} &  & \\
 &  & \homz (R_{k},\integers )\otimes \homz (R_{1},\integers )\otimes \cdots \otimes \homz (R_{1},\integers )))
\end{eqnarray*}

\end{proof}

\subsection{The operad $\mathfrak{S}$\label{sec:sfrakdef}}

We will study the symmetric construct, mentioned in \ref{def:sfrakfirstmention}.
It is the result of applying the unreduced bar construction to the (somewhat
trivial) operad \( \mathfrak{S}_{0} \) defined in \ref{def:mathfrakS0}. Since
its components are unreduced bar constructions (or bar resolutions of \( \integers  \)
over \( \zs{{n}} \), we can exploit the Cartan Theory of constructions in \cite{Cartan}\index{Cartan theory of constructions}\index{Cartan@\textsc{Henri Cartan}}.
The basic idea of the Cartan Theory of Constructions is the following:

\begin{lem}
\label{lem:cartanconstruction}Let \( M_{i} \), \( i=1,2 \) be DGA-modules,
where:
\begin{enumerate}
\item \( M_{1}=A_{1}\otimes N_{1} \), where \( N_{1} \) is \( \integers  \)-free
and \( A_{1} \) is a DGA-algebra (so \( M_{1} \), merely regarded as a DGA-algebra,
is free on a basis equal to a \( \integers  \)-basis of \( N_{1} \))
\item \( M_{2} \) is a left DGA-module over a DGA-algebra \( A_{2} \), possessing 

\begin{enumerate}
\item a sub-DG-module, \( N_{2}\subset M_{2} \), such that \( \partial _{M_{2}}|N_{2} \)
is injective, 
\item a contracting chain-homotopy \( \varphi :M_{2}\rightarrow M_{2} \) whose image
lies in \( N_{2}\subset M_{2} \).
\end{enumerate}
\end{enumerate}
Suppose we are given a chain-map \( f_{0}:M_{1}\rightarrow M_{2} \) in dimension
\( 0 \) with \( f_{0}(N_{1})\subseteq N_{2} \) and want to extend it to a
chain-map from \( M_{1} \) to \( M_{2} \), subject to the conditions:

\begin{itemize}
\item \( f(N_{1})\subseteq N_{2} \)
\item \( f(a\otimes n)=g(a)\cdot f(n) \), where \( g:A_{1}\rightarrow A_{2} \) is
some morphism of DG-modules such that \( a\otimes n\mapsto g(a)\cdot f(n) \)
is a chain-map.
\end{itemize}
Then the extension \( f:M_{1}\rightarrow M_{2} \) exists and is unique.

\end{lem}
\begin{rem}
In applications of this result, the morphism \( g \) will often be a morphism
of DGA-algebras, but this is not necessary.

The \emph{existence} of \( f \) immediately follows from basic homological
algebra; the interesting aspect of it is \emph{}its \emph{uniqueness} (not merely
uniqueness up to a chain-homotopy). We will use it repeatedly to prove associativity
conditions by showing that two apparently different maps satisfying the hypotheses
must be \emph{identical}.

This elementary but powerful result was due to Henri Cartan and formed the cornerstone
of his theory of Eilenberg-MacLane spaces.
\end{rem}
\begin{proof}
This result is due to Cartan. The uniqueness of \( f \) follows by induction
and: 
\begin{enumerate}
\item \( f \) is determined by its values on \( N_{1} \) 
\item the image of the contracting chain-homotopy, \( \varphi  \), lies in \( N_{2}\subset M_{2} \).
\item the boundary map of \( M_{2} \) is \emph{injective} on \( N_{2} \) (which
implies that there is a \emph{unique} lift of \( f \) into the next higher
dimension).
\end{enumerate}
\end{proof}
This will turn out to be an \( E_{\infty } \)-operad such that any simplicial
set's chain-complex is a natural m-coalgebra over it (see~\ref{def:m-coalgebra}).
This will form the basis of our topological results.

We define the operad structure of \( \mathfrak{S} \) in terms of compositions,
as in \ref{def:operadcomps} using \ref{lem:cartanconstruction}. Recall that
the \( \{\rs{n}\} \) are simplicial sets for all \( n>0 \) with face and degeneracy
operations given by: 

\begin{equation}
\label{eq:faceops}
\mathrm{F}_{i}[a_{1}|\dots |a_{m}]=\left\{ \begin{array}{ll}
{a_{1}[a_{2}|\dots |a_{m}]} & \mathrm{if}\, i=0\\
{[a_{1}|\dots |a_{i}\cdot a_{i+1}|\dots |a_{m}]} & \mathrm{if}\, 0<i<m\\
{[a_{1}|\dots |a_{m-1}]} & \mathrm{if}\, i=m
\end{array}\right. 
\end{equation}
 and, as usual, the differential is defined as an alternating sum of these face-operations.

It also has a well-known co-associative coalgebra structure given by \( \Delta _{R}(a)=\sum _{i=0}^{n}\lface ^{i}(a)\otimes \fface ^{n-i}(a) \),
where \( \lface  \) is the last face operator. With the face operations given
above, this amounts to \( \Delta _{R}:[a_{1}|\dots |a_{n}]=\sum _{i=0}^{n}[a_{1}|\dots |a_{i}]\otimes a_{1}\cdots a_{i}[a_{i+1}|\dots |a_{n}] \).
In addition, \label{pag:asnref}

\begin{enumerate}
\item \( \rs{{n}}=\zs{{n}}\otimes A(S_{n},1) \), as \( \zs{{n}} \)-modules, where
\( A(S_{n},1) \), in dimension \( k \), is the submodule generated by elements
of the form \( 1[a_{1}|\dots |a_{k}] \).
\item There exists a contracting cochain, \( \varphi :\rs{{n}}\rightarrow \rs{{n}} \),
is given in dimension \( k \) by
\begin{eqnarray}
\varphi (a[a_{1}|\dots |a_{k}]) & = & 1[a|a_{1}|\dots |a_{k}]\in A(S_{1},1)_{k+1}\label{eq:contractcochain} \\
\varphi (1[a_{1}|\dots |a_{k}]) & = & 0\label{eq:contractcochain2} 
\end{eqnarray}

\end{enumerate}
Clearly, the twisted differential (the alternating sum of face-operators in
\ref{eq:faceops}) is \emph{injective} on \( 1\otimes A(S_{n},1)\subset \rs{{n}} \)
--- in fact, it defines an isomorphism \( A(S_{n},1)_{k}\cong (\rs{n})_{k-1} \). 

Now we define the compositions \( \rs{{n}}\circ _{i}\rs{{m}}\rightarrow \rs{{n+m-1}} \),
requiring them to satisfy the condition that
\begin{equation}
\label{eq:carcancondition}
(1\otimes A(S_{n},1))\circ _{i}(1\otimes A(S_{m},1))\subseteq 1\otimes A(S_{n+m-1},1)
\end{equation}

We set \( M_{1}=\rs{{n}}\otimes \rs{{m}} \), \( M_{2}=\rs{{n+m-1}} \), \( N_{2}=A(S_{n+m-1},1) \).
At this point, a slight problem arises. We would like to regard \( M_{1} \)
as a free \( \zs{{n}}\otimes \zs{m} \)-module, but it is equal to \( \zs{{n}}\otimes A(S_{n},1)\otimes \zs{{m}}\otimes A(S_{m},1) \).
We can simply rearrange factors and redefine the differential.

\begin{cor}
\label{cor:mathfrakSoperad}There exists a unique composition operation \( \rs{{n}}\circ _{i}\rs{{m}}\rightarrow \rs{{n+m-1}} \)
satisfying the equivariance conditions in \ref{def:operadcomps} and the condition
in equation \ref{eq:carcancondition}. In addition, these \( \{\circ _{i}\} \)
satisfy the associativity and commutativity conditions.
\end{cor}
\begin{rem}
This is a unitary operad. Its identity element is the generator \( [\, ]\in \{\rn 1\}_{0} \).
Note that an element \( [g_{1}|\dots |g_{k}] \) is not determined in \( \mathfrak{S} \)
by the symmetric-group elements \( g_{1},\dots ,g_{k} \) --- one must also
specify which summand, \( \rs{{n}} \), the element resides in. We therefore
distinguish between the \( S_{n} \) and their isomorphic images in \( S_{N} \),
where \( N>n \). 
\end{rem}
\begin{proof}
Existence and uniqueness follow from \ref{lem:cartanconstruction}, setting
\( A_{1}=\zs{{n}}\otimes \zs{{m}} \) and \( N_{1}=A(S_{n},1)\otimes A(S_{m},1) \)
and the fact that, in dimension \( 0 \) \( N_{1}=N_{2}=[\, ]\otimes [\, ]\cong \integers  \),
so any chain maps between them that commute with augmentation must be equal.
In dimension \( 0 \)
\[
[\, ]\circ _{i}[\, ]=[\, ]\]
for all \( i \). In higher dimensions, set 
\[
a\circ _{i}b=\varphi (\partial a\circ _{i}b+(-1)^{k}a\circ _{i}\partial b)\]
for \( a\in 1\otimes A(S_{n},1)_{k} \), \( b\in 1\otimes A(S_{m,}1) \) and
for all \( i \). The boundaries of \( a \) and \( b \) can be arbitrary elements
of \( \rs{{n}} \) and \( \rs{m} \), respectively (one dimension lower), so
we must compute the composition of the boundaries via the formulas:

\[
a\circ _{\sigma (i)}\sigma b=\tunderi {n}{i}(\sigma )\cdot (a\circ _{i}b)\]
and
\[
\tau a\circ _{i}b=\underbrace{{1\oplus \cdots \oplus \tau \oplus \cdots \oplus 1}}_{i^{\mathrm{th}}\, \mathrm{position}}(a\circ _{i}b)\]
 given by the equivariance conditions in \ref{def:operadcomps}.

Associativity and commutativity follow from the \emph{same} construction carried
out with\( M_{1}=\rs{{n}}\otimes \rs{{m}}\otimes \rs{{t}} \) and \( M_{2}=\rs{{n+m+t-2}} \).
In this case, we compare the maps
\[
f_{1}=(*\circ _{i}*)\circ _{j}*:\rs{{n}}\otimes \rs{{m}}\otimes \rs{{t}}\rightarrow \rs{{n+m+t-2}}\]
and
\[
f_{2}=*\circ _{i+j-1}(*\circ _{j}*):\rs{{n}}\otimes \rs{{m}}\otimes \rs{{t}}\rightarrow \rs{{n+m+t-2}}\]
Both compositions agree in dimension \( 0 \) and 
\[
f_{i}(1\otimes A(S_{n}1)\otimes 1\otimes A(S_{m},1)\otimes 1\otimes A(S_{t},1))\subseteq 1\otimes A(S_{n+m-2},1)\]
Furthermore, the equations for extending the \( f_{i} \) from \( 1\otimes A(S_{n}1)\otimes 1\otimes A(S_{m},1)\otimes 1\otimes A(S_{t},1) \)
to all of \( M_{1} \) must agree, since they are given by the composition operations
on the operad \( \mathfrak{S}_{0} \) --- see \ref{def:mathfrakS0}. Consequently,
\( f_{1}=f_{2} \) in \emph{all} dimensions.

A similar argument (with only one factor, obviously) implies that \( \circledast  \)
preserves coproducts. See \cite{Smith:1994} for an explicit formula for the
composition operations of \( \mathfrak{S} \). The compositions given in that
paper \emph{must} agree with the ones proved to exist in \ref{cor:mathfrakSoperad},
because they satisfy the defining conditions and because such compositions are
\emph{unique}.
\end{proof}
I originally computed these compositions these --- as chain-maps \( \rs{{n}}\otimes \rs{{m}}\rightarrow \rs{{n+m-1}} \)
satisfying the equivariance conditions. In fact, I first defined \( A(S_{n},1)\circ _{i}[\, ] \)
and \( [\, ]\circ _{i}A(S_{m},1) \) and had to find a way to combine them together,
so I derived the twisted shuffle operation (using the Cartan Theory of Constructions)
and \emph{multiplied} them. I was reassured of the correctness of doing so by
the fact that the twisted shuffle product is the unique one (satisfying natural
conditions). Then I proved associativity and commutativity by strenuous brute-force
computations (unaware that I could have used the Cartan Theory \emph{throughout}).

\begin{thm}
\label{th:soperad}Equipped with the composition operations defined above, \( \mathfrak{S} \)
is an \( E_{\infty } \)-operad-coalgebra with coproduct
\[
\Delta :\mathfrak{S}\rightarrow \mathfrak{S}\otimes \mathfrak{S}\]

\end{thm}
\begin{proof}
We have already proved that \( \mathfrak{S} \) is an operad --- see \ref{cor:mathfrakSoperad}.
To see that the coproduct is an operad morphism, observe that the composition
operations of \( \mathfrak{S} \) are \emph{coalgebra} morphisms. This follows
from the Cartan Theory of Constructions once again, observing that the \( \Delta  \)
can be \emph{defined} by this theory --- as the \emph{unique} map
\[
\Delta :\rs{{n}}\rightarrow \rs{{n}}\otimes \rs{{n}}\]
with the properties:
\begin{enumerate}
\item \( \Delta (1\otimes [\, ])=1\otimes [\, ]\otimes 1\otimes [\, ] \)
\item \( \Delta (1\otimes A(S_{n},1))\subset 1\otimes A(S_{n},1)\otimes \rs{{n}}=N_{2} \)
\item \( \Delta  \) is a morphism of \( \zs{{n}} \)-modules (where \( \rs{{n}}\otimes \rs{{n}} \)
is equipped with the diagonal \( S_{n} \)-action).
\item \( \Delta (1\otimes A)=(\varphi \otimes 1+\epsilon \otimes \varphi )(\Delta \partial A) \),
where \( A\in A(S_{n},1) \), \( \varphi  \) is the contracting cochain in
\ref{eq:contractcochain2}, and \( \epsilon :\rs{{n}}\rightarrow \integers  \)
is the augmentation (so \( \varphi \otimes 1+\epsilon \otimes \varphi  \) is
a contracting cochain of \( \rs{{n}}\otimes \rs{{n}} \)).
\end{enumerate}
Having done this, we would only need observe that composing \( \Delta  \) with
the \( \{\circ _{i}\} \) preserves these defining conditions.

\end{proof}
We conclude this section with an application of the operad \index{C@$\cf{*}$-functor!definition}\index{definition!$\cf{*}$-functor}\( \mathfrak{S} \):

\begin{thm}
\label{th:scrCdef}There exists a functor 
\[
\mathscr{C}:\mathrm{SS}\rightarrow \mathrm{m}-\mathrm{Coalgebras}\, \, \mathrm{over}\, \, \mathfrak{S}\]
such that the underlying chain-complex of \( \mathscr{C}(X) \) is the simplicial
chain-complex, where \( \mathrm{SS} \) is the category of pointed, simply-connected,
2-reduced simplicial sets .
\end{thm}
\begin{proof}
For all \( k \), let \( s^{k} \) denote the standard \( k \)-simplex, whose
vertices are \( \{[0],\dots ,[k]\} \) and whose \( j \)-faces are \( \{[i_{0},\dots ,i_{j}]\} \),
with \( i_{1}<\cdots <i_{j} \), \( j\le k \). We define \( \mathscr{C} \)
to be a free functor on models that are simplices and use the Cartan Theory
of Constructions to define maps:
\[
f_{n}:\rs{{n}}\otimes C(s^{k})\rightarrow {C(s^{k})}^{n}\]
where \( s^{k} \) is a \( k \)-simplex with chain complex \( C(s^{k}) \),
and \( S_{n} \) acts on \( {C(s^{k})}^{n} \) by permuting factors. These maps
must make the following diagram commute: \textit{\emph{$$\xymatrix@d{{\rs{n}\otimes\rs{m}\otimes{C}}\ar[r]^-{\circ_i}&{\rs{n+m-1}\otimes{C}}\ar[r]^-{f_{n+m-1}}&{C^{n+m-1}}\\{\rs{n}\otimes\rs{m}\otimes{C}}\ar[r]_-{1\otimes{f_m}}\ar@{=}[u]&{\rs{n}\otimes{C^m}}\ar[r]_-{V_{i-1}}&{C^{i-1}\otimes\rs{n}\otimes{C}\otimes{C^{m-i}}}\ar[u]_-{1\otimes\cdots\otimes{f_n}\otimes\cdots\otimes1}}\label{dia:coalgebradefined}$$}}
\textit{\emph{and for all \( n,m\geq 1 \) and \( 1\leq i\leq m \), where \( C=C(s^{k}) \)}}
\textit{\emph{and \( V:\rs{{n}}\otimes C^{m}\rightarrow C^{i-1}\otimes \rs{{n}}\otimes C\otimes C^{m-i} \)
is the map that shuffles the factor \( \rs{{n}} \) to the right of \( i-1 \)
factors of \( C \). We define a contracting cochain on \( C \) via
\[
\mathfrak{s}([i_{1},\dots ,i_{t}]=\left\{ \begin{array}{cc}
(-1)^{t}[i_{1},\dots ,i_{t},k] & \mathrm{if}\, i_{t}\ne k\\
0 & \mathrm{if}\, i_{t}=k
\end{array}\right. \]
where \( [i_{1},\dots ,i_{t}] \), \( i_{1}<i_{2}<\dots <k \), denotes a \( t-1 \)-dimensional
face of \( s^{k}=[0,\dots ,k] \). We can define a corresponding contracting
homotopy on \( C^{n} \) via \( \Phi =\mathfrak{s}\otimes 1\otimes \cdots \otimes 1+\epsilon \otimes \mathfrak{s}\otimes \cdots \otimes 1+\cdots +\epsilon \otimes \cdots \otimes \epsilon \otimes \mathfrak{s} \),
where 
\[
\epsilon ([i_{1},\dots ,i_{t}])=\left\{ \begin{array}{cc}
[k] & \mathrm{if}\, [i_{1},\dots ,i_{t}]=[k]\\
0 & \mathrm{otherwise}
\end{array}\right. \]
}}Above dimension \( 0 \), \( \Phi  \) is effectively equal to \( \mathfrak{s}\otimes 1\otimes \cdots \otimes 1 \)\textit{\emph{.
Now set \( M_{2}=C^{n} \) and \( N_{2}=\im (\Phi ) \). In dimension \( 0 \),
we define \( f_{n} \) for all \( n \) via:
\[
f_{n}(A\otimes [0])=\left\{ \begin{array}{ll}
[0]\otimes \cdots \otimes [0] & \mathrm{if}\, A=[\, ]\\
0 & \mathrm{if}\, \dim A>0
\end{array}\right. \]
This clearly makes \( \mathscr{C}(s^{0}) \) a coalgebra over \( \mathfrak{S} \).}}

\textit{\emph{Suppose that the \( f_{n} \) are defined below dimension \( k \).
Then, by acyclic models, \( \mathscr{C}(\partial s^{k}) \) is well-defined
and satisfies the conclusions of this theorem. We define \( f_{n}(a[a_{1}|\dots |a_{j}]\otimes [0,\dots ,k]) \)
by induction on \( j \), requiring that the following invariant condition be
satisfied:}}

\[
f_{n}(A(S_{n},1)\otimes s^{k})\subseteq [i_{1},\dots ,k]\otimes \mathrm{other}\, \mathrm{terms}\]
--- in other words, the leftmost factor must be in \( \im \mathfrak{s} \).
Now we simply set
\begin{eqnarray*}
f_{n}(A\otimes s^{k}) & = & \Phi \circ f_{n}(\partial A\otimes s^{k})\\
 & + & (-1)^{\dim A}\Phi \circ f_{n}(A\otimes \partial s^{k})
\end{eqnarray*}
where \( A\in A(S_{n},1) \).

The upper line follows from induction on the dimension of \( A \) and the lower
line follows from acyclic models and the induction on \( k \). The Cartan Theory
of Constructions implies that this map is uniquely determined by the invariant
condition and the contracting homotopy \( \Phi  \). It follows that diagram
\ref{dia:coalgebradefined} must commute since:
\begin{enumerate}
\item any composite of \( f_{n} \)-maps will continue to satisfy the invariant condition.
\item \( \circ _{i}(1\otimes A(S_{n},1)\otimes \cdots \otimes A(S_{m},1))\subseteq 1\otimes A(S_{n+m-1},1) \)
so that composing an \( f_{n} \)-map with \( \circ _{i} \) results in a map
that still satisfies the invariant condition.
\item the diagram commutes in lower dimensions (by induction on \( k \) and acyclic
models)
\end{enumerate}
\end{proof}
We conclude this section with a definition that will be important in the sequel:

\begin{defn}
\label{def:cartanmcoalgebra}Let \( C \) be a DGA complex that is an acyclic
m-coalgebra over \( \mathfrak{S} \) satisfying the conditions:
\begin{enumerate}
\item there exists a self-annihilating contracting homotopy \( \varphi :C\rightarrow C \)
\item the adjoint of the m-structure of \( C \), \( f_{n}:\mathfrak{S}\otimes C\rightarrow C^{n} \),
has the property that\emph{
\begin{equation}
\label{eq:c2invarcond}
f_{n}(A(S_{n},1)\otimes \varphi (C))\subseteq \Phi (C)
\end{equation}
}where: 

\begin{enumerate}
\item \( A(S_{n},1)\subset \rs{{n}} \) is as defined on page~\pageref{pag:asnref};
\item \textit{\emph{\( \Phi =\varphi \otimes 1\otimes \cdots \otimes 1+\epsilon \otimes \varphi \otimes \cdots \otimes 1+\cdots +\epsilon \otimes \cdots \otimes \epsilon \otimes \varphi  \)
is the contracting homotopy of \( C^{n} \) induced by \( \varphi  \), where
\( \epsilon :C\rightarrow \integers  \) is the augmentation. }}
\end{enumerate}
\end{enumerate}
Then \( C \) will be said to be \emph{Cartan, with contracting homotopy} \( \varphi  \).

\end{defn}
\begin{rem}
For instance, the m-coalgebra of any \emph{simplex} is clearly Cartan, by \ref{th:scrCdef}.
\end{rem}

\section{m-coalgebras and homotopy\label{sec:mcoalgehomotopy}}

\subsection{A category of fractions}

In this section, we will define four categories that will be important in the
sequel: 

\begin{itemize}
\item \( \mathfrak{M}_{0} \) 
\item \( \mathfrak{M}=\mathfrak{M}_{0}[S^{-1}] \)
\item \( \mathfrak{L}_{0} \) and
\item \( \mathfrak{L}=\mathfrak{L}_{0}[T^{-1}] \)
\end{itemize}
The \emph{objects} of these categories will be called \emph{m-coalgebras} ---
they are DG-coalgebras over \( E_{\infty } \)-operads, and the morphisms will
be those of DG-coalgebras, extended slightly to allow a change of underlying
operad.

The categories \( \mathfrak{M} \) and \( \mathfrak{L} \) are categories of
fractions of \( \mathfrak{M}_{0} \) and \( \mathfrak{L}_{0} \), respectively,
in the sense of Gabriel and Zisman in \index{Gabriel@\textsc{P. Gabriel}}\index{Zisman@\textsc{M. Zisman}}\cite{Gabriel-Zisman:1967}.
The classes of inverted morphisms, \( \{S\} \) and \( \{T\} \), are known
as \emph{elementary equivalences}, defined in \ref{def:emorph}, and \( \mathfrak{M} \)
and \( \mathfrak{L} \) are interesting because they contain sub-categories
equivalent to the homotopy category of pointed, simply-connected CW complexes.

Gabriel and Zisman presented a theory of categories of fractions and used it
to study various homotopy categories. Although categories of fractions exist
and are well-defined in great generality, their structure is also usually very
complex.

Gabriel and Zisman showed that, under certain (rigid) conditions, a category
and subcategory \emph{admits a calculus of left or right fractions}. When this
happens, the category of fractions has a particularly simple structure. 

Unfortunately, as \index{Quillen@\textsc{D. Quillen}}Quillen observed in \cite{Quillen:1967},
interesting topological categories rarely admit a calculus of left or right
fractions. Quillen replaced these \emph{algebraic} conditions by topologically
inspired ones (related to fibrations and cofibrations) and defined \emph{model
categories}. Model categories have properties that facilitate the study of homotopy
types without explicitly using a calculus of fractions.

It turns out that our category \( \mathfrak{M}_{0} \) and \( \mathfrak{L}_{0} \)
are \emph{not} model categories (except up to equivalence) so that we cannot
use Quillen's theory. They \emph{almost} admit a calculus of left fractions,
but fail to meet one of Gabriel and Zisman's conditions. On the other hand,
they have other algebraic properties that are almost as good as having a calculus
of fractions, so that we get a tractable homotopy theory .

\begin{defn}
\label{def:m-coalgebra} The category, \( \mathfrak{M}_{0} \) is defined \index{M0@$\mathfrak{M}_0$!definition}as
follows:
\begin{enumerate}
\item Its objects are \emph{m-coalgebras}, \index{definition!m-coalgebras}\index{m-coalgebra!definition}where
an m-coalgebra, \( C \), is an equivalence class of DG-coalgebras over \( E_{\infty } \)-operads
with\footnote{%
By abuse of notation, we use the same term for the underlying chain-complex
as for the m-coalgebra.
} \( C_{0}=\integers  \) and \( C_{1}=0 \). Two such DG-coalgebras, \( C_{1} \)
and \( C_{2} \), with structure morphisms \begin{equation}\xymatrix{{\mathfrak{R}_1}\ar[r]^-{m_1}&{\coend(C_1)}\\{\mathfrak{R}_2}\ar[r]_-{m_2}&{\coend(C_2)}}\label{eq:mcoalgdef}\end{equation}
define the \emph{same} m-coalgebra if \( C_{1}=C_{2} \) (as DG-modules), and
there exists an operad morphism \( f:\mathfrak{R}_{1}\rightarrow \mathfrak{R}_{2} \)
that makes \[\xymatrix{{\mathfrak{R}_1}\ar[r]^-{m_1}\ar[d]_-{f}&{\coend(C_1)}\ar@{=}[d]\\{\mathfrak{R}_2}\ar[r]_-{m_2}&{\coend(C_2)}}\] commute.
\item Morphisms \index{m-coalgebra!morphism}\index{M0@$\mathfrak{M}_0$!morphism}\index{morphism!m-coalgebra}induced
by those of DG-coalgebras over \( E_{\infty } \)-operads that preserve homology
in dimensions \( 0 \) and \( 1 \). In other words, a morphism between two
m-coalgebras, \( C_{1} \) and \( C_{2} \) over the operad \( \mathfrak{R} \),
with representative structure maps as in \ref{eq:mcoalgdef} consists of a morphism
of DG-coalgebras over \( \mathfrak{R} \): \( C_{1}\rightarrow C_{2} \). We
consider two morphisms to be the same if their underlying \emph{chain-maps}
are the same. 
\end{enumerate}
\end{defn}
\begin{rem}
m-coalgebras were first defined in \cite{Smith:1994}. The present definition
corresponds to the term \emph{weakly coherent} m-coalgebras in that paper. 

The category \( \mathfrak{M}_{0} \)  \emph{fails} to be Abelian --- cokernels
do not always exist (their underlying chain-complexes may have \( \integers  \)-torsion).

Note that \( C=\integers  \), concentrated in dimension \( 0 \) is an m-coalgebra
over any \( E_{\infty } \)-operad. This is the initial and final object in
\( \mathfrak{M}_{0} \) since it maps to any m-coalgebra and any m-coalgebra
maps to it.
\end{rem}
\begin{defn}
\label{def:L0category}Let \( \mathfrak{L}_{0} \) be the full subcategory of
\( \mathfrak{M}_{0} \) of m-coalgebras over the operad \( \mathfrak{S} \).
Define \( \mathfrak{F}_{0}\subset \mathfrak{L}_{0} \) to be the subcategory
of finitely generated m-coalgebras over \( \mathfrak{S} \).
\end{defn}
\begin{prop}
\label{prop:mcoalgprod} Let \( C_{1} \) and \( C_{2} \) be chain-complexes.
Then there exist natural transformations of functors 
\[
\mathfrak{E}_{n}:\homz (C_{1},C_{1}^{n})\otimes \homz (C_{2},C_{2}^{n})\rightarrow \homz (C_{1}\otimes C_{2},(C_{1}\otimes C_{2})^{n})\]
 for all \( n \). They induce operad morphisms
\[
\mathfrak{E}:\coend (C_{1})\otimes \coend (C_{2})\rightarrow \coend (C_{1}\otimes C_{2})\]

\end{prop}
\begin{rem}
If \( u\in \homz (C_{1},C_{1}^{n}) \), \( v\in \homz (C_{2},C_{2}^{n}) \),
then \( \mathfrak{E}_{n}(u\otimes v) \) sends \( c_{1}\otimes c_{2}\in C_{1}\otimes C_{2} \)
to \( V_{n}((u\otimes v)(c_{1}\otimes c_{2}) \), where \( V_{n}:C_{1}^{n}\otimes C_{2}^{n}\rightarrow (C_{1}\otimes C_{2})^{n} \)
is the map that shuffles the factors together. 
\end{rem}
We can now define tensor products of m-coalgebras:

\begin{defn}
\label{def:einftycomoduleprod} Let \( \mathrf {C}=\mathfrak{M}_{0},\mathfrak{L}_{0}\, \mathrm{or}\, \mathfrak{F}_{0} \)
and let \( C_{1} \) and \( C_{2} \) be objects of \( \mathrf {C} \) with
structure maps
\[
f_{i}:\mathfrak{R}_{i}\rightarrow \coend (C_{i})\]
 \( i=1,2 \). Their \emph{tensor product} is a \index{tensor product of m-coalgebras}\index{m-coalgebra!tensor product}DG-coalgebra
over \( \mathfrak{R}_{1}\otimes \mathfrak{R}_{2} \) with underlying chain complex
\( C_{1}\otimes C_{2} \) and structure map 
\[
\mathfrak{E}\circ f_{1}\otimes f_{2}:\mathfrak{R}_{1}\otimes \mathfrak{R}_{2}\rightarrow \coend (C_{1}\otimes C_{2})\]
 It is straightforward to verify that equivalent DG-coalgebras give rise to
equivalent tensor products so that this operation is well-defined for objects
of \( \mathfrak{M}_{0} \).

Given two objects \( C_{1} \) and \( C_{2} \) of \( \mathfrak{L}_{0} \) or
\( \mathfrak{F}_{0} \), form the tensor product as above. The result is an
m-coalgebra over the \( E_{\infty } \)-operad \( \mathfrak{S}\otimes \mathfrak{S} \).
Now pull this m-structure back over the diagonal map (see \ref{th:soperad})
\[
\mathfrak{S}\rightarrow \mathfrak{S}\otimes \mathfrak{S}\]
 to get an m-coalgebra over \( \mathfrak{S} \), hence an object of \( \mathfrak{L}_{0} \)
or \( \mathfrak{F}_{0} \), respectively.
\end{defn}
\begin{prop}
\label{prop:pushout}Let \( \mathrf {C}_{0}=\mathfrak{M}_{0} \), \( \mathfrak{L}_{0} \),
or \( \mathfrak{F}_{0} \) and let \[\xymatrix{{A}\ar[r]^{f}\ar[d]_{g}&{B}\\{C}&}\]
be a diagram in \( \mathrf {C}_{0} \) of m-coalgebras over an operad \( \mathfrak{R} \),
such that the quotient (of Abelian groups) \( (B\oplus C)/(f,-g)A \) is \( \integers  \)-torsion
free (for instance, if either \( f \) or \( g \) are split injections). Then
we may form the push out \[\xymatrix{{A}\ar[r]^{f}\ar[d]_{g}&{B}\ar@{.>}[d]^{u}\\{C}\ar@{.>}[r]_{v}&{D}}\]
in \( \mathrf {C}_{0} \).

This push out is natural with respect to morphisms of diagrams: \[\xymatrix{{.}\ar[r]\ar[d]\ar[ddr]&{.}\ar[ddr]&&\\{.}\ar[ddr]&&&\\&{.}\ar[r]\ar[d]&{.}&\\&{.}&&}\]

\end{prop}
\begin{proof}
We form the \emph{push out\index{push out}\index{M0@$\mathfrak{M}_0$!push out}}
of coalgebras over \( \mathfrak{R}_{D} \):
\[
D=\mathrm{coker}(f,-g):A\rightarrow B\oplus C\]
 which is possible because of proposition~\ref{prop:comodulequotient}. This
gives us the diagram \begin{equation}\xymatrix{{A}\ar[r]^{f}\ar[d]_{g}&{B}\ar[d]^{u}\\{C}\ar[r]_{v}&{D}}\label{dia:pushoutcommu}\end{equation} where
\( u \) and \( v \) are induced by the inclusions. They are clearly m-coalgebra
morphisms.
\end{proof}
\begin{defn}
\label{def:emorph} Let \( C \) and \( D \) be coalgebras over \( E_{\infty } \)-operads.
An injective morphism \( \iota :C\rightarrow D \) will be called an \emph{\index{elementary equivalence}\index{M0@$\mathfrak{M}_0$!elementary equivalence}\index{equivalence!elementary}elementary
equivalence} if its cokernel is a projective, acyclic chain complex.
\end{defn}
\begin{rem}
Note that, by \ref{prop:pushout}, the cokernel will possess an \emph{m-structure,}
and the quotient map \( D\rightarrow D/\iota (C) \) will be a morphism of m-coalgebras.
\end{rem}
\begin{defn}
\label{def:mcat} The \emph{homotopy category} \index{m-coalgebra!homotopy category}\index{category!$\mathfrak{M}$}of
m-coalgebras, denoted \( \mathfrak{M}=\mathfrak{M}_{0}[S^{-1}] \), is defined
to be the category of fractions of \( \mathfrak{M}_{0} \) by the subcategory,
\( S \), of elementary equivalences --- see  \cite{Gabriel-Zisman:1967}. Two
morphisms \( f_{1},f_{2}:A\rightarrow B \) in \( \mathfrak{M}_{0} \) will
be called \emph{homotopic} if they define the \emph{same} morphism in \( \mathfrak{M}=\mathfrak{M}_{0}[S^{-1}] \),
under the canonical functor \( \mathfrak{M}_{0}\rightarrow \mathfrak{M}_{0}[S^{-1}] \).
The \emph{saturation} of \( S \), denoted \( \hat{S} \) is defined to be all
morphisms in \( \mathfrak{M}_{0} \) that become invertible in \( \mathfrak{M} \).

In like fashion, we define the category \( \mathfrak{L}=\mathfrak{L}_{0}[T^{-1}] \)
and \( \mathfrak{F}=\mathfrak{F}_{0}[R^{-1}] \), where \( T \) and \( R \)
are the sub-categories of elementary equivalences in \( \mathfrak{L}_{0} \)
and \( \mathfrak{F}_{0} \), respectively.
\end{defn}
\begin{rem}
The canonical functors 
\begin{eqnarray*}
\mathfrak{M}_{0} & \rightarrow  & \mathfrak{M}_{0}[S^{-1}]\\
\mathfrak{L}_{0} & \rightarrow  & \mathfrak{L}_{0}[T^{-1}]\\
\mathfrak{F}_{0} & \rightarrow  & \mathfrak{F}_{0}[R^{-1}]
\end{eqnarray*}
 are, in general, not faithful. There exist an obvious functors
\begin{eqnarray*}
\mathfrak{f}:\mathfrak{L} & \rightarrow  & \mathfrak{M}\\
\mathfrak{g}:\mathfrak{F} & \rightarrow  & \mathfrak{L}
\end{eqnarray*}

It is easy to see that morphisms preserve coalgebra structures up to a chain-homotopy.
\end{rem}
\begin{defn}
\label{prop:hammockprop}Let \( A \) and \( B \) be m-coalgebras. A morphism
in \( \mathfrak{M}=\mathfrak{M}_{0}[S^{-1}] \) (respectively \( \mathfrak{L} \),
or \( \mathfrak{F} \)) is a diagram of the form \[\xymatrix{{...}\ar[r]^{m_1}&{.}&{.}\ar[l]_{s_1}\ar[r]^{m_2}&{.}&{.}\ar[l]_{s_2}\ar[r]^{}&{...}}\]
where the \( \{m_{i}\} \) are morphisms of m-coalgebras (respectively, m-coalgebras
over \( \mathfrak{S} \) or finitely generated m-coalgebras over \( \mathfrak{S} \))
and the \( \{s_{k}\} \) are elementary equivalences defined in \ref{def:emorph}.
A \emph{commutative hammock} \index{hammock!commutative}\index{commutative hammock}is
a diagram: \begin{equation}\xymatrix@R=10pt{{}&{X_1}\ar@{-}[r]\ar[dd]_{v_1}&{X_2}\ar@{-}[r]\ar[dd]^{v_2}&{\dots}\ar@{-}[r]&{X_k}\ar[rd]\ar[dd]^{v_k}&\\ {A}\ar[ur]\ar[dr]&&&&&{B}\\ {}&{Y_1}\ar@{-}[r]&{Y_2}\ar@{-}[r]&{\dots}\ar@{-}[r]&{Y_k}\ar[ur]&}\label{dia:hammock1}\end{equation}

where:
\begin{enumerate}
\item Horizontal maps in the same column (i.e., one above the other) go in the same
direction.
\item Vertical maps and horizontal maps to the right are morphisms of m-coalgebras.
\item Horizontal maps to the \emph{left} are elementary equivalences.
\end{enumerate}
\end{defn}
\begin{rem}
In any diagram like the above, we can clearly consolidate adjacent horizontal
maps that go in the same direction. Consequently, we may assume that adjacent
maps go in opposite directions.
\end{rem}
\begin{prop}
Let \( \mathrf {C}_{0}=\mathfrak{M}_{0} \), \( \mathfrak{L}_{0} \), or \( \mathfrak{F}_{0} \)
and let \( \mathrf {C} \) denote \( \mathfrak{M} \), \( \mathfrak{L} \),
or \( \mathfrak{F} \), respectively. The existence of a commutative hammock
in \( \mathrf {C}_{0} \) implies that its upper and lower rows are equal in
\( \mathrf {C} \).
\end{prop}
\begin{proof}
We will prove this for \( \mathfrak{M}_{0} \) --- the arguments for \( \mathfrak{L}_{0} \)
and \( \mathfrak{F}_{0} \) are identical. The claim follows by induction on
the number of columns. We show that the following equation (in \( \mathfrak{M}=\mathfrak{M}_{0}[S^{-1}] \))
is satisfied at column \( i \) of diagram \ref{dia:hammock1} for all \( i \):
\begin{equation}
\label{eq:diagraminvariant}
v_{i}\circ (\mathrm{upper}\, \mathrm{row})=(\mathrm{lower}\, \mathrm{row})
\end{equation}
 Every commutative subdiagram of the form:\[\xymatrix{{X_i}\ar[r]^{f_i}\ar[d]_{v_i}&{X_{i+1}}\ar[d]^{v_{i+1}}\\ {Y_i}\ar[r]_{g_i}&{Y_{i+1}} }\]results
in an equation 
\[
v_{i+1}f_{i}=g_{i+1}v_{i}\]
 and any commutative subdiagram of the form $$\xymatrix{{X_i}\ar[d]_{v_i}&{X_{i+1}}\ar[l]_{s_{i+1}}\ar[d]^{v_{i+1}}\\ {Y_i}&{Y_{i+1}} \ar[l]_{t_{i+1}} }$$
where \( s_{i+1} \) and \( t_{i+1} \) are elementary equivalences, results
in an equation of the form
\[
t_{i+1}^{-1}v_{i}=v_{i+1}s^{-1}_{i+1}\]
so that the validity of equation \ref{eq:diagraminvariant} in column \( i \)
implies its validity in column \( i+1 \). The conclusion follows.

The following converse was stated but not proved by Dwyer, Hirshhorn and Kan
in \cite{DHK:1997}\index{Dwyer@\textsc{W. Dwyer}}\index{Hirschhorn@\textsc{P. Hirschhorn}}\index{Kan@\textsc{D. Kan}}:
\end{proof}
\begin{prop}
Given two objects \( C \) and \( D \) in the category \( \mathfrak{M}_{0} \)
(respectively, \( \mathfrak{L}_{0} \) or \( \mathfrak{F}_{0} \)), the morphisms
\( \mathrm{Hom}_{\mathfrak{M}_{0}[S^{-1}]}(C,D) \) (respectively, \( \mathrm{Hom}_{\mathfrak{L}_{0}[T^{-1}]}(C,D) \)
or \( \mathrm{Hom}_{\mathfrak{F}_{0}[R^{-1}]}(C,D) \)) correspond to the components
of a 1-dimensional simplicial set whose 0-simplices are morphisms in \( \mathfrak{M}=\mathfrak{M}_{0}[S^{-1}] \)
(respectively, \( \mathfrak{L}=\mathfrak{L}_{0}[T^{-1}] \) or \( \mathfrak{F}=\mathfrak{F}_{0}[R^{-1}] \))
and whose 1-simplices are commutative hammocks.
\end{prop}
\begin{rem}
In other words two morphisms are equivalent \emph{if and only if} they can be
connected by a sequence of commutative hammocks.
\begin{proof}
Again, we will prove the result for \( \mathfrak{M}_{0} \) --- the arguments
for \( \mathfrak{L}_{0} \) and \( \mathfrak{F}_{0} \) are identical.

We use the description of localized categories in \cite{Gabriel-Zisman:1967}:
they define it via ``generators and relations''. The only relations that exist
in \( \mathfrak{M}=\mathfrak{M}_{0}[S^{-1}] \) are:
\end{proof}
\begin{enumerate}
\item The relations of \( \mathfrak{M}_{0} \) --- this is expressed by all formulas
\( a=bc \) that are true in \( \mathfrak{M}_{0} \), where \( a \), \( b \),
and \( c \) are morphisms such that \( b \) and \( c \) can be composed;
\item The relations that create a new set of formal morphisms, \( \{t_{i}\} \) and
define \( s_{i}t_{i}=1 \) and \( t_{i}s_{i}=1 \), where the \( s_{i} \) are
the elementary equivalences (we will use the notation \( s_{i}^{-1} \) for
\( t_{i} \) throughout the remainder of this argument).
\end{enumerate}
If \begin{equation} \xymatrix{{C}\ar[r]^{f_1}&{C_1}&{C_2}\ar[l]_{s_1}\ar[r]^{f_2}&{\dots}\ar[r]^{f_k}&D}\label{dia:fcomposite}\end{equation}
is a morphism in \( \mathfrak{M}=\mathfrak{M}_{0}[S^{-1}] \), then all morphisms
equivalent to it can be obtained by 

\begin{enumerate}
\item Factoring a morphism, \( a \), in \ref{dia:fcomposite} into a composite \( bc \)
whenever this is valid in \( \mathfrak{M}_{0} \), or the reverse --- composing
two adjacent morphisms that go in the same direction, or
\item Inserting pairs of morphisms \[\xymatrix{{.}\ar[r]^{s_i}&{.}&{.}\ar[l]_{s_i}}\] into
the diagram wherever it is possible to do so, or removing such pairs, or
\item Inserting pairs of morphisms \[\xymatrix{{.}&{.}\ar[l]_{s_i}\ar[r]^{s_i}&{.}}\] into
the diagram wherever it is possible to do so, or removing such pairs.
\end{enumerate}
It is straightforward to verify that all of these operations can be described
in terms of commutative hammocks. For instance, transformation 1 above corresponds
to a hammock that maps every term in \ref{dia:fcomposite} via the identity
map except for a central portion that looks like \[\xymatrix{{.}\ar[r]^{c}\ar[d]_{1}&{.}\ar[d]^{b}\ar[r]^{b}&{.}\\{.}\ar[r]_{a}&{.}\ar[r]_{1}&{.}}\] 

Transformation 2 corresponds to a hammock that maps every term of \ref{dia:fcomposite} via
the identity, except for a sub-diagram that looks like:\[\xymatrix{&{.}\ar[l]_{s}\ar[r]^{s}\ar[d]^{s}&\\&{.}\ar[l]^{1}\ar[r]_{1}&}\]
It is equally easy to see that transformation 3 can be realized by a commutative
hammock. The conclusion follows.

\end{rem}
Now we restrict our attention to the category \( \mathfrak{M}_{0} \) and its
category of fractions \( \mathfrak{M} \). The following result shows that \( \mathfrak{M} \)
satisfies \emph{one} of Gabriel and Zisman's conditions (in \cite{Gabriel-Zisman:1967})
for having a calculus of left fractions\index{calculus of left fractions}:

\begin{prop}
\label{prop:leftfraction}Let \( \mathrf {C}_{0}=\mathfrak{M}_{0} \), \( \mathfrak{L}_{0} \),
or \( \mathfrak{F}_{0} \) and let \( \mathrf {C} \) denote \( \mathfrak{M} \),
\( \mathfrak{L} \), or \( \mathfrak{F} \), respectively. If \[\xymatrix{{A}\ar[r]^{f}\ar[d]_{s}&{B}\\{C}&}\]
is a diagram in \( \mathrf {C} \), with \( s \) an elementary equivalence,
then there exists a diagram \begin{equation}\xymatrix{{A}\ar[r]^{f}\ar[d]_{s}&{B}\ar@{.>}[d]^{t}\\{C}\ar@{.>}[r]_{g}&{D}}\label{dia:leftfraction}\end{equation}
in \( \mathrf {C}_{0} \) with \( t \) an elementary equivalence. In particular,
for any m-coalgebra morphism \( f \) and elementary equivalence \( s, \) there
exist an m-coalgebra morphism \( g \) and elementary equivalence \( t \) such
that the relation \( fs^{-1}=t^{-1}g \) holds in \( \mathrf {C} \).
\end{prop}
\begin{proof}
Again, we only prove the statement for \( \mathfrak{M}_{0} \). This is an immediate
consequence of \ref{prop:pushout}. The statement that \( t \) is an elementary
equivalence follows from the fact that, as a chain-complex, \( D=B\oplus C/s(A) \),
where \( C/s(A) \) is \emph{acyclic} and \( \integers  \)-torsion free.

The equation in \( \mathfrak{M}_{0}[S^{-1}] \) follows from the commutativity
of \ref{dia:leftfraction}, which implies that \( tf=gs \). The invertibility
of \( s \) and \( t \) (in \( \mathfrak{M}_{0}[S^{-1}] \)) implies that we
can rewrite this equation as \( t^{-1}g=fs^{-1} \).
\end{proof}
\begin{cor}
\label{cor:reducehammock}Let \( \mathrf {C}_{0}=\mathfrak{M}_{0} \), \( \mathfrak{L}_{0} \),
or \( \mathfrak{F}_{0} \) and let \( \mathrf {C} \) denote \( \mathfrak{M} \),
\( \mathfrak{L} \), or \( \mathfrak{F} \), respectively. If \( m:A\rightarrow B \)
is a morphism in \( \mathrf {C} \), then \( m \) has a canonical representative
diagram of the form $\xymatrix{{A}\ar[r]^{f}&{C}&{b}\ar[l]_{s}}$, where \( f \)
is a morphism of m-coalgebras in \( \mathrf {C}_{0} \) and \( s \) is an elementary
equivalence in \( \mathrf {C}_{0} \). In addition, any commutative hammock
\begin{equation}\xymatrix@R=10pt{{}&{X_1}\ar@{-}[r]\ar[dd]_{v_1}&{X_2}\ar@{-}[r]\ar[dd]^{v_2}&{\dots}\ar@{-}[r]&{X_k}\ar[rd]\ar[dd]^{v_k}&\\ {A}\ar[ur]\ar[dr]&&&&&{B}\\ {}&{Y_1}\ar@{-}[r]&{Y_2}\ar@{-}[r]&{\dots}\ar@{-}[r]&{Y_k}\ar[ur]&}\end{equation} can
be simplified to a hammock of the form \begin{equation}\xymatrix@R=10pt{{}&{X}\ar[dd]_{v}&{}\\{A}\ar[ur]^{f_1}\ar[dr]_{f_2}&&{B}\ar[ul]_{s_1}\ar[dl]^{s_2}\\{}&{Y}&{}}\end{equation}
\end{cor}
\begin{rem}
We will call this canonical representative diagram, the \emph{canonical form}
\index{canonical form of a morphism in $\mathfrak{M}$}\index{morphism!canonical form}\index{M@$\mathfrak{M}$!morphism!canonical form}of
the morphism \( m \).
\end{rem}
\begin{proof}
The first statement follows immediately from \ref{prop:pushout}, which implies
that we can permute leftward elementary equivalences with the rightward maps
and compose --- thus simplifying the morphism. The second statement follows
from the \emph{naturality} of the push-outs used in the simplification: we can
simplify the upper and lower rows of a hammock in such a way that there is a
map of the simplifications.
\end{proof}
\begin{defn}
\label{def:ememtaryhomotopy}Let \( \mathrf {C}_{0}=\mathfrak{M}_{0} \), \( \mathfrak{L}_{0} \),
or \( \mathfrak{F}_{0} \) and let \( \mathrf {C} \) denote \( \mathfrak{M} \),
\( \mathfrak{L} \), or \( \mathfrak{F} \), respectively. If $\xymatrix{{A}\ar[r]^{f_i}&{X_i}&{B}\ar[l]_{s_i}}$
\( i=1,2 \) are two morphisms in \( \mathrf {C} \), a commutative diagram
in \( \mathrf {C}_{0} \) \[\xymatrix@R=10pt{{}&{X}\ar[dd]_{v}&{}\\{A}\ar[ur]^{f_1}\ar[dr]_{f_2}&&{B}\ar[ul]_{s_1}\ar[dl]^{s_2}\\{}&{Y}&{}}\label{dia:reducehammock}\]
will be called an \emph{\index{elementary homotopy}\index{M@$\mathfrak{M}$!elementary homotopy}elementary
homotopy}.
\end{defn}
\begin{cor}
\label{cor:elementaryhomotopy}Let \( \mathrf {C}_{0}=\mathfrak{M}_{0} \),
\( \mathfrak{L}_{0} \), or \( \mathfrak{F}_{0} \) and let \( \mathrf {C} \)
denote \( \mathfrak{M} \), \( \mathfrak{L} \), or \( \mathfrak{F} \), respectively.
Two morphisms of \( \mathrf {C} \) in canonical form $\xymatrix{{A}\ar[r]^{f_i}&{X_i}&{B}\ar[l]_{s_i}}$,
\( i=1,2 \) are equivalent in \( \mathrf {C} \) if and only if there exists
a sequence of elementary homotopies \begin{equation}\xymatrix@R=10pt{&{X_1}\ar@{-}[d]&\\&{Y_1}\ar@{-}[d]&\\{A}\ar[uur]|-{f_1}\ar[ur]|-{g_1}\ar[r]|-{g_k}\ar[dr]|-{g_n}\ar[ddr]|-{f_2}&{\vdots}\ar@{-}[d]&{B}\ar[uul]|-{s_1}\ar[ul]|-{t_1}\ar[l]|-{t_k}\ar[dl]|-{t_n}\ar[ddl]|-{s_2}\\&{Y_n}\ar@{-}[d]&\\&{X_2}&}\label{dia:elementaryhomotopyseq}\end{equation} between
them, in \( \mathrf {C}_{0} \). 
\end{cor}
\begin{rem}
Here all vertical maps are upward or downward morphisms of m-coalgebras, the
\( \{g_{i}\} \), \( i=1,\dots ,n \) are morphisms of m-coalgebras and the
\( \{t_{i}\} \) are elementary equivalences.
\end{rem}
\begin{proof}
This is an immediate consequence of \ref{prop:hammockprop} and \ref{cor:reducehammock}.
\end{proof}
It is interesting to see what happens if we can ``abelianize'' the category
\( \mathfrak{M}_{0} \):

\begin{thm}
\label{th:abelianization}Let \( \mathfrak{V} \) be an Abelian category and
let \( \mathfrak{f}:\mathfrak{M}_{0}\rightarrow \mathfrak{V} \) be a full functor.
In addition, let \( \mathfrak{V}'=\mathfrak{V}[\mathfrak{f}(S)^{-1}] \) be
the corresponding category of fractions. Then: 
\begin{enumerate}
\item any morphism in \( \mathfrak{V}' \) can be expressed as a diagram $$\xymatrix{{\mathfrak{f}(A)}\ar[r]^{\mathfrak{f}(f)}&{\mathfrak{f}(X)}&{\mathfrak{f}(B)}\ar[l]_{\mathfrak{f}(s)}}$$ 
\item two morphisms $\xymatrix{{\mathfrak{f}(A)}\ar[r]^{\mathfrak{f}(f_i)}&{\mathfrak{f}(X_i)}&{\mathfrak{f}(B)}\ar[l]_{\mathfrak{f}(s_i)}}$,
\( i=1,2 \) are equivalent (in \( \mathfrak{V}' \)), if and only if there
exists a commutative diagram in \( \mathfrak{V} \): \begin{equation}\xymatrix{&{\mathfrak{f}(X_1)}\ar[d]^{u}&\\{\mathfrak{f}(A)}\ar[ru]^{\mathfrak{f}(f_1)}\ar[rd]_{\mathfrak{f}(f_2)}&Z&{\mathfrak{f}(B)}\ar[lu]_{\mathfrak{f}(s_1)}\ar[ld]^{\mathfrak{f}(s_2)}\\&{\mathfrak{f}(X_2)}\ar[u]_{v}&}\label{dia:smallhomotopy}\end{equation} where
\( u \) and \( v \) are morphisms of \( \mathfrak{V} \) that are invertible
in \( \mathfrak{V}' \).
\end{enumerate}
In particular, if \( Z \) is the saturation of \( \mathfrak{f}(S) \), then
\( \mathfrak{V}'=\mathfrak{V}[Z^{-1}] \) admits a calculus of right fractions.

\end{thm}
\begin{rem}
For an example where this theorem applies, let \( \mathfrak{U} \) be the category
of \emph{rational} DG-coalgebras over \( E_{\infty } \)-operads, and otherwise
defined like \( \mathfrak{M}_{0} \). There is an obvious functor \( \mathfrak{M}_{0}\rightarrow \mathfrak{U} \)
and we define \( \mathfrak{V} \) to be the image of this functor. This is an
abelian category, and~\ref{th:abelianization} implies that its homotopy theory
has a particularly simple form.

Another example is the category, \( \mathfrak{N}_{0} \) of \emph{m-coalgebras
with torsion}. This is defined like \( \mathfrak{M}_{0} \) except that we allow
m-coalgebras to have \( \integers  \)-torsion. 

The definition of homotopy implied by a calculus of right fractions \index{calculus of right fractions}resembles
Quillen's definition of right homotopy in \cite{Quillen:1967}\index{Quillen@\textsc{D. Quillen}},
but is considerably simpler.
\end{rem}
\begin{proof}
The first statement follows from the fact that \( \mathfrak{f} \) is full.
The second follows from the fact that, in an abelian category, we can \emph{always}
form push-outs. We, consequently, take a diagram of elementary homotopies (like~\ref{dia:elementaryhomotopyseq})
and, whenever we encounter a subdiagram of the form \begin{equation}\xymatrix@R=10pt{&{\vdots}\ar@{-}[d]&\\&{Y_{i-1}}&\\{\mathfrak{f}(A)}\ar[uur]\ar[ur]|-{g_{i-1}}\ar[r]|-{g_i}\ar[dr]|-{g_{i+1}}\ar[ddr]&{Y_i}\ar[d]_-{w}\ar[u]^-{r}&{\mathfrak{f}(B)}\ar[uul]\ar[ul]|-{t_{i-1}}\ar[l]|-{t_k}\ar[dl]|-{t_{i+1}}\ar[ddl]\\&{Y_{i+1}}\ar@{-}[d]\\&{\vdots}&}\label{dia:pushout2}\end{equation} and
we form the push out of \[\xymatrix{{Y_i}\ar[r]^{r}\ar[d]_{w}&{Y_{i-1}}\\ {Y_{i+1}}&}\]

It is straightforward to check that this push-out can be inserted into diagram~\ref{dia:pushout2}.
After many cycles of consolidating morphisms going in the same direction and
repeating this push-out construction, we end up with a diagram that has the
required form. The final claim follows from page 12 of~\cite{Gabriel-Zisman:1967}.
\end{proof}

\subsection{Morphisms}

\begin{defn}
Let \( \mathrf {C}_{0}=\mathfrak{M}_{0} \), \( \mathfrak{L}_{0} \), or \( \mathfrak{F}_{0} \)
and let \( \mathrf {C} \) denote \( \mathfrak{M} \), \( \mathfrak{L} \),
or \( \mathfrak{F} \), respectively. A morphism \( f:C_{1}\rightarrow C_{2} \)
in \( \mathrf {C}_{0} \) is defined to be a \emph{cofibration} \index{cofibration}if
the push out exists in \( \mathrf {C}_{0} \) \[\xymatrix{{C_1}\ar[r]^{f}\ar[d]&{C_2}\ar@{.>}[d]\\{\integers}\ar@{.>}[r]&{D}}\]
in which case \( D \) is called its \emph{cofiber}. We extend this to \( \mathrf {C} \)
via the canonical functor \( \mathrf {C}_{0}\rightarrow \mathrf {C} \).
\end{defn}
Any map is homotopic to a cofibration, as we can see from the construction:

\begin{defn}
\label{def:mapcyl}Let \( \mathrf {C}_{0}=\mathfrak{M}_{0} \), \( \mathfrak{L}_{0} \),
or \( \mathfrak{F}_{0} \) and let \( \mathrf {C} \) denote \( \mathfrak{M} \),
\( \mathfrak{L} \), or \( \mathfrak{F} \), respectively. We define the\index{algebraic mapping cylinder}
\emph{algebraic mapping cylinder}, \( \mathscr{M}(f) \), of a morphism, \( f:C\rightarrow D \),
in \( \mathrf {C}_{0} \) to be the push out \begin{equation}\xymatrix{{C}\ar[r]^{f}\ar[d]_{1\otimes p_0}&{D}\ar@{.>}[d]^{\iota_D}\\{C \otimes I/ C_0 \otimes I }\ar@{.>}[r]_-{\eta_C}&{\mathscr{M}(f)}}\label{dia:mapcylpushout}\end{equation}

If \( C \) and \( D \) are m-coalgebras over \( \mathfrak{R} \), the underlying
operad of \( \mathscr{M}(f) \) is \( \mathfrak{R}\otimes \mathfrak{S} \).
When \( \mathrf {C}_{0}= \) \( \mathfrak{L}_{0} \) or \( \mathfrak{F}_{0} \),
the algebraic mapping cylinder in \( \mathrf {C}_{0} \) is that of \( \mathfrak{M}_{0} \)
with structure map pulled back over the diagonal morphism
\[
\mathfrak{S}\rightarrow \mathfrak{S}\otimes \mathfrak{S}\]
--- see \ref{th:soperad}.

Given a morphism \( f=s^{-1}g:C\rightarrow D \) in \( \mathrf {C} \), where
\( g:C\rightarrow X \) is a morphism in \( \mathrf {C}_{0} \) and \( s:D\rightarrow X \)
is an elementary equivalence, define \( \mathscr{M}(f)=\mathscr{M}(g) \), with
the canonical map \( \iota _{D}=\iota _{X}\circ s:D\rightarrow \mathscr{M}(f) \).
\end{defn}
\begin{prop}
The canonical map
\[
\iota _{D}:D\rightarrow \mathscr{M}(f)\]
 in diagram \ref{dia:mapcylpushout} is an elementary equivalence (in the sense
of \ref{def:emorph}).
\end{prop}
\begin{proof}
If \( f \) is a morphism in \( \mathfrak{M}_{0} \), \( \mathfrak{L}_{0} \),
or \( \mathfrak{F}_{0} \), this is clear: the inclusion is a direct summand
(as a chain-complex) and a homology equivalence. In the general case (i.e.,
in \( \mathfrak{M} \), \( \mathfrak{L} \), or \( \mathfrak{F} \)) \( \iota _{D} \)
is a composite of two elementary equivalences (in opposite directions).
\end{proof}
\begin{cor}
\label{cor:hequivanequiv}Let \( f:C_{1}\rightarrow C_{2} \) be a morphism
in \( \mathfrak{M} \), \( \mathfrak{L} \), or \( \mathfrak{F} \). Then \( f \)
is an equivalence if and only if its underlying map induces homology isomorphisms.
\end{cor}
\begin{proof}
The only-if part of this is clear. The if part follows immediately from the
\emph{existence} of algebraic mapping cylinders.
\end{proof}
\begin{cor}
\label{cor:singularequivalent}Let \( X \) be a pointed, simply-connected,
2-reduced simplicial set, let \( \ddelta (X) \) be the 2-reduced form of the
singular complex, and let \( i:X\rightarrow \ddelta (X) \) be the inclusion.
Then 
\[
\cf{i}:\cf{X}\rightarrow \cf{\ddelta (X)}\]
 is an elementary equivalence in \( \mathfrak{L}_{0} \), so that \( \cf{{X}} \)
and \( \cf{\ddelta (X)} \) define the same object of \( \mathfrak{L} \).
\end{cor}
\begin{rem}
In \cite{Smith:1994}, \( \cf{{\ddelta (X)}} \) was called \( \mathfrak{C}(X) \).
\end{rem}
\begin{thm}
\label{th:cfhomotopmfrak}The functor \index{functor}\( \cf{{*}}:\mathrm{SS}\rightarrow \mathfrak{L}_{0} \)
induces a functor
\[
\cf{{*}}:\underline{\mathrm{Homotop}}_{0}\rightarrow \mathfrak{L}\]
where \( \underline{\mathrm{Homotop}}_{0} \) is the homotopy category of pointed,
simply-connected, 2-reduced simplicial sets.
\end{thm}
\begin{rem}
By abuse of notation, we continue to denote this functor as \( \cf{{*}} \).
\end{rem}
\begin{proof}
This follows from results in \cite{Quillen:1967}, which characterizes the homotopy
category of simplicial sets as a localization of the category of simplicial
sets by weak equivalences. Proposition 4 on page 3.19 of \cite{Quillen:1967}
shows that a morphism \( f:X\rightarrow Y \) is an isomorphism in the homotopy
category if and only if:
\begin{enumerate}
\item \( f^{*}:H^{0}(Y,S)\rightarrow H^{0}(X,S) \) is an isomorphism for any set
\( S \)
\item \( f^{*}:H^{1}(Y,G)\rightarrow H^{1}(X,G) \) is an isomorphism for any group
\( G \)
\item \( f^{*}:H^{q}(Y,\mathbb{L})\rightarrow H^{q}(X,f^{*}\mathbb{L}) \) for an
local coefficient system \( \mathbb{L} \) of abelian groups on \( Y \) and
any \( q\geq 0 \)
\end{enumerate}
In our case (pointed, simply-connected, 2-reduced simplicial sets), these conditions
reduce to \( f \) inducing homology isomorphisms with \( \integers  \) coefficients.
Since, \( \underline{\mathrm{Homotop}}_{0}=\mathrm{SS}[V^{-1}] \), where \( V \)
is the morphisms inducing homology isomorphisms, it follows from \ref{cor:hequivanequiv}
that \( \cf{{*}} \) carries elements of \( V \) to the saturation of \( T \)
in \( \mathfrak{L}=\mathfrak{L}_{0}[S^{-1}] \), hence extending (\emph{uniquely})
to a functor
\[
\cf{{*}}:\underline{\mathrm{Homotop}}_{0}\rightarrow \mathfrak{L}\]

\end{proof}
\begin{defn}
\label{def:algebraicmapcone}If \( f:C\rightarrow D \) is a morphism in \( \mathfrak{M} \),
\( \mathfrak{L} \), or \( \mathfrak{F} \), define the \emph{algebraic mapping
cone}\index{algebraic mapping cone}, \( \mathscr{A}(f) \) to be the cofiber
of the cofibration \( \eta _{C}:C\rightarrow \mathscr{M}(f) \) defined in diagram
\ref{dia:mapcylpushout}. It comes equipped with a canonical morphism \( \iota _{D}:D\rightarrow \mathscr{A}(f) \).
\end{defn}
\begin{prop}
Let \( \mathrf {C}=\mathfrak{M} \), \( \mathfrak{L} \), or \( \mathfrak{F} \).
The algebraic mapping cone defines a functor from the category of \( \mathrf {C} \)-
morphisms to \( \mathrf {C} \).
\end{prop}
\begin{rem}
In other words, equivalent morphisms give rise to equivalent algebraic mapping
cones.
\end{rem}
\begin{proof}
Consider an elementary homotopy \[\xymatrix@R=10pt{{}&{X}\ar[dd]_{v}&{}\\{A}\ar[ur]^{f_1}\ar[dr]_{f_2}&&{B}\ar[ul]_{s_1}\ar[dl]^{s_2}\\{}&{Y}&{}}\] It
induces \[\xymatrix@R=10pt{{}&{X\otimes I}\ar[dd]_{v\otimes 1}&{}\\{A}\ar[ur]^{f_1}\ar[dr]_{f_2}&&{B}\ar[ul]_{s_1}\ar[dl]^{s_2}\\{}&{Y\otimes I}&{}}\] where
the map \( v\otimes 1 \) is an equivalence since \( v=s_{2}s_{1}^{-1} \) in
\( \mathrf {C} \). This factors through to the algebraic mapping cones.
\end{proof}
\begin{defn}
\label{def:suspension}If \( C \) is an m-coalgebra over an \( E_{\infty } \)-operad
\( \mathfrak{R} \), then the \emph{suspension}, \( \Sigma C \) , is \index{suspension}defined
via
\[
\Sigma C=C\otimes I/(C^{+}\otimes p_{0}\oplus C^{+}\otimes p_{1}\oplus C_{0}\otimes I)\]
where \( C^{+} \) is the portion of \( C \) above dimension \( 0 \) and \( C_{0}=\integers  \).
This is an m-coalgebra over \( \mathfrak{R}\otimes \mathfrak{S} \).
\end{defn}
\begin{rem}
It is not hard to see that this defines a functor \( \Sigma :\mathfrak{M}\rightarrow \mathfrak{M} \).
\end{rem}
\begin{lem}
\label{lem:algmappingcylextend}Let \( C_{1} \) and \( C_{2} \) be m-coalgebras
over \( \mathfrak{S} \) that are both acyclic and Cartan (see \ref{def:cartanmcoalgebra}
for the definition of the property of being Cartan) and let \( g:C_{1}\rightarrow C_{2} \)
be a chain map. Then there exists an m-structure on the algebraic mapping cylinder,
\( \mathscr{M}(g) \), extending those of \( C_{1} \) and \( C_{2} \).
\end{lem}
\begin{rem}
Note that the map \( f \) need \emph{not} be a morphism of m-coalgebras. The
m-structure on \( \mathscr{M}(g) \) measures the extent to which \( g \) \emph{fails}
to be a morphism of m-coalgebras.
\end{rem}
\begin{proof}
By abuse of notation, we will denote the algebraic mapping cylinder by \( \mathscr{M}(g) \)
(even though \( g \) is not necessarily a morphism). Let the chain complex
of \( I \) (the unit interval) be given by: 
\end{proof}
\begin{enumerate}
\item \( C_{0}(I)=\integers \cdot p_{0}\oplus \integers \cdot p_{1} \) (we have chosen
canonical generators, \( p_{1} \) and \( p_{2} \), of \( \integers ^{2} \))
and 
\item \( C_{1}(I)=\integers \cdot q \), with \( \partial q=p_{1}-p_{2} \).
\end{enumerate}
Let \( \varphi _{i}:C_{i}\rightarrow C_{i} \) be self-annihilating, contracting
homotopies associated with \( C_{1} \) and \( C_{2} \) (see \ref{def:cartanmcoalgebra}).
The following is a self-annihilating, contracting homotopy of \( \mathscr{M}(g) \):
\begin{eqnarray*}
\psi  & | & C_{1}\otimes p_{0}=\varphi _{1}\\
\psi  & | & C_{2}\otimes p_{1}=\varphi _{2}\\
\psi  & | & C_{1}\otimes q=\varphi _{1}\otimes q
\end{eqnarray*}
Now we construct an m-structure on \( \mathscr{M}(g) \) by a \emph{relativized}
form of the argument used in \ref{th:scrCdef}. Define
\[
\Psi _{n}=\psi \otimes 1\otimes \cdots \otimes 1+\epsilon \otimes \psi \otimes \cdots \otimes 1+\cdots +\epsilon \otimes \cdots \otimes \epsilon \otimes \psi \]
where each product has \( n \) factors. This is a self-annihilating contracting
homotopy of \( \mathscr{M}(g)^{n} \). We will build an m-structure over \( \mathfrak{S} \)
on \( \mathscr{M}(g) \) with adjoint map
\[
F_{n}:\rs{{n}}\otimes \mathscr{M}(g)\rightarrow \mathscr{M}(g)^{n}\]
that satisfies the invariant condition\emph{
\begin{equation}
\label{eq:mfinvariantcondition}
F_{n}(A(S_{n},1)\otimes \mathscr{M}(g))\subseteq \Psi _{n}(\mathscr{M}(g))
\end{equation}
}for all \( n>1 \). This m-structure is already defined on \( C_{1}\otimes p_{0} \)
and \( C_{2}\otimes p_{1} \); we need only extend it to \( C_{1}\otimes q \). 

Let \( 1\otimes A\in 1\otimes A(S_{n},1)\subset \rs{{n}}\subset \mathfrak{S} \)
and let \( z\in \im \psi \subset \mathfrak{M}(g) \) (i.e., \( z \) is either
in \( \im \varphi \otimes p_{1}\subset C_{2}\otimes p_{1} \) or \( C_{1}\otimes q \)).
We perform induction on the dimension of \( 1\otimes A\otimes z \) by setting
\[
F_{n}(\sigma \otimes A\otimes m)=\sigma \cdot F_{n}(A\otimes m)\]
 where \( \sigma \in S_{n} \), and \( m\in \mathscr{M}(g) \), and
\[
F_{n}(1\otimes A\otimes z)=\Psi _{n}\circ F_{n}\circ \partial (1\otimes A\otimes z)\]
 There is a unique chain-map satisfying these conditions because \( \psi ^{2}=0 \),
so that \( \partial |\im \psi  \) is \emph{injective}.

The proof that compositions in \( \mathfrak{S} \) map into compositions of
the \( F_{n} \)-maps proceeds \emph{exactly} as in \ref{th:scrCdef}, using
the invariant condition in \ref{eq:mfinvariantcondition} and a commutative
diagram like \ref{dia:coalgebradefined}.

\begin{cor}
\label{lem:homologyequival}Let \( \mathscr{T} \) be a category and let \( C(*) \)
and \( D(*) \) be functors from \( \mathscr{T} \) to \( \mathfrak{L}_{0} \).
Suppose that \( C(*) \) is free on some set of models \( \{m_{\alpha }\} \)
in \( \mathscr{T} \) and \( C(m_{\alpha }) \) and \( D(m_{\alpha }) \) are
acyclic \index{acyclic models}with Cartan m-structures (see \ref{def:cartanmcoalgebra})
for all models \( m_{\alpha } \). In addition, let \( f(\mathrm{ob}):C(\mathrm{ob})\rightarrow D(\mathrm{ob}) \)
be a natural transformation of functors for all \( \mathrm{ob}\in \mathscr{T} \).
Then \( f(\mathrm{ob}) \) is the underlying map of a morphism \( F(\mathrm{ob}):C(\mathrm{ob})\rightarrow D(\mathrm{ob}) \)
in \( \mathfrak{L} \) for all \( \mathrm{ob}\in \mathscr{T} \).

If, in addition, \( f(\mathrm{ob}) \) induces homology isomorphisms for all
\( \mathrm{ob}\in \mathscr{T} \), then it is the underlying map of an equivalence
in \( \mathfrak{L} \) for all \( \mathrm{ob}\in \mathscr{T} \).

\end{cor}
\begin{proof}
We  construct \( Z(\mathrm{o}b) \) to have an underlying chain complex that
is the \emph{algebraic mapping cylinder} of \( f(\mathrm{ob}) \). Since \( f \)
is a natural transformation, and since \( C \) is free on the \( \{m_{\alpha }\} \),
it follows that \( Z(*) \) is constructed from \( Z(f|C(m_{\alpha })) \).
Consequently, \ref{lem:homologyequival} implies that we can extend the m-structures
of \( C(m_{\alpha }) \) and \( D(m_{\alpha }) \) to an m-structure on \( Z(f|C(m_{\alpha })) \).
These m-structures induce an m-structure on on all of \( Z(f) \). 

With this m-structure, the inclusion of \( C(\mathrm{ob}) \) and \( D(\mathrm{ob}) \)
into the ends of this algebraic mapping cylinder are morphisms of m-coalgebras,
with the inclusion of \( D(\mathrm{ob}) \) an equivalence. The conclusion follows.
\end{proof}
We get a version of the Eilenberg-Zilber \index{Zilber@\textsc{J. A. Zilber}}\index{Eilenberg@\textsc{S. Eilenberg}}\index{Eilenberg-Zilber theorem}theorem:

\begin{cor}
\label{cor:eilenbergzilber}Let \( X \) and \( Y \) be pointed, simply-connected
simplicial sets. Then \( \cf{{X\times Y}}=\cf{{X}}\otimes \cf{{Y}}\in \mathfrak{L} \).
\end{cor}
\begin{itemize}
\item [\dbend] It is important to note that the Eilenberg-Zilber maps are \emph{not}
morphisms of coalgebras over \( \mathfrak{S} \) --- they are \emph{only} morphisms
in the localized category \( \mathfrak{L} \). 

This has a number of consequences --- not the least of which is the fact that
topological groups usually \emph{do not} define m-Hopf-algebras\index{m-Hopf algebra}.
Although their product map is associative, it does not preserve m-structures.
This means, for instance, that the reduced bar construction of \( \cf{{G}} \),
for \( G \) a simplicial group, does \emph{not} have a well-defined m-structure.

\end{itemize}
\begin{proof}
It is only necessary to note that these are functors from the category of \emph{ordered}
\emph{pairs} of simplicial sets and that the respective functors are free and
acyclic on models composed of pairs of simplices.
\end{proof}

\subsection{Homotopy sets and groups}

We can also define homotopy sets:

\begin{defn}
Let \( C \) and \( D \) be m-coalgebras. Then \( \mathfrak{L}_{0}(C,D) \)
and \( [C,D]=\mathfrak{L}(C,D) \) are defined to be the sets of morphisms in
the categories \( \mathfrak{L}_{0} \) and \( \mathfrak{L} \), respectively. 
\end{defn}
\begin{itemize}
\item [\dbend] One might be inclined to think of \( \mathfrak{L}(C,D) \) as a kind
of \emph{quotient} of \( \mathfrak{L}_{0}(C,D) \) (at least the author was
so inclined for some time). It really turns out to be a quotient of a \emph{proper
superset} of \( \mathfrak{L}_{0}(C,D) \). Indeed, it is \emph{easily} possible
for \( \mathfrak{L}_{0}(C,D)=0 \) and for \( [C,D]\neq 0 \). Furthermore,
an equivalence \( f:C\rightarrow C' \) does \emph{not} necessarily induce an
isomorphism \( \mathfrak{L}_{0}(C',D)\rightarrow \mathfrak{L}_{0}(C,D) \).
\end{itemize}
Here is an example of this phenomenon:

Let \( B^{n}=H_{*}(\cf{{S^{n}}},\integers ) \), for all \( n>1 \). This is
a chain-complex concentrated in dimensions \( 0 \) and \( n \), where its
chain-modules are \( \integers  \). It is not hard to see that the m-structure
on \( \cf{{S^{n}}} \) induces a (mostly trivial) m-structure on \( B^{n} \)
and that there exists an m-coalgebra morphism 
\begin{equation}
\label{eq:enmaps}
e_{n}:\cf{{S^{n}}}\rightarrow B^{n}
\end{equation}
Corollary~\ref{cor:hequivanequiv} implies that this is an \emph{equivalence}
in \( \mathfrak{M} \). 

Now we consider morphisms from \( B^{3} \) to \( B^{2} \). It is easy to see
that \( \mathfrak{M}_{0}(B^{3},B^{2})=0 \) since the chain-module of \( B^{3} \)
vanishes in dimension 2. 

On the other hand, let: 

\begin{enumerate}
\item \( Z_{3}=\mathscr{M}(e_{3}) \) --- the algebraic mapping cylinder of \( e_{3} \)
defined in \ref{eq:enmaps},
\item \( Z_{2}=\mathscr{M}(e_{2}) \) --- the algebraic mapping cylinder of 2 defined
in \ref{eq:enmaps},
\item \( V(h)=\mathscr{M}(h) \), where \( h:\cf{{S^{2}}}\rightarrow \cf{{S^{2}}} \)
is induced by a geometric map (the Hopf map, for instance).
\end{enumerate}
Now define \( K(h)=Z_{3}\cup _{\cf{{S^{3}}}}V(h)\cup _{\cf{{S^{2}}}}Z_{2} \).
The inclusion of \( B^{2} \) in \( K(h) \) is an elementary equivalence and
the diagram \[\xymatrix{{B^3}\ar[r]&{K(h)}&{B^2}\ar[l]}\] represents an element
of \( [B^{3},B^{2}] \). Two such diagrams represent the same element if they
fit into a diagram \[\xymatrix@R=10pt{&{K(h_1)}\ar[dd]^{g}&\\{B^3}\ar[ru]^{i_1}\ar[rd]_{i_2}&&{B^2}\ar[lu]_{j_1}\ar[ld]^{j_2}\\&{K(h_2)}&}\label{dia:equivmaphomotop}\] where
the \( i \) and \( j \)-maps are inclusions (and the \( j \)-maps are equivalences).

Now consider the relative cohomology \( H^{*}(K(h_{i}),B^{3}) \): this is given
by
\[
H^{i}(K(h_{i}),B^{3})=\left\{ \begin{array}{ll}
\integers  & i=0\\
\integers  & i=2\\
\integers  & i=4\\
0 & \mathrm{otherwise}
\end{array}\right. \]
where the generator, \( \alpha _{i} \), of \( H^{2}(K(h_{i}),B^{3}) \) is
induced by the canonical generator of \( H^{2}(B^{2}) \). Now \( \alpha _{i}\cup \alpha _{i}\in H^{4}(K(h_{i}),B^{3}) \)
is some multiple, \( w(h_{i}) \), of the generator, \( \beta _{i} \), of \( H^{4}(K(h_{i}),B^{3}) \). 

The commutativity of diagram~\ref{dia:equivmaphomotop} (and its \emph{m-structures})
implies that \( w(h_{1})=w(h_{2}) \), at least setting \( g^{*}(\beta _{2})=\beta _{1}\in H^{4}(K(h_{1}),B^{3}) \),
and a diagram chase involving the long exact sequence in cohomology shows that
\( w(h) \) coincides with the \emph{Hopf invariant} \index{Hopf invariant}of
the map \( h \). If \( h \) is a \emph{trivial} map (i.e., any map induced
from \( \mathfrak{M}_{0}(B^{3},B^{2}) \)), then \( w(h)=0 \), since \( K(h)/B^{3} \)
is essentially \( S^{4}\vee S^{2} \). 

On the other hand, the Hopf map \( S^{3}\rightarrow S^{2} \) is well-known
(see \index{Hopf@\textsc{H. Hopf}}\cite{Hopf:1935}) to have a Hopf invariant
of \( \pm 1 \). With suitable simplicial decompositions of \( S^{3} \) and
\( S^{2} \) (for instance, the singular complexes), we can compute the Hopf
invariant on the chain-level. 

This shows that \( \mathfrak{M}(B^{3},B^{2})\neq 0 \). With a little more work,
one can show that \( \mathfrak{M}(B^{3},B^{2})=\integers  \) (see chapter~\ref{ch:equivcategories}).

\begin{defn}
\label{def:mcoalgeonepointunion}Let \( C \) and \( D \) be m-coalgebras.
Define their \emph{one-point union}, \( C\vee D \) \index{m-coalgebra!one-point union}to
have an underlying chain-complex given by:
\begin{enumerate}
\item \( (C\vee D)_{0}=\integers  \)
\item \( (C\vee D)_{i}=C_{i}\oplus D_{i} \) for \( i>0 \)
\end{enumerate}
with m-structure induced from that of \( C \) and \( D \) via the canonical
projection 
\[
p:(C_{i}\oplus D_{i})^{n}\rightarrow C_{i}^{n}\oplus D_{i}^{n}\]

\end{defn}
\begin{prop}
\label{prop:mspherecohspace}For all \( n>1 \) there exists a morphism 
\[
d_{n}:B^{n}\rightarrow B^{n}\vee B^{n}\]
sending \( 1\in B^{n}_{n}=\integers  \) to \( 1\oplus 1\in B^{n}_{n}\oplus B^{n}_{n} \).
This is co-commutative and co-associative and, therefore, induces a map
\[
[B^{n},C]\times [B^{n},C]\rightarrow [B^{n},C]\]
making \( [B^{n},C] \) into an abelian group for any m-coalgebra, \( C \).
\end{prop}
\begin{rem}
Note that homotopy sets generally \emph{do not} have an abelian group structure:
the sum of two morphisms is usually \emph{not} a morphism due to the \emph{nonlinearity}
of m-structures.
\end{rem}
\begin{proof}
The key fact here is that the m-structure of \( B^{n} \) is \emph{trivial},
thereby eliminating any obstruction to the existence of \( d_{n} \). The statements
about co-commutativity and co-associativity follow immediately from the definition
of \( d_{n} \) and the remaining statements are clear.
\end{proof}
\begin{defn}
\label{def:mcoalghomotopygroup}Let \( C \) be an m-coalgebra and let \( n>1 \)
be an integer. Then \( \pi _{n}(C)=[B^{n},C] \) is defined to be the \( n^{\mathrm{th}} \)\emph{homotopy
group} of \( C \).\index{m-coalgebra!homotopy groups}\index{homotopy group of an m-coalgebra}
\end{defn}

\section{$\ainfty$-structures\label{sec:ainfinity}}

\subsection{Definitions}

In this section we define non-\( \Sigma  \) operads that can be used to define
``algebras and coalgebras that are associative up to a homotopy.'' These operads
were first described by Stasheff\index{Stasheff@\textsc{J. Stasheff}} in \cite{Stasheff:1963}
as a way of algebraically describing \( H \)-spaces whose product operation
was only \emph{homotopy associative}. 

\begin{defn}
\label{def:freeainfty} The free \( A_{\infty } \)-operad, \index{A@$\ainfty$ structure}\index{A@$\mathfrak{A}$!definition}\index{definition!$\mathfrak{A}$}\( \mathfrak{A} \),
is a non-\( \Sigma  \) operad equal to \( \mathfrak{A}=\hat{\afr }/\mathscr{I} \)
where: 
\begin{enumerate}
\item \( \hat{\afr } \) is the free, non-\( \Sigma  \) composition algebra (see~\ref{def:operadcomps})
generated over \( \mathbb{Z} \) by formal compositions of elements \( \{\mathfrak{D}_{k}\} \),
where \( \dim (\mathfrak{D}_{k})=k-2 \), \( \deg (\mathfrak{D}_{k})=k \).
This means that its \( \integers  \)-generators are words \( \{\D _{i_{1}}\circ _{j_{1}}\cdots \circ _{k_{k-1}}\D _{i_{k}}\} \)
subject only to the two relations 

\begin{description}
\item [Associativity]\( (a\circ _{i}b)\circ _{j}c=a\circ _{i+j-1}(b\circ _{j}c) \)
\item [Commutativity]\( a\circ _{i+m-1}(b\circ _{j}c)=(-1)^{mn}b\circ _{j}(a\circ _{i}c) \)
\end{description}
where \( a,b,c\in \hat{\afr } \).

\item \( \mathscr{I}\subset \hat{\afr } \) is the ideal generated by 
\begin{equation}
\label{eq:ainfinityidentity}
\sum _{k=1}^{n}\sum _{\lambda =0}^{n-k}(-1)^{k+\lambda +k\lambda }\mathfrak{D}_{k}\circ _{n-k-\lambda +1}\mathfrak{D}_{n-k+1}
\end{equation}
 for all \( n\geq 1 \). 
\end{enumerate}
The identity in \ref{eq:ainfinityidentity} implies that \( \mathfrak{D}_{1}^{2}=0 \).
We make \( \mathfrak{A} \) a sequence of DG-modules by setting boundary operators
to \( \mathfrak{D}_{1} \).
\begin{defn}
\label{def:ainfinitystructure} Let \( \mathfrak{R} \) be an operad. An \( A_{\infty } \)-\emph{structure}
on\index{A@$\ainfty$-structure} \( \mathfrak{R} \) is a morphism of (non-\( \Sigma  \))
operads: 
\[
\mathfrak{A}\rightarrow \mathfrak{R}\]
 If \( C \) is a DG-module, 
\end{defn}
\begin{enumerate}
\item an \( A_{\infty } \)-algebra structure on \( C \) is an \( A_{\infty } \)-structure
on the endomorphism operad, i.e., an operad morphism 
\[
\mathfrak{A}\rightarrow \zend (C)\]

\item an \( A_{\infty } \)-coalgebra structure on \( C \) is an \( A_{\infty } \)-structure
on the co-endomorphism operad, i.e., an operad morphism 
\[
\mathfrak{A}\rightarrow \coend (C)\]

\end{enumerate}
\end{defn}
\begin{rem}
We briefly consider what \( A_{\infty } \)-algebras and coalgebras look like.
Suppose \( A \) is an \( A_{\infty } \)-algebra with structure map 
\[
f:\mathfrak{A}\rightarrow \zend (C)\]
 Then \( A \) is equipped with maps 
\[
\highprod{i}=f(\mathfrak{D}_{i}):A^{i}\rightarrow A\]
 of degree \( i-2 \) satisfying the identity 
\[
\sum _{k=1}^{n}\sum _{\lambda =1}^{n-k}\highprod{n-k+1}\circ (1^{\lambda }\otimes \highprod{k}\otimes 1^{n-k-\lambda })=0\]
 where \( \highprod{1}=\partial :A\rightarrow A \). Similarly, an \( A_{\infty } \)-coalgebra
with structure map 
\[
f:\mathfrak{A}\rightarrow \coend (C)\]
 comes equipped with a system of degree \( k-2 \)-maps 
\[
f(\mathfrak{D}_{k})=\delta _{k}:C\rightarrow C^{k}\]
 satisfying the identity 
\[
\sum _{k=1}^{n}\sum _{\lambda =1}^{n-k}(-1)^{k+\lambda +k\lambda }\delta _{n-k+1}\circ (1^{\lambda }\otimes \delta _{k}\otimes 1^{n-k-\lambda })=0\]

\end{rem}
An \( A_{\infty } \)-structure on an operad trivially induces corresponding
structures on algebras and coalgebras over the operad:

\begin{prop}
\label{prop:ainfinitythings} Let \( \mathfrak{R} \) be an operad equipped
with an \( A_{\infty } \)-structure. Then 
\begin{enumerate}
\item A coalgebra over \( \mathfrak{R} \) is an \( A_{\infty } \)-coalgebra. 
\item A algebra over \( \mathfrak{R} \) is an \( A_{\infty } \)-algebra. 
\end{enumerate}
\end{prop}
These two constructs are dual to each other, as the following result shows:

\begin{prop}
\label{prop:ainftyalgcoalgdual}Let \( C \) be an \( \ainfty  \)-coalgebra
and let \( A \) be a ring. Then \( \homz (C,A) \) has a natural \( \ainfty  \)-algebra
structure.
\end{prop}
\begin{proof}
This follows from the Duality Theorem (\ref{lem:operadduality}), which states
the existence of an operad morphism:
\[
\coend (C)\rightarrow \zend _{\{a_{n}'\}}(\homz (C,A))\]
Now we compose the structure map of \( C \) with this to get an operad morphism
\[
\mathfrak{A}\rightarrow \zend _{\{a_{n}\}}(\homz (C,A))\]
The conclusion follows by noting that, as a non-\( \Sigma  \) operad \( \zend _{\{a_{n}\}}(\homz (C,A))=\zend (\homz (C,A)) \).
\end{proof}

\subsection{Twisted tensor products\label{sec:twistedtensoralgebraic}}

Now we develop the important algebraic concept of a twisting cochain and twisted
tensor product. We begin with the ``categorical definition'' given in \cite{Gugenheim:1972}\index{Gugenheim@\textsc{Victor K. A. M. Gugenheim}}:

\begin{defn}
\label{def:twistedtensorcategorical}Let \index{definition!twisted tensor product}\index{twisted tensor product!definition}\( C \)
be a DGA-coalgebra with coproduct \( \Delta :C\rightarrow C\otimes C, \) and
let \( A \) be a DGA-algebra with product \( \mu :A\otimes A\rightarrow A \).
A \emph{principal twisted tensor product,} \( Z \), of \( C \) by \( A \)
is a chain-complex satisfying:
\begin{enumerate}
\item \( Z \) is isomorphic, as a graded module, to \( C\otimes A \)
\item \( Z \) is a left, DG \( C \)-comodule --- i.e., there exists a morphism of
DG-modules
\[
\delta :Z\rightarrow C\otimes Z\]
 making \[\xymatrix@C+10pt{{Z}\ar@{=}[d] \ar[r]^-{\delta}&{C\otimes Z}\ar[r]^-{1\otimes \delta \circ \delta}&{C\otimes C\otimes Z}\ar@{=}[d]\\ {Z}\ar[r]_-{\delta}&{C\otimes Z}\ar[r]_-{\Delta \otimes 1}&{C\otimes C\otimes Z}}\] commute;
\item \( Z \) is a right, DG \( A \)-module --- i.e., there exists a morphism of
DG-modules
\[
\alpha :Z\otimes A\rightarrow A\]
making \[\xymatrix@C+10pt{{Z\otimes A\otimes A}\ar@{=}[d]\ar[r]^-{1\otimes \alpha}&{Z\otimes A}\ar[r]^-{\alpha}&{Z}\ar@{=}[d]\\ {Z\otimes A\otimes A}\ar[r]_-{1 \otimes \mu}&{Z\otimes A}\ar[r]_-{\alpha}&{Z}}\]
commute.
\end{enumerate}
\end{defn}
\begin{rem}
This is very similar to Cartan's\index{Cartan@\textsc{Henri Cartan}} definition
of a Construction: the essential difference is the present requirement that
\( Z \) be a \emph{comodule} over \( C \).
\end{rem}
Given this definition, one can readily determine the types of \emph{differentials}
that can exist on a twisted tensor product:

\begin{prop}
Let \( Z=C\otimes A \) be \index{twisted tensor product!differential}a principal
twisted tensor product, as defined in \ref{def:twistedtensorcategorical}. Then
the differential of \( Z \) is given by
\[
\partial _{Z}=\partial _{C\otimes A}+x\cap *\]
where \( x:C\rightarrow A \) is a cochain (called the twisting cochain) satisfying
the identity
\[
\partial x+x\cup x=0\]
This twisted tensor product is written \( C\otimes _{x}A \).
\end{prop}
\begin{proof}
See \cite{Gugenheim:1972}\index{Gugenheim@\textsc{Victor K. A. M. Gugenheim}}
for a proof. This involves considering \( \partial _{Z}|C\otimes 1\cap \ker \delta  \)
and applying the identities that modules and comodules must satisfy.
\end{proof}
Definition~\ref{def:twistedtensorcategorical} easily generalizes to the case
where \( C \) is an \( \ainfty  \)-coalgebra, \( A \) is an \( \ainfty  \)-algebra,
and \( Z \) is an \( \ainfty  \)-comodule over \( C \) and an \( \ainfty  \)-module
over \( A \). 

In the case where \( C \) is \( \ainfty  \) and \( A \) is a DGA-algebra,
we get:

\begin{defn}
\label{def:twistingcochaincobar}Let \( A \) be a DGA-algebra, let \( C \)
be an \( \ainfty  \)-coalgebra, and let \( x:C\rightarrow A \) be a map of
degree \( -1 \). Then this map is a \emph{twisting cochain} if
\[
\partial _{A}\circ x+\sum _{i=1}^{\infty }\mu ^{i-1}\circ \underbrace{{x\otimes \cdots \otimes x}}_{i\, \mathrm{factors}}\circ \Delta _{i}=0\]
 where \( \mu :A\otimes A\rightarrow A \) is the multiplication and \( \Delta _{i}:C\rightarrow C^{i} \)
defines the \( \ainfty  \)-structure. The map \( \partial _{x}:C\otimes A\rightarrow C\otimes A \)
given by
\[
\partial _{x}=1\otimes \partial _{A}+\sum _{i=1}^{\infty }(1\otimes \mu ^{i-1})\circ (1\otimes x^{i-1}\otimes 1)\circ (\Delta _{i}\otimes 1)\]
is a differential. The chain-complex, \( C\otimes A \), equipped with this
differential, is a \emph{principal twisted tensor product\index{twisted tensor product!principal}}
and denoted \( C\otimes _{x}A \).
\end{defn}
Another variant of twisted tensor products is:

\begin{defn}
\label{d:ainftyalgtwistedtensor} Let \( C \) be a DGA-coalgebra and let \( F \)
be an \( A_{\infty } \)-algebra with structure morphism 
\[
f:\mathfrak{A}\rightarrow \zend (F)\]

\begin{enumerate}
\item A map \( x:C\rightarrow F \), of degree \( -1 \) will be called a \emph{twisting
cochain} if it satisfies the condition 
\[
x\circ \partial _{C}+\sum _{i=1}^{\infty }f(\mathfrak{D}_{i})\circ \underbrace{{x\otimes \cdots \otimes x}}_{i\, \mathrm{factors}}\circ \Delta ^{i-1}=0\]

\item Given a twisting cochain, \( x \), the map 
\[
b(x)=\epsilon +\sum _{i=1}^{\infty }\underbrace{{x\otimes \cdots \otimes x}}_{i\, \mathrm{factors}}\circ \Delta ^{i-1}:C\rightarrow \barcs (F)\]
 is a homomorphism of DGA-coalgebras. 
\item Given a twisting cochain, \( x \), the map \begin{multline*}\partial_{x} = 1\otimes\partial_{C} + \sum_{i=1}^{\infty}( f(\mathfrak{D}_{i})\otimes1)\circ1\otimes\underbrace{x\otimes\cdots\otimes x}_{\text{\(i-1\) factors}}\otimes1\circ\Delta^{i-1}\otimes1 \\ C\otimes F \to C\otimes F \end{multline*}
is self-annihilating. 
\end{enumerate}
The chain-complex, \( \barcs (F)\otimes C \), equipped with the differential
\( \partial _{x} \) will be called the \emph{\( x \)-twisted tensor product}
and denoted \( \barcs (F)\otimes _{x}C \).
\end{defn}
\begin{rem}
As in the previous definition of a twisted tensor product, we have incorporated
the differential of one of the chain-complexes (namely \( A \)) via an \( A_{\infty } \)-structure.
Here \( f(\mathfrak{D}_{1})=\partial _{F} \). See \cite[\S~3]{Proute:1}\index{Proute@\textsc{Alain Prout{\'e}}}
for definitions and the proof that a map \( b:C\rightarrow \barcs (F) \) is
a homomorphism of DGA-coalgebras if and only if it is of the form \( b(x) \)
for some twisting cochain \( x \). Note that \( x(C_{0})=0 \) so that the
0-dimensional components of the coproduct, \( \Delta  \), have no effect on
the formula for \( b(x) \) or the definition of a twisting cochain. 

We can generalize this definition of twisted tensor product somewhat. Let \( M \)
be a left \( A_{\infty } \)-algebra over \( A \), i.e. suppose there exist
maps \( \{\ainfaction {k}\} \) of degree \( k-2 \), \( \ainfaction {k}:A\otimes \cdots \otimes A\otimes M\rightarrow M \)
(\( k \)-factors in all), otherwise satisfying the same identities as an \( A_{\infty } \)-product.
This means that the equation in \ref{def:ainfinitystructure} is satisfied with
the rightmost copy of \( \highprod{k} \) replaced by \( \ainfaction {k} \).
In this case we can define a twisted tensor product \( C\otimes _{\ainfaction {k}}M \).
\end{rem}

\subsection{m-Hopf algebras}

\begin{defn}
\label{def:mhopfalg}Let: 
\end{defn}
\begin{enumerate}
\item \( \mathfrak{H}_{0} \) denote the category of \emph{m-Hopf algebras\index{m-Hopf algebra!definition}\index{definition!m-Hopf algebra}}
in \( \mathfrak{M}_{0} \). Its objects are the monoid objects of \( \mathfrak{M}_{0} \)
that satisfy the additional conditions:

\begin{enumerate}
\item underlying operads are operad-coalgebras (see \ref{def:operadcoproduct})
\item the \emph{product} operation of an m-Hopf algebra, \( C, \) \( \mu :C\otimes C\rightarrow C \)
covers the underlying operad's \emph{coproduct} \( \Delta :\mathfrak{R}\rightarrow \mathfrak{R}\otimes \mathfrak{R} \).
\end{enumerate}
Morphisms are required to preserve \emph{all structures} (i.e. the product operation
and the underlying operad's coproduct).

\item \( \mathfrak{H} \) denote the category of m-Hopf algebras in \( \mathfrak{M} \).
Its objects are objects of \( \mathfrak{M} \) over \( E_{\infty } \)-operad-coalgebras
that have a product (a morphism of \( \mathfrak{M} \), now) \( \mu :C\otimes C\rightarrow C \)
that covers the underlying operad's coproduct.
\item \( \mathfrak{G} \) denote the category of group objects in \( \mathfrak{M} \).
\end{enumerate}
\begin{rem}
If \( C \) is a monoid object of \( \mathfrak{M}_{0} \) over the \( E_{\infty } \)-operad
\( \mathfrak{R} \), the existence of a multiplication on \( C \) implies (by
\ref{def:m-coalgebra} and \ref{def:einftycomoduleprod}) the existence of an
\( E_{\infty } \)-operad \( \mathfrak{R}' \) and morphisms \( \mathfrak{R}'\rightarrow \mathfrak{R} \)
and \( \mathfrak{R}'\rightarrow \mathfrak{R}\otimes \mathfrak{R} \). It might
seem unnecessarily restrictive to require that \( \mathfrak{R}'=\mathfrak{R} \)
but appears to be necessary to define the bar construction (since we must take
\emph{iterated products} of elements).
\end{rem}
The DGA-algebras, \( A(M,n) \), defined by Eilenberg\index{Eilenberg@\textsc{S. Eilenberg}}
and MacLane\index{MacLane@\textsc{S. MacLane}} in \cite{EM1} and \cite{EM2}
turn out to have natural m-Hopf algebra structures over the operad \( \mathfrak{S} \)
defined in \S~\ref{sec:sfrakdef}. These DGA-algebras were originally used to
compute the homology of Eilenberg-MacLane spaces --- and in our theory, they
act very much like Eilenberg-MacLane spaces.

\section{The bar and cobar constructions\label{sec:barcobarconstructions}}

\subsection{Introduction}

In this section, we define two important constructs in algebraic homotopy: the
bar and cobar constructions. 

The \emph{bar construction} \index{construction!bar}was first defined by Eilenberg
\index{Eilenberg@\textsc{S. Eilenberg}}and MacLane \index{MacLane@\textsc{S. MacLane}}in
\cite{EM1} and \cite{EM2}. Its defining property is that the bar construction,
\( \barcs (C(G)) \), of the chain complex of a simplicial group, \( G \),
gives the chain complex of the \emph{classifying space} of the group. Eilenberg
and MacLane used it to compute the homology of Eilenberg-MacLane spaces. 

The \emph{cobar construction}\index{construction!cobar}, \( \cobar C \), was
defined by Adams in \index{Adams@\textsc{J. F. Adams}}\cite{Adams:1956} in
a manner dual to the definition of the bar construction. Its characteristic
property is that the cobar construction of the chain complex of a pointed, simply-connected
simplicial set is the chain-complex of the \emph{loop space}. Together the bar
and cobar constructions allow one to understand fibrations at the chain-complex
level.

In \S~\ref{subsec:mapsequencecobar} and \ref{subsec:mapseqbar}, we will give
slightly nonstandard definitions of the bar and cobar construction in terms
of mapping telescopes (see \S~\ref{sec:mappingseqteles}) and prove a duality
theorem:

\emph{Proposition (\ref{prop:bardualcobar}): Let \( C \) be an \( \ainfty  \)-coalgebra
with structure morphism 
\[
f:\mathfrak{A}\rightarrow \coend (C)\]
and let \( A \) be a DGA-algebra. Then the \( \ainfty  \)-coalgebra structure
on \( C \) induces an \( \ainfty  \)-algebra structure on the hyper-chain
complex \( \homz (C,A) \) (via the cup-product map in \ref{prop:ainftyalgcoalgdual})
such that there exists a natural morphism of chain complexes
\[
\barcs (\homz (C,A))\rightarrow \homz (\cobar (C),A)\]
}

This will play an important part in our computation of the m-structure of the
cobar construction in \S~\ref{sec:cobar}.

\subsection{Mapping sequences and telescopes\label{sec:mappingseqteles}}

In this section we will define \emph{mapping telescopes} --- a simple generalization
of algebraic mapping cones --- that we will use in throughout the remainder
of this paper. We recall the familiar definition of algebraic mapping cone:

\begin{defn}
\label{def:rightalgmappingcone}Given a chain-map \( f:C\rightarrow D \) between
chain-complexes, the \emph{right algebraic mapping cone}\index{right algebraic mapping cone},
\( C_{R}(f) \), is the chain-complex \( C\oplus \Sigma ^{-1}D \), equipped
with the differential 
\[
\left( \begin{array}{cc}
\partial _{C} & 0\\
\desusp \circ f & -\desusp \circ \partial _{D}\circ \susp 
\end{array}\right) :C\oplus \Sigma ^{-1}D\rightarrow C\oplus \Sigma ^{-1}D\]
This comes equipped with a canonical chain map \( C_{R}(f)\rightarrow C \)
and a canonical inclusion \( \Sigma ^{-1}D\subset C_{R}(f) \).

The \emph{left algebraic mapping cone}\index{left algebraic mapping cone},
\( C_{L}(f) \), is just \( \Sigma C_{L}(f) \) with boundary map
\[
\left( \begin{array}{cc}
-\susp \circ \partial _{C}\circ \desusp  & 0\\
f\circ \desusp  & \partial _{D}
\end{array}\right) :\Sigma C\oplus D\rightarrow \Sigma C\oplus D\]
This comes with a canonical inclusion \( D\subset C_{L}(f) \) and projection
\( C_{L}(f)\rightarrow \Sigma C \)
\end{defn}
\begin{rem}
Compare the Koszul Convention and our assumption that \( \susp  \) and \( \desusp  \)
are chain-maps, which implies that the boundary of \( \Sigma ^{-1}C \) is \( -\desusp \circ \partial _{C}\circ \susp  \).
\end{rem}
The following result is an immediate consequence of the definitions:

\begin{prop}
Given \( f:C\rightarrow D \) with right algebraic mapping cone \( C_{R}(f) \)
and a chain-complex \( E \) with chain-map \( e:E\rightarrow C \), we can
define a map 
\[
e\oplus \desusp \circ g:E\rightarrow C_{R}(f)=C\oplus \Sigma ^{-1}D\]
if and only if \( g \) is a null-homotopy of \( f\circ e \). If this condition
is satisfied, we may form an iterated algebraic mapping cone
\[
C_{R}(e\oplus g)=E\oplus \Sigma ^{-1}C\oplus \Sigma ^{-2}D\]

\end{prop}
\begin{defn}
\label{def:rightmappingsequencetele}A \emph{right mapping sequence}\index{right mapping sequence}\index{mapping sequence!right},
\( \{C_{i}\}_{\rightarrow } \), is an infinite sequence of chain-complexes
and maps \[\xymatrix{{C_0}\ar[r]^{\{f_{0,*}\}}&{C_1}\ar[r]^{\{f_{1,*}\}} &{C_2}\ar[r]^{\{f_{2,*}\}}&{\cdots}}\] where:
\begin{enumerate}
\item \( f_{i,0}:C_{i}\rightarrow C_{i+1} \) is a chain map with right mapping cone
\( C_{R}(f_{i,0})=F_{i,1} \)
\item \( f_{i,1}:C_{i}\rightarrow C_{i+2} \) is a nullhomotopy of \( f_{i+1,0}\circ f_{i,0} \)
--- inducing a chain-map \( f_{i,0}\oplus \desusp \circ f_{i,1}:C_{i}\rightarrow F_{i+1,1} \)
\item the map \( f_{i,j}:C_{i}\rightarrow C_{i+j+1} \) has the property that 
\[
(f_{i,1}\oplus \desusp \circ f_{i,2,}\oplus \cdots \oplus \desusp ^{j-2}\circ f_{i,j-1})\oplus \desusp ^{j-1}\circ f_{i,j}:C_{i}\rightarrow F_{i+2,j-1}\]
is a nullhomotopy of 
\[
(f_{i+1,0}\oplus \desusp \circ f_{i+1,1}\oplus \cdots \oplus \desusp ^{j-1}\circ f_{i+1,j-1})\circ f_{i,0}:C_{i}\rightarrow F_{i+2,j-1}\]
inducing a chain-map \( f_{i,0}\oplus \desusp \circ f_{i,1}\oplus \cdots \oplus \desusp ^{j}\circ f_{i,j}:C_{i}\rightarrow F_{i+1,j} \)
with right mapping cone \( C_{R}(f_{i,0}\oplus \desusp \circ f_{i,1}\oplus \cdots \oplus \desusp ^{j}\circ f_{i,j})=F_{i,j+1} \).
\end{enumerate}
The \( \{F_{i,j}\} \) come equipped with canonical surjections \( F_{i,j}\leftarrow F_{i,j'} \)
for all \( i \) and all \( j'>j \). Define the associated \emph{mapping telescope}
\index{mapping telescope}via \( T\{C_{i}\}_{\rightarrow }=\varprojlim F_{0,*} \).
Let \( \mathfrak{T}_{R} \) denote the category of right mapping sequences.

A \emph{morphism} of mapping sequences is a sequence of morphisms of corresponding
chain-complexes that commute with all of the \( \{f_{i,j}\} \) (whenever composites
exist).

\end{defn}
\begin{rem}
Given a ``chain-complex of chain-complexes'' --- a sequence of chain-complexes
and chain-maps whose composites are zero --- we can form an ``iterated algebraic
mapping cone''. This is nothing but the hyperhomology complex. 

A right sequence as defined above is a ``\emph{homotopy chain-complex of chain-complexes}'':
two successive chain-maps might have composites that are not zero but \emph{null-homotopic}.
We equip the structure with whatever it needs to make the notion of iterated
algebraic mapping cone well-defined. This turns out to be precisely the ``higher
null-homotopies'' \( \{f_{i,j}\} \). The mapping telescope is just the hyperhomology
complex.

In practice, the \( \{f_{i,j}\} \) will have the property that \( f_{i,j}=0 \)
for all but a finite number of values of \( j \): this causes the inverse limit
to be convergent (so the mapping telescope is countably generated).
\end{rem}
\begin{prop}
\label{prop:maptelescopeformula}Given the definitions in \ref{def:rightmappingsequencetele},
the boundary map of \( F_{i,j} \) is a \( j+1\times j+1 \) array:
\begin{equation}
\label{eq:mapsequenceboundary}
\left( \begin{array}{cccc}
\partial _{i} & 0 & \cdots  & 0\\
\desusp f_{i,0} & -\desusp \partial _{i+1}\susp  & \cdots  & 0\\
\desusp ^{2}f_{i,1} & -\desusp ^{2}f_{i+1,0}\susp  & \cdots  & 0\\
\vdots  & \vdots  & \ddots  & \vdots \\
\desusp ^{j}f_{i,j-1} & -\desusp ^{j}f_{i+1,j-2}\susp  & \cdots  & (-1)^{j}\desusp ^{j}\partial _{i+j}\susp ^{j}
\end{array}\right) 
\end{equation}
 where \( \partial _{i} \) is the boundary map of \( C_{i} \). The general
element is given by 
\begin{equation}
\label{eq:mappingsequenceboundary2}
b_{\alpha ,\beta }=\left\{ \begin{array}{cl}
0 & \mathrm{if}\, \alpha <\beta \\
(-1)^{\alpha -1}\desusp ^{\alpha -1}\partial _{i+\alpha -1}\susp ^{\alpha -1} & \mathrm{if}\, \alpha =\beta \\
(-1)^{\beta -1}\desusp ^{\alpha -1}f_{i+\beta -1,\alpha -\beta -1}\susp ^{\beta -1} & \mathrm{if}\, \alpha >\beta 
\end{array}\right. 
\end{equation}
and the condition on \( f_{i,j} \) can be written as:
\begin{eqnarray*}
\lefteqn {\sum _{\alpha =2}^{j}(-1)^{\alpha }f_{i+\alpha ,j-\alpha }\circ f_{i,\alpha -1}+(-1)^{j-1}\partial _{i+j+1}\circ f_{i,j}} &  & \\
 &  & \hspace {1in}+f_{i,j}\circ \partial _{i}=f_{i+1,j-1}\circ f_{i,0}
\end{eqnarray*}

\end{prop}
\begin{proof}
The equations \ref{eq:mapsequenceboundary} and \ref{eq:mappingsequenceboundary2}
follow by an inductive application of the boundary formula in \ref{def:rightalgmappingcone}.
The formula for \( f_{i,j} \) follows by a straightforward but tedious computation.
\end{proof}
It is straightforward to dualize \ref{def:rightmappingsequencetele} to get
the dual concept of \emph{left mapping sequence} and telescope:

\begin{defn}
\label{def:leftmappingsequencetele}A \emph{left mapping sequenc}e\index{left mapping sequence}\index{mapping sequence!left},
\( \{C_{i}\}_{\leftarrow } \), is an infinite sequence of chain-complexes and
chain-maps \[\xymatrix{{\cdots}& {C_2}\ar[r]^{\{f_{2,*}\}}&{C_1}\ar[r]^{\{f_{1,*}\}}
&{C_0}}\] where:
\end{defn}
\begin{enumerate}
\item \( f_{i,0}:C_{i+1}\rightarrow C_{i} \) is a chain map with left mapping cone
\( C_{L}(f_{i,0})=F_{i,1} \)
\item \( f_{i,1}:C_{i+2}\rightarrow C_{i} \) is a nullhomotopy of \( f_{i,0}\circ f_{i+1,0} \)
--- inducing a chain-map \( \left( \begin{array}{c}
f_{i,1}\circ \desusp \\
f_{i,0}
\end{array}\right) :F_{i+1,1}\rightarrow C_{i} \)
\item the map \( f_{i,j}:C_{i+j+1}\rightarrow C_{i} \) has the property that 
\[
\left[ \begin{array}{c}
f_{i,j}\circ \desusp ^{j-1}\\
f_{i,j-1}\circ \desusp ^{j-2}\\
\vdots \\
f_{i,1}
\end{array}\right] :F_{i+2,j-1}\rightarrow C_{i}\]
is a nullhomotopy of 
\[
f_{i,0}\circ \left[ \begin{array}{c}
f_{i+1,j-1}\circ \desusp ^{j-1}\\
\vdots \\
f_{i+1,1}\circ \desusp \\
f_{i+1,0}
\end{array}\right] :F_{i+2,j-1}\rightarrow C_{i}\]
inducing a chain-map \( \left[ \begin{array}{c}
f_{i,j}\circ \desusp ^{j}\\
\vdots \\
f_{i,1}\circ \desusp \\
f_{i,0}
\end{array}\right] :F_{i+1,j}\rightarrow C_{i} \) with left mapping cone \( C_{L}=F_{i,j+1} \).

The \( \{F_{i,j}\} \) come with inclusions \( F_{i,j}\subset F_{i,j'} \) for
all \( i \) and \( j'>j \). Define the \emph{mapping telescope} via \( T\{C_{i}\}_{\leftarrow }=\varinjlim F_{i,i} \).

\end{enumerate}
\begin{defn}
Let \( \mathfrak{T}_{L} \) denote the category of left mapping sequences.
\end{defn}
The following result is straightforward:

\begin{prop}
Morphisms of (left or right) mapping sequences induce morphisms of their corresponding
mapping telescopes (as chain-complexes).

Given a mapping sequence \( \{C_{i}\}_{*} \), where \( *=\rightarrow  \) or
\( \leftarrow  \), and a chain-complex \( D \), one may define the induced
mapping sequence \( \mathrm{Hom}(D,\{C_{i}\}_{*})=\{\homz (D,C_{i})\}_{*} \). 
\end{prop}
We also have the following relation between right and left mapping sequences:

\begin{prop}
\label{prop:mappingtelesdual}Let \( \{C_{i}\}_{\rightarrow } \) be a right
mapping sequence and let \( D \) be a DGA-algebra. Then \( \mathrm{Hom}(\{C_{i}\}_{\rightarrow },D)=\{\homz (C_{i},D\}_{\leftarrow } \)
has the structure of a left mapping sequence and there exists a chain-map 
\[
T\{\homz (C_{i},D\}_{\leftarrow }\rightarrow \homz (T\{C_{i}\}_{\rightarrow },D)\]
that is natural with respect to morphisms of mapping sequences.
\end{prop}
\begin{proof}
Most of the statements follow from the functoriality of \( \homz (*,*) \).
The last statement follows by considering the \( F_{i,j} \) and passing to
the limits: it is not hard to see that \( F_{i,j}(\{\homz (C_{*},D\}_{\leftarrow }\})\cong \homz (F_{i,j(}T\{C_{*}\}_{\rightarrow }),D) \)
for finite \( i,j \). When we pass to the limits, there still exists a map,
but it is not necessarily an isomorphism. 
\end{proof}
We can define \emph{tensor products} of mapping sequences:

\begin{prop}
\label{prop:maptelescopeproduct}Let \( \{C_{i},f_{*,*}\} \) and \( \{D_{j},g_{*,*}\} \)
be two right mapping sequences. If we define their tensor product, \( \{C_{i}\}\otimes \{D_{j}\} \)
to be the right mapping sequence 
\begin{equation}
\label{eq:maptelescomplex}
\{(\{C_{i}\}\otimes \{D_{j}\})_{k}=\bigoplus _{i+j=k}C_{i}\otimes D_{j},h_{*,*}\}
\end{equation}
where 
\begin{equation}
\label{eq:mapteleprodmap}
h_{i,j}|C_{\alpha }\otimes D_{\beta }=f_{\alpha ,j}\otimes 1+(-1)^{\alpha (j+1)}1\otimes g_{\beta ,j}
\end{equation}
then \( \iota :T\{C_{i}\}_{\rightarrow }\otimes \{D_{j}\}_{\rightarrow }\rightarrow T\{C_{i}\}\otimes T\{D_{j}\}_{\rightarrow } \)
is an isomorphism, where \( \iota |\Sigma ^{-i}C_{i}\otimes \Sigma ^{-j}D_{j}=\desusp ^{i+j}_{C\otimes D}\circ \susp ^{i}_{C}\otimes \susp ^{j}_{D}:\Sigma ^{-i}C_{i}\otimes \Sigma ^{-j}D_{j}\rightarrow \Sigma ^{-(i+j)}(C\otimes D) \)
(see \ref{prop:suspisos}).
\end{prop}
\begin{rem}
Note that \ref{eq:maptelescomplex} implies that \( i=\alpha +\beta  \) in
equation \ref{eq:mapteleprodmap}.
\end{rem}
\begin{proof}
Recall that, in the notation of \ref{def:rightmappingsequencetele}, \( T\{*\}_{\rightarrow }=F_{0,\infty } \).
The boundary of \( T\{C_{i}\}_{\rightarrow }\otimes \{D_{j}\}_{\rightarrow } \)
is \( \partial _{1}\otimes 1+1\otimes \partial _{2} \), where \( \partial _{1} \)
is given, on \( \Sigma ^{-\alpha }C_{\alpha } \) by (see \ref{eq:mappingsequenceboundary2}
)
\[
\sum _{\alpha '\geq \alpha }(-1)^{\alpha }\desusp ^{\alpha '}f_{\alpha ,\alpha '-\alpha -1}\susp ^{\alpha }:\Sigma ^{-\alpha }C_{\alpha }\rightarrow \bigoplus _{\alpha '\geq \alpha }\Sigma ^{-\alpha '}C_{\alpha '}\]
and \( \partial _{2} \) is given by 
\[
\sum _{\beta '\geq \beta }(-1)^{\beta }\desusp ^{\beta '}g_{\beta ,\beta '-\beta -1}\susp ^{\beta }:\Sigma ^{-\beta }C_{\beta }\rightarrow \bigoplus _{\beta '\geq \beta }\Sigma ^{-\beta '}C_{\beta '}\]
 The diagram \[\xymatrix@C+40pt{{\Sigma^{-\alpha}C_{\alpha}\otimes \Sigma^{-\beta}D_{\beta}}\ar[d]|-{(-1)^{\alpha}\desusp_C^{\alpha'}f_{\alpha,\alpha'-\alpha-1}\susp_C^{\alpha}\otimes 1} \ar[r]^{\desusp^{\alpha+\beta }\susp_C^{\alpha}\otimes \susp_D^{\beta} }&{\Sigma^{-(\alpha+\beta)}(C_{\alpha} \otimes D_{\beta})} \ar[d]|-{(-1)^{\alpha+\beta}\desusp^{\alpha'+\beta}f_{\alpha,\alpha'-\alpha+1}\otimes 1\susp^{\alpha+\beta}} \\ {\Sigma^{-\alpha'}C_{\alpha'}\otimes \Sigma^{-\beta}D_{\beta}}\ar[r]_{\desusp^{\alpha'+\beta}\susp_C^{\alpha'}\otimes\susp_D^{\beta}} & {\Sigma^{-(\alpha'+\beta)}(C_{\alpha'}\otimes D_{\beta})} }\]
commutes because the composite \( \desusp _{C}^{\alpha '}f_{\alpha ,\alpha '-\alpha +1}\susp _{C}^{\alpha }\otimes 1 \)
is of degree \( -1 \), and the diagram

\[\xymatrix@C+40pt{{\Sigma^{-\alpha}C_{\alpha}\otimes \Sigma^{-\beta}D_{\beta}}\ar[d]|-{(-1)^{\beta}1\otimes \desusp_D^{\beta'}g_{\beta,\beta'-\beta-1}\susp_D^{\beta}} \ar[r]^{\desusp^{\alpha+\beta }\susp_C^{\alpha}\otimes \susp_D^{\beta} }&{\Sigma^{-(\alpha+\beta)}(C_{\alpha} \otimes D_{\beta})}\ar[d]|-{(-1)^{\alpha(\beta'-\beta)}(-1)^{\alpha+\beta}\desusp^{\alpha+\beta'}1\otimes g_{\beta,\beta'-\beta-1}\susp^{\alpha+\beta}}\\ { \Sigma^{-\alpha}C_{\alpha}\otimes \Sigma^{-\beta'}D_{\beta'}}\ar[r]_{\desusp^{\alpha+\beta'}\susp_C^{\alpha}\otimes\susp_D^{\beta'}} & {\Sigma^{-(\alpha+\beta')}(C_{\alpha}\otimes D_{\beta'})} }\]
is also seen to commute, after we permute \( \susp _{C}^{\alpha }\otimes \susp _{D}^{\beta } \)
with \( 1\otimes g_{\beta ,\beta '-\beta -1} \).

he conclusion follows.
\end{proof}
\begin{defn}
A \emph{group-object} in the category \( \mathfrak{T}_{R} \) is a right mapping
sequence \( \{C_{i},f_{*,*}\} \) equipped with a morphism 
\[
\mu :\{C_{i},f_{*,*}\}\otimes \{C_{i},f_{*,*}\}\rightarrow \{C_{i},f_{*,*}\}\]

In the case where \( C_{i}=C^{i} \)(\( i \)-fold iterated tensor product),
and the morphism \( \mu  \) just sends each \( C^{i}\otimes C^{j} \) isomorphically
to \( C^{i+j} \) we will call \( \{C_{i},f_{*,*}\} \) a \emph{free group-object}
of \( \mathfrak{T}_{R} \). It is not hard to see that \( \{f_{*,*}\} \) of
a free group object is determined by \( \{f_{1,*}\} \). 
\end{defn}
The following is immediate:

\begin{prop}
The mapping telescope of a group-object of \( \mathfrak{T}_{R} \) is a DGA-algebra.
\end{prop}
Here are some examples of mapping telescopes and group-objects:

\begin{example}
et \( C \) be a chain-complex and let \( \{f_{*,*}\} \) all be zero. Then
the sequence \( \{C^{i},f_{*,*}\} \) is a free group-object whose associated
mapping telescope \( T\{C_{i}\}_{\rightarrow } \) is nothing but the \emph{tensor
algebra} of \( \Sigma ^{-1}C \).
\end{example}
We also have the following interesting variation on this example:

\begin{prop}
\label{prop:ainftycobartelescope}Let \( C \) be an \( \ainfty  \)-coalgebra
with \( C=\integers \oplus \cbar  \) (where \( C_{0}=\integers ) \) and define
\[
f_{1,i}=(-1)^{(i+2)(i+1)/2}\Delta _{i+2}:\cbar \rightarrow \cbar ^{i+2}\]
where \( \Delta _{j}:C\rightarrow C^{j} \) is the \( \ainfty  \)-structure.
Then \( \{\cbar ^{i},f_{*,*}\} \) is a free group-object.
\end{prop}
\begin{rem}
The factor \( (-1)^{(i+2)(i+1)/2} \) is equal to the composition of
\[
\underbrace{{\susp \otimes \cdots \otimes \susp }}_{i+2\, \mathrm{factors}}\]
with 
\[
\underbrace{{\desusp \otimes \cdots \otimes \desusp }}_{i+2\, \mathrm{factors}}\]

\end{rem}
\begin{proof}
We must verify that the identities in \ref{def:rightmappingsequencetele} are
satisfied. The statement that \( \{C^{i},f_{*,*}\} \) is a free group-object
implies that 
\begin{eqnarray}
f_{i,j} & = & (-1)^{(j+2)(j+1)/2}\label{eq:cobartelescopeformula} \\
 &  & \quad \cdot \sum _{\beta =1}^{i}(-1)^{(\beta -1)(j+1)}\underbrace{{1\otimes \cdots \otimes \Delta _{j+2}\otimes \cdots \otimes 1}}_{\mathrm{position}\, \beta }\nonumber \\
 & = & \sum _{\beta =1}^{i}(-1)^{(2\beta +j)(j+1)/2}\underbrace{{1\otimes \cdots \otimes \Delta _{j+2}\otimes \cdots \otimes 1}}_{\mathrm{position}\, \beta }\nonumber 
\end{eqnarray}
 by a recursive application of \ref{prop:maptelescopeproduct}. The identity
in \ref{prop:maptelescopeformula} gives us:

\begin{eqnarray*}
\lefteqn {\sum _{\alpha =2}^{j}(-1)^{\alpha }f_{i+\alpha ,j-\alpha }\circ f_{i,\alpha -1}+(-1)^{j-1}\partial _{i+j+1}\circ f_{i,j}} &  & \\
 &  & \hspace {1in}+f_{i,j}\circ \partial _{i}=f_{i+1,j-1}\circ f_{i,0}
\end{eqnarray*}

\end{proof}
If we compute \( f_{i+\alpha ,j-\alpha }\circ f_{i,\alpha -1} \) in the case
where \( i=1 \), we get 
\begin{eqnarray*}
f_{\alpha +1,j-\alpha }\circ f_{1,\alpha -1} & = & AB\sum _{\beta =1}^{\alpha +1}(-1)^{(\beta -1)(j+1-\alpha )}\cdot \\
 &  & \hspace {.5in}\underbrace{{1\otimes \cdots \otimes \Delta _{j+2-\alpha }\otimes \cdots \otimes 1}}_{\mathrm{position}\, \beta }\circ \Delta _{\alpha -1}\\
 & = & AB\sum _{\beta =1}^{\alpha +1}(-1)^{(\beta -1)(j+1-\alpha )}\Delta _{j+2-\alpha }\circ _{\beta }\Delta _{\alpha -1}\\
 & = & (-1)^{nk+n+k}AB\cdot \\
 &  & \hspace {.5in}\sum _{\lambda =1}^{n-k}(-1)^{k\lambda +k+\lambda }\Delta _{k}\circ _{n-k-\lambda +1}\Delta _{n-k-1}
\end{eqnarray*}
where 
\begin{eqnarray*}
\beta  & = & n-k-\lambda +1\\
\alpha  & = & n-k\\
j & = & n-2\\
A & = & (-1)^{(j+2-\alpha )(j+1-\alpha )/2}\\
 & = & (-1)^{k(k-1)/2}\\
B & = & (-1)^{(\alpha -1)(\alpha -2)/2}\\
 & = & (-1)^{(n-k-1)(n-k-2)/2}
\end{eqnarray*}
or
\begin{eqnarray*}
k & = & j+2-\alpha \\
\lambda  & = & \alpha -\beta +1\\
n & = & j+2
\end{eqnarray*}
This implies that
\begin{eqnarray}
\hspace {.4in}\lefteqn {(-1)^{\alpha }f_{\alpha ,j-\alpha }\circ f_{0,\alpha -1}} &  & \label{eq:fijcompdelta} \\
 &  & \hspace {.6in}=(-1)^{nk}AB\sum _{\lambda =1}^{n-k}(-1)^{k\lambda +\lambda +k}\Delta _{k}\circ _{n-k-\lambda +1}\Delta _{n-k-1}\nonumber 
\end{eqnarray}
Now we compute
\begin{eqnarray*}
AB & = & (-1)^{(k^{2}-2k+1)/2+(n^{2}+k^{2}-2nk-2n+2k+2)/2}\\
 & = & (-1)^{(2k^{2}+n^{2}-2nk-2n+2k+3)/2}\\
 & = & (-1)^{k^{2}+k-nk+1+(n^{2}-2n+1)/2}\\
 & = & (-1)^{nk+1+n(n-1)/2}
\end{eqnarray*}
since \( k^{2}\equiv k\, (\mathrm{mod}\, 2) \). Substituting this into \ref{eq:fijcompdelta}
gives:
\begin{eqnarray*}
\lefteqn {(-1)^{\alpha }f_{\alpha ,j-\alpha }\circ f_{0,\alpha -1}=} &  & \\
 &  & \hspace {.7in}-(-1)^{n(n-1)/2}\sum _{\lambda =1}^{n-k}(-1)^{k\lambda +\lambda +k}\Delta _{k}\circ _{n-k-\lambda +1}\Delta _{n-k-1}
\end{eqnarray*}
and the conclusion follows from \ref{eq:ainfinityidentity} and \ref{prop:maptelescopeformula}.

We also make the dual definition:

\begin{defn}
A \emph{cogroup-object} in the category \( \mathfrak{T}_{L} \) is a \emph{left}
mapping sequence \( \{C_{i},f_{*,*}\} \) equipped with a morphism 
\[
\delta :\{C_{i},f_{*,*}\}\rightarrow \{C_{i},f_{*,*}\}\otimes \{C_{i},f_{*,*}\}\]

In the case where \( C_{i}=C^{i} \) (\( i \)-fold iterated tensor product),
and the morphism \( \delta  \) just sends \( C^{i}\otimes C^{j} \) isomorphically
to \( C^{i+j} \) we will call \( \{C_{i},f_{*,*}\} \) a \emph{free cogroup-object}
of \( \mathfrak{T}_{L} \). It is not hard to see that \( \{f_{*,*}\} \) of
a free group object is determined by \( \{f_{1,*}\} \).
\end{defn}
The following result is the dual to \ref{prop:ainftycobartelescope}:

\begin{prop}
\label{prop:ainftybartelescope}Let \( C \) be a \( \ainfty  \)-algebra with
structure map
\[
g:\mathfrak{A}\rightarrow \zend (C)\]
 with \( C=\integers \oplus \cbar  \) and \( C_{0}=\integers  \). Then \begin{equation}B=\xymatrix{{\cdots}\ar[r]&{\cbar\otimes \cbar}\ar[r]^-{\{f_{1,*}\}}&{\cbar}\ar[r]^{0}&{\integers}}\label{eq:barmappingsequence}\end{equation} is
a free cogroup object in \( \mathfrak{T}_{L} \), where the \( \{f_{i,j}\} \)
are given by
\[
f_{1,j}=(-1)^{(j+2)(j+1)/2}g(\mathfrak{D}_{j+2}):\cbar ^{j+2}\rightarrow \cbar \]

\end{prop}
\begin{proof}
The proof is very similar to \ref{prop:ainftycobartelescope} (we take the duals
of all of the maps there).
\end{proof}

\subsection{The cobar construction\label{subsec:mapsequencecobar}}

The cobar construction is a functor from the category of DGA-coalgebras to that
of DGA-algebras first defined by Adams in \cite{Adams:1956}\index{Adams@\textsc{J. F. Adams}},
who proved that it gave the chain complex of the loop space of a pointed simply-connected
space. One of our main results will be that it is a functor 
\[
\mathscr{F}:\mathfrak{M}\rightarrow \mathfrak{H}\]
 --- in other words, the cobar construction of an m-coalgebra comes equipped
with an m-structure. This makes it possible to iterate the cobar construction,
solving a question first posed by Adams.

We will give a somewhat nonstandard definition of the cobar construction that
will be useful in constructing its m-structure. 

\begin{cor}
\label{def:mapsequencecobar} Let \( C \) be a \( \ainfty  \)-coalgebra with
structure map
\[
g:\mathfrak{A}\rightarrow \coend (C)\]
 with \( C=\integers \oplus \cbar  \), where \( C_{0}=\integers  \). Then
(by \ref{prop:ainftycobartelescope}) \begin{equation}\xymatrix{{\integers}\ar[r]^{0}&{\cbar}\ar[r]^-{\{f_{1,*}\}}&{\cbar\otimes \cbar}\ar[r]&{\cdots}}\label{eq:cobarmappingsequence}\end{equation} is
a right mapping sequence and \( \{\cbar ^{i},f_{i,j}\} \)is a free group-object
in \( \mathfrak{T}_{R} \) and we define its associated mapping telescope \( T\{\cbar ^{i},f_{i,j}\}_{\rightarrow } \)
to be the cobar construction of \( C \). This is denoted \( \cobar (C) \).
\end{cor}
\begin{rem}
Since \( \{\cbar ^{i},f_{i,j}\} \) is a free group-object, it is clear that
\( \cobar (C) \) is a free DGA-algebra.

It is not hard to see our definition of the cobar construction produces something
isomorphic to the standard cobar construction (as defined in \cite{Adams:1956}\index{Adams@\textsc{J. F. Adams}},
\index{Baues@\textsc{Hans Baues}}\cite{Baues:1980}, or \cite{Smith:1994}\index{Smith@\textsc{J. Smith}})
in which \( \cobar (C)=T(C) \) (the tensor algebra), equipped with the differential
\[
\partial _{\cobar (C)}|\Sigma ^{-1}C=\partial _{(\Sigma ^{-1}\cbar )^{\otimes }}+\sum _{i=2}^{\infty }\underbrace{{\desusp \otimes \cdots \otimes \desusp }}_{i\, \mathrm{factors}}\circ f(\mathfrak{D}_{i})\circ \susp \]
 and required to be a DGA-algebra. The isomorphism in question is
\[
(-1)^{n(n-1)/2}\underbrace{{\desusp \otimes \cdots \otimes \desusp }}_{n\, \mathrm{factors}}\circ \susp ^{n}:\Sigma ^{-n}\cbar ^{n}\rightarrow (\Sigma ^{-1}\cbar )^{n}\]
on this direct summand, where
\[
(-1)^{n(n-1)/2}\underbrace{{\desusp \otimes \cdots \otimes \desusp }}_{n\, \mathrm{factors}}=\left( \underbrace{{\susp \otimes \cdots \otimes \susp }}_{n\, \mathrm{factors}}\right) ^{-1}\]

\end{rem}
The cobar construction comes equipped with surjections to the \emph{truncated}
cobar constructions: 
\begin{equation}
\label{eq:truncatedcobar}
\cobar (C)\rightarrow \cobar _{n}(C)=F_{1,n}
\end{equation}
 (see \ref{def:rightmappingsequencetele} for the definition of \( F_{1,n} \)).
The \( \cobar _{n}(C) \) form an inverse system with \( \cobar (C) \) as inverse
limit.

\begin{claim}
\label{claim:cobarnegative}Although the cobar construction was defined for
an m-coalgebra, it \emph{can} be defined for an \emph{arbitrary} coalgebra,
\( C \), over an operad. There is only a slight problem if \( C_{1}\neq 0 \)
--- in this case, \( (\cobar C)_{0} \) will be infinitely generated. In particular,
if \( C \) finitely-generated in each dimension and is concentrated in dimensions
\( \leq 0 \), \( \cobar C \) is perfectly well-defined and finitely generated
in each (negative) dimension.
\end{claim}
\begin{defn}
The map \( \ell :C\rightarrow \cobar (C) \) that sends \( c\in C \) to \( \Sigma ^{-1}c\in \cobar (C) \)
is a twisting cochain. The corresponding twisted tensor product, \( C\otimes _{\ell }\cobar (C) \)
is well-known to be acyclic and is called the \emph{canonical twisted tensor
product} of \( \cobar (C) \).

We can expand on \ref{def:mapsequencecobar} to get a mapping sequence for the
canonical twisted tensor product:
\end{defn}
\begin{prop}
\label{prop:mappingsequencetensor}Let \( C \) be a \( \ainfty  \)-coalgebra
with structure map
\[
g:\mathfrak{A}\rightarrow \coend (C)\]
 and cobar construction \( \cobar (C) \). Then \begin{equation}V=\xymatrix{{C}\ar[r]^-{z_1}&{C\otimes \cbar}\ar[r]^-{\{f_{1,*}\}}&{C\otimes \cbar^2}\ar[r]&{\cdots}}\label{eq:cobarprodmappingsequence}\end{equation} is
a right mapping sequence, where \( C=\integers \oplus \cbar  \), \( p:C\rightarrow \cbar  \)
is the projection, and the \( \{f_{i,j}\} \) are given by 
\[
f_{0,j}=(-1)^{j(j+1)/2}1\otimes p_{j+1}\circ \Delta _{j+2}:C\rightarrow C\otimes \cbar ^{j+1}\]
and, for \( j>0 \) 
\begin{eqnarray*}
f_{i,j} & = & \sum _{\alpha =1}^{i}(-1)^{(2\alpha -j)(j+1)/2}1\otimes \underbrace{{1\otimes \cdots \otimes \Delta _{j+2}\otimes \cdots \otimes 1}}_{\mathrm{position}\, \alpha }\\
 &  & +(-1)^{(2i-j)(j+1)/2}1\otimes p_{j+1}\circ \Delta _{j+2}\otimes 1
\end{eqnarray*}
where \( p_{j+1}=\underbrace{{p\otimes \cdots \otimes p}}_{j+1\, \mathrm{times}}:C^{j+1}\rightarrow \cbar ^{j+1} \).
The associated mapping telescope, \( T\{C^{i}\}_{\rightarrow } \), is isomorphic
to \( C\otimes _{\ell }\cobar (C) \). 
\end{prop}
\begin{rem}
The composite \( p_{j+1}\circ \Delta _{j+2} \) is effectively \( \Delta _{j+2} \)
unless \( j=0 \), in which case its kernel is \( \integers  \). For \( j>0 \),
the the image of \( \Delta _{j+2} \) does not contain a factor of \( 1\in \integers \subset C \).

The mapping sequence \ref{eq:cobarprodmappingsequence} is the result of perturbing
the tensor product, \( C\otimes \{\cbar ^{i},f_{i,j}\} \), where \( \{\cbar ^{i},f_{i,j}\} \)
is the mapping sequence for \( \cobar (C) \). Because \( \{\cbar ^{i},f_{i,j}\} \)
is a group-object in \( \mathfrak{T}_{R} \) and the boundary maps defined above
coincide with those of \( \{\cbar ^{i},f_{i,j}\} \) when restricted to this
factor, it follows that there exists a morphism
\[
V\otimes \{\cbar ^{i},f_{i,j}\}\rightarrow V\]
so that \( V \) is a right-module object in \( \mathfrak{T}_{R} \).
\end{rem}
\begin{proof}
This is very similar to the proof of \ref{prop:ainftycobartelescope} and \ref{def:mapsequencecobar}
--- we have formed the mapping sequence representing the tensor product \( C\otimes \cobar (C) \)
using \ref{prop:maptelescopeproduct} (regarding \( C \) as a mapping sequence
concentrated in the \( 0^{\mathrm{th}} \) term) and added a term representing
the twisting cochain.
\end{proof}
Now we define the notion of equivalence of \( \ainfty  \)-structures:

\begin{defn}
\label{def:ainftyequivalence} Let \( \mathfrak{R} \) be an operad that possesses
an \( A_{\infty } \)-structure. Two \( A_{\infty } \) structures on \( \mathfrak{R} \)
\[
f_{1},f_{2}:\mathfrak{A}\rightarrow \mathfrak{R}\]
 will be called \emph{equivalent} if there exists a map 
\[
\Phi :\mathfrak{A}\rightarrow \mathfrak{R}\]
 of degree \( +1 \) with the following property: For every coalgebra \( C \)
with \( \mathfrak{R} \)-action 
\[
g:\mathfrak{R}\rightarrow \coend (C)\]
 the map 
\begin{equation}
\label{eq:ainfinitymorphism}
1+\sum _{n=2}^{\infty }\underbrace{{\desusp \otimes \cdots \otimes \desusp }}_{n\, \mathrm{factors}}\circ g(\Phi (\mathfrak{D}_{n}))\circ \susp :\Sigma ^{-1}C\rightarrow (\Sigma ^{-1}C)^{\otimes }
\end{equation}
 induces an isomorphism of DGA-algebras 
\[
\cobar _{1}(C)\rightarrow \cobar _{2}(C)\]
 Here, \( \cobar _{i}(C) \) are the cobar constructions formed using \( f_{i} \),
\( i=1,2 \), respectively. 
\end{defn}
\begin{prop}
\label{prop:equivprop} Under the assumptions of \ref{def:ainftyequivalence},
the \( A_{\infty } \)-structures \( f_{1} \) and \( f_{2} \) are equivalent
if and only if \begin{multline}\Phi(\mathfrak{D}_{q})\circ\partial+   \partial\circ\Phi(\mathfrak{D}_{q}) \\
=\sum_{j=2}^{q}\sum_{\lambda=0}^{q-j} 
(-1)^{j+\lambda+j\lambda}f_{2}(\mathfrak{D}_{j})
\circ_{q-\lambda-j+1}\Phi(\mathfrak{D}_{q-j+1})\\
 -\sum_{j=2}^{q}\sum_{k_{1}+\dots+k_{j}=q} 
 (-1)^{q+\sum_{1\le\alpha<\beta\le j}(k_{\alpha}+1)k_{\beta}}\\
 \cdot\Phi(\mathfrak{D}_{k_{1}})\circ_{1}\dots
 \Phi(\mathfrak{D}_{k_{j}})\circ_{j} f_{1}(\mathfrak{D}_{j}) 
 \label{eq:ainfinityequiv}\end{multline}
\end{prop}
\begin{rem}
Although we've appealed to the published proof of this result in \cite[D\'efinition 3.4]{Proute:1},
we could have regarded the equivalence of \( \ainfty  \)-structures as a \emph{homotopy
of the identity map} of the mapping sequence for the cobar construction.
\end{rem}
\begin{proof}
A lengthy computation (see \cite[D\'efinition 3.4]{Proute:1}\index{Proute@\textsc{Prout\'e, Alain}})
shows that equivalence of \( A_{\infty } \)-structures implies an equation
like~\ref{eq:ainfinityequiv} \emph{in} \( \coend (C) \). The corresponding
equation \emph{in} \( \mathfrak{R} \) will hold if there exists a coalgebra,
\( C \), over \( \mathfrak{R} \) such that the structure map 
\[
\mathfrak{R}\rightarrow \coend (C)\]
is \emph{injective}. But this follows from~\ref{prop:everyoperad}. 
\end{proof}
\begin{prop}
\label{prop:einfinityainfinity} Let \( \mathfrak{R} \) be an \( E_{\infty } \)-operad.
Then \( \mathfrak{R} \) possesses an \( A_{\infty } \)-structure, and any
two such structures are equivalent. 
\end{prop}
\begin{proof}
This follows from~\ref{prop:equivprop} and the acyclicity of an \( E_{\infty } \)-operad's
components. 
\end{proof}

\subsection{The bar construction\label{subsec:mapseqbar}}

The bar construction was defined by Eilenberg and Mac Lane in a series of papers
beginning with \cite{EM1}. They used it to compute chain-complexes of Eilenberg-Mac
Lane spaces. They showed that the \( n \)-fold iterated bar construction of
the ring \( \zpi  \) is chain-homotopy equivalent to the chain-complex of the
Eilenberg-MacLane space \( K(\pi ,n) \).

The bar construction is dual to the cobar construction in a certain sense (and
historically preceded it) and is defined over \( A_{\infty } \)-algebras.

Recall that an \( A_{\infty } \)-algebra is a algebra over \( \mathfrak{A} \)
--- see~\ref{def:ainfinitystructure} and~\ref{prop:ainfinitythings}.

\begin{defn}
\label{def:barconst} Let \( C \) be an \( A_{\infty } \)-algebra with structure
\[
f:\mathfrak{A}\rightarrow \zend (C)\]
 with \( C=\integers \oplus \cbar  \) and \( C_{0}=\integers  \) and consider
the left mapping sequence defined in \ref{def:mapsequencecobar}: \begin{equation}B=\xymatrix{{\cdots}\ar[r]&{\cbar\otimes \cbar}\ar[r]^-{\{f_{1,*}\}}&{\cbar}\ar[r]^{0}&{\integers}}\end{equation} where
is a free cogroup object in \( \mathfrak{T}_{L} \), where the \( \{f_{i,j}\} \)
are given by
\[
f_{1,j}=(-1)^{(j+2)(j+1)/2}g(\mathfrak{D}_{j+2}):\cbar ^{j+2}\rightarrow \cbar \]
We will define the \emph{bar construction, \( \barcs (C) \)}, of \( C \) to
be the left mapping telescope of \( B \).
\end{defn}
\begin{rem}
The more traditional definition defines the differential via
\[
p\circ \partial _{\barcs }=\susp \circ \partial _{C}\circ \desusp +\sum _{n=2}^{\infty }\susp \circ f(\mathfrak{D}_{n})\circ \underbrace{{\desusp \otimes \cdots \otimes \desusp }}_{n\, \mathrm{factors}}:(\Sigma A)^{\otimes }\rightarrow \Sigma C\subset \barcs (C)\]
 where \( p:(\Sigma A)^{\otimes }\rightarrow \Sigma A \) and extended to all
of \( \barcs (C) \) by requiring it to be a free DGA-coalgebra. 

The defining identity for an \( A_{\infty } \)-algebra is equivalent to the
differential on the associated bar construction being self-annihilating ---
see \cite[\S~3]{Proute:1}.
\end{rem}
The following result shows how the bar and cobar constructions are mutually
dual:

\begin{prop}
\label{prop:bardualcobar}Let \( C \) be an \( \ainfty  \)-coalgebra with
structure morphism\index{bar-cobar duality}\index{duality!the bar and cobar constructions}
\[
f:\mathfrak{A}\rightarrow \coend (C)\]
and let \( A \) be a DGA-algebra. Then the \( \ainfty  \)-coalgebra structure
on \( C \) induces an \( \ainfty  \)-algebra structure on the hyper-chain
complex \( \homz (C,A) \) (via the cup-product map in \ref{prop:ainftyalgcoalgdual})
such that there exists a natural morphism of chain complexes
\[
\barcs (\homz (C,A))\rightarrow \homz (\cobar (C),A)\]

\end{prop}
\begin{rem}
This is not an isomorphism unless \( C \) and \( A \) are finitely generated
n all dimensions, and \( A \) has no zero-divisors.
\end{rem}
\begin{proof}
We appeal to \ref{prop:ainftyalgcoalgdual} to get an \( \ainfty  \)-module
structure on \( \homz (C,A) \). The conclusion follows from \ref{prop:mappingtelesdual},
\ref{prop:ainftycobartelescope}, \ref{def:mapsequencecobar}, and \ref{def:barconst}.
\end{proof}
It is possible to define morphisms of \( A_{\infty } \)-algebras in a manner
that is entirely dual to the definition of morphisms of \( A_{\infty } \)-coalgebras.
We get:

\begin{defn}
\label{d:ainftymorphism} Given two \( A_{\infty } \)-algebras \( (A,\{\highprod{i}\}) \),
\( (A',\{m_{i}'\}) \), a \emph{morphism} from \( A \) to \( A' \) is a family
of maps \( \{f_{i}\} \), where \( f_{i} \) is of degree \( i-1 \), satisfying:
\begin{multline*}\sum_{j=1}^{q}\sum_{k_{1}+\dots+k_{j}=q} (-1)^{\sum_{1\le\alpha<\beta\le j}(k_{\alpha}+1)k_{\beta}}m'_{j}\circ(f_{k_{1}}\otimes\dots\otimes f_{k_{j}}) \\ = \sum_{j=1}^{q}\sum_{\lambda=0}^{q-j} (-1)^{q+j+\lambda+j\lambda}f_{q-j+1}\circ(1_{\lambda}\otimes m_{j}\otimes1_{q-j-\lambda}) \end{multline*}
\end{defn}
\begin{rem}
See \cite[\S~3]{Proute:1}. Morphisms of \( A_{\infty } \)-algebras induce
DGA-coalgebra morphisms of the associated bar-constructions. This can, in fact,
be taken as their definition.
\end{rem}

\section{Operad actions on the cobar construction\label{sec:cobar}}

Given an m-coalgebra \( C \), we may form the following chain-complexes:
\begin{eqnarray*}
 & \cobar (C) & \\
 & C\otimes _{\ell }\cobar (C) & 
\end{eqnarray*}

\fancyhead[RO,LE]{\thesection. THE COBAR CONSTRUCTION} \fancyhead[RE,LO]{Justin R. Smith}We
will show that both of these construct come equipped with canonical m-coalgebra
structures (that will turn out to be geometrically valid). Since our construction
is fairly complex, we begin with an overview:

Let, \( C \), be an m-coalgebra with structure morphism 
\begin{equation}
\label{eq:m-coalgdef2}
\mathfrak{R}\rightarrow \coend (C)
\end{equation}
 According to \ref{def:mapsequencecobar}, we can regard \( \cobar (C) \) as
the mapping telescope of a mapping sequence \( \{C^{i},f_{i,j}\} \), where
\( f_{i,j} \) is a linear combination of terms \( 1\otimes \cdots \otimes \Delta _{j}\otimes \cdots \otimes 1 \).
Proposition~\ref{prop:maptelescopeproduct} implies that we can make a similar
statement about \( \cobar (C)^{n} \). Consequently, we can express \( \homz (C,\cobar (C)^{n}) \)
as a mapping telescope of a mapping sequence \( \{\homz (C,C^{i}),\homz (1,f_{i,j})\} \)
where the \( f_{i,j} \) are maps of the form \( \homz (1,1\otimes \cdots \otimes \Delta _{j}\otimes \cdots \otimes 1) \).

The m-structure of \( C \) defines maps
\[
\mathfrak{R}_{n}\rightarrow \homz (C,C^{i})\]

We use these maps to define a mapping sequence of the \( \{\mathfrak{R}_{n}\} \)and
a morphism from it to the defining mapping sequence of the cobar construction.
This gives rise to a morphism of \emph{mapping telescopes:}
\[
\iota :Z_{n}(\mathfrak{R})\rightarrow \homz (C,{(\cobar (C))}^{n})\]
 that represents the \emph{first} stage of our construction.

Now, we dualize the \( \ainfty  \)-coalgebra structure of \( C \) in \( \homz (C,{(\cobar (C))}^{n}) \)
to get an \( \ainfty  \)-\emph{algebra} structure on the cochains (using the
product-operation of \( {(\cobar (C))}^{n} \) and the results of appendix \ref{app:cobarduality}),
compute a product-operation on \( Z_{n}(\mathfrak{R}) \) such that the map
\( \iota  \) respects products and higher products. We use \ref{prop:everyoperadalgebra}
to conclude that the product operation on \( Z_{n}(\mathfrak{R}) \) \emph{actually}
makes it an \( \ainfty  \)-algebra. Now we appeal to \ref{prop:bardualcobar}
to get
\[
\barcs (Z_{n}(\mathfrak{R}))\rightarrow \barcs (\homz (C,{(\cobar (C))}^{n})\rightarrow \homz (\cobar (C),{(\cobar (C))}^{n})\]
 which implies the conclusion.

We carry out a similar construction for the functor \( C\otimes _{\alpha }\cobar (C) \)
--- we derive a functor \( Y(\mathfrak{R}) \) with a canonical map (see \ref{def:that})
\[
Y_{n}(\coend (C))\rightarrow \homz (C,(C\otimes _{\lambda }\cobar (C))^{n})\]
 This \( Y \)-functor is defined very much like the \( Z \)-functor described
above except that it is \emph{not} an \( \ainfty  \)-algebra --- instead, it
is an \( \ainfty  \)-\emph{module} over \( L(*) \) (see \ref{cor:znainftyalg}).
We will then show (see \ref{cor:bigstarmaps}) that there exists a canonical
map 
\[
Y_{n}(\coend (C))\bigstar _{\rho }L(\coend (C))\rightarrow \coend (C\otimes _{\lambda }\cobar (C))\]
where \( \bigstar  \)- is just a twisted tensor product written as \( \mathrm{Fiber}\times \mathrm{Base} \)
rather than \( \mathrm{Base}\times \mathrm{Fiber} \) (this is the reason for
using this notation). We, consequently, get a structure map
\[
Y_{n}(\mathfrak{R})\bigstar _{\rho }L(\mathfrak{R})\rightarrow \coend (C\otimes _{\lambda }\cobar (C))\]
Note that the m-structures of \( \cobar C \) and \( C\otimes _{\alpha }\cobar C \)
(the canonical acyclic twisted tensor product) are \emph{directly} \emph{induced}
from the m-structure of \( C \) in the following sense:

\begin{enumerate}
\item We regard \( \cobar C \) as a direct sum of tensor products of desuspension
of copies of \( C \).
\item In phase 1, we compute a first approximation
\begin{equation}
\label{eq:structuremapfirstapprox1}
f:Z\rightarrow \homz (C,\cobar C^{n})
\end{equation}
 to the structure map of \( \cobar C \) by taking direct sums of desuspensions
of the structure map of \( C \). For instance, on the direct summand
\[
\homz (C,(\Sigma ^{-1}C)^{2}\otimes (\Sigma ^{-1}C)^{3})\subset \homz (C,\cobar C^{2})\]
we note that there are a total of 5 copies of \( C \) and desuspend the structure
map
\[
g_{5}:\mathfrak{R}_{5}\rightarrow \homz (C,C^{5})\]
(for the unsubtle reason that it maps to a tensor product of \emph{5} copies
of \( C \)).

Of course, this wreaks havoc with the group-action of \( S_{5} \). It turns
out that this procedure makes sense only if \( Z \) in \ref{eq:structuremapfirstapprox1}
is a direct sum of desuspensions of components of \( \mathfrak{R} \) --- equipped
with a \emph{redefined} symmetric group-action. In the example above, we would
use \( \Sigma ^{-5}\mathfrak{R}_{5} \) and \emph{index} this as \( (2,3) \)
(expressing that it really maps into \( \homz (C,(\Sigma ^{-1}C)^{2}\otimes (\Sigma ^{-1}C)^{3}) \).
We drop the (now meaningless) action of \( S_{5} \) on \( \Sigma ^{-5}\mathfrak{R}_{5} \)
and replace it by an action of \( S_{2} \) that carries \( \Sigma ^{-5}\mathfrak{R}_{5} \)
indexed as \( (2,3) \) into \emph{another copy} of \( \Sigma ^{-5}\mathfrak{R}_{5} \)
indexed as \( (3,2) \).

We \emph{finally get} \( Z_{n}(\mathfrak{R}) \), which, for every \( k \),
contains \( \beta (k,n) \) copies of \( \Sigma ^{-k}\mathfrak{R}_{k} \). Here,
\( \beta (k,n) \) is the number of \emph{partitions} of \( k \) into \( n \)
parts. This is equipped with an \( S_{n} \)-action.

\item We're not done however --- even with our first approximation. The boundary map
of \( \cobar C \) is different from that of \( T(\Sigma ^{-1}C) \) (the tensor
algebra). We solve this by perturbing the natural differential of \( Z_{n}(\mathfrak{R}) \)
(i.e., the one inherited from \( \mathfrak{R} \)) via \emph{compositions} (using
the operad structure of \( \mathfrak{R} \)) involving the \emph{coproduct}
of \( C \).
\item In phase 2, we use the duality theorem to get a map
\[
\barcs Z_{n}(\mathfrak{R})\rightarrow \homz (\cobar C,\cobar C^{n})\]
for all \( n>1 \).
\end{enumerate}
\begin{prop}
\label{prop:zdefmappingsequence}Let \begin{equation}\xymatrix{{\integers}\ar[r]^-{0}&{ \bigoplus_{|\alpha|=1}{\cbar}^\alpha}\ar[r]^-{\{f_{1,*}\}}&{\bigoplus_{|\alpha|=2}{\cbar}^\alpha}\ar[r]&{\cdots}}\label{dia:mapseq1}\end{equation}
be the \( n \)-fold tensor product of the mapping sequence that defines the
cobar construction (see \ref{eq:cobarmappingsequence} and \ref{prop:maptelescopeproduct}),
where 
\begin{enumerate}
\item \( C=\integers \oplus \cbar  \)
\item \( \alpha  \) is a length-\( n \) sequence of nonnegative integers with \( |\alpha |=\sum _{i=1}^{n}\alpha _{i} \). 
\item \( {\cbar }^{\alpha }=\bigotimes _{i=1}^{n}{\cbar }^{\alpha _{i}}={\cbar }^{|\alpha |} \)---
we distinguish terms corresponding to distinct length-\( n \) sequences, \( \alpha  \).
\end{enumerate}
Form the induced mapping sequence by applying the functor \( \homz (C,*) \)
to \ref{dia:mapseq1}). 

We may fill this mapping sequence out to a diagram \[\xymatrix{{}&{\bigoplus_{|\alpha|=i}\mathfrak{R}_\alpha}\ar[d]_{g_i}&{} \\ {\cdots}\ar[r]_-{\{\bar{f}_{i-1,*}\}}&{\homz(C,\bigoplus_{|\alpha|=i}{\cbar}^\alpha)}\ar[r]_-{\{\bar{f}_{i,*}\}}&{\cdots}}\] where: 
\begin{enumerate}
\item \( \bar{f}_{i,*}=\homz (1,f_{i,*}) \), 
\item \( \mathfrak{R}_{\alpha } \) denotes \( \mathfrak{R}_{|\alpha |} \), but we
distinguish copies of \( \mathfrak{R}_{|\alpha |} \) with distinct sequences,
\( \alpha  \) adding up to the same \( |\alpha | \). Summation in \( \bigoplus _{|\alpha |=k}\mathfrak{R}_{\alpha } \)
is over all length-\( n \) sequences whose total is \( k \).
\item all vertical maps are structure maps of \( C \) (as an m-coalgebra). The summand
\( \mathfrak{R}_{\alpha }\subset \bigoplus _{|\alpha |=k}\mathfrak{R}_{\alpha } \),
with \( \alpha =\{\alpha _{1},\dots ,\alpha _{n}\} \) maps to the summand \( {\cbar }^{\alpha } \)
in the \( n \)-fold tensor product of the mapping sequence for \( \cobar (C) \).
\end{enumerate}
Then, there exist maps (represented by dotted arrows) that make every square
in the diagram \[\xymatrix{{\cdots}\ar@{.>}[r]^{z_{i-1,*}}&{\bigoplus_{|\alpha|=i}\mathfrak{R}_\alpha}\ar[d]_{g_i}\ar@{.>}[r]^{z_{i,*}}&{\cdots} \\ {\cdots}\ar[r]_-{\{\bar{f}_{i-1,*}\}}&{\homz(C,\bigoplus_{|\alpha|=i}{\cbar}^\alpha)}\ar[r]_-{\{\bar{f}_{i,*}\}}&{\cdots}}\label{dia:cobarzdef1}\]
commute.
\end{prop}
\begin{proof}
The only thing to be proved is that the maps \( \bar{f}_{i,j} \) pull back
to the upper row as maps of the \( \{\mathfrak{R}_{i,j}\} \). But this is an
immediate consequence of: 
\end{proof}
\begin{enumerate}
\item the definition of the \( \bar{f}_{i,j} \) as suitable compositions (see \ref{eq:cobarmappingsequence}
and \ref{prop:maptelescopeproduct}),
\item the definition of the action of \( \mathfrak{R} \) on \( C \).

In fact, it is not hard to write down a formula for the \( z_{i,j} \):
\[
z_{i,j}=\phi _{i,j,1}+\phi _{i,j,2}+\cdots +\phi _{i,j,n}\]
with \( n \) terms, where
\[
\phi _{i,j,k}=(-1)^{(k-1)n}\sum _{\alpha =1}^{i}(-1)^{\alpha +1}\Delta _{j+2}\circ _{\alpha +(k-1)n}*\]
The expression for \( \phi _{i,j,1} \) is clearly a direct translation of the
equation for \( f_{i,j} \) in \ref{def:mapsequencecobar} into compositions
in \( \{\mathfrak{R}_{n}\} \). The term \( \phi _{i,j,k} \) represents a translation
of \( \underbrace{{1\otimes \cdots \otimes f_{i,j}\otimes \cdots \otimes 1}}_{k^{\mathrm{th}}\, \mathrm{term}} \). 

\end{enumerate}
We have the parallel result (with a virtually identical proof):

\begin{prop}
\label{prop:ydefmappingsequence}Let {\smaller\begin{equation}\xymatrix{{\bigoplus_{|\alpha|=0}{\cbar}^{\alpha'}}\ar[r]^-{0}&{ \bigoplus_{|\alpha|=1}{\cbar}^{\alpha'}}\ar[r]^-{\{v_{1,*}\}}&{ \bigoplus_{|\alpha|=2}{\cbar}^{\alpha'}}\ar[r]&{\cdots}}\label{dia:mapsequencetensor2}\end{equation}}
be the \( n \)-fold tensor product of the mapping sequence that defines \( C\otimes _{\ell }\cobar (C) \)
(see \ref{prop:mappingsequencetensor} and \ref{prop:maptelescopeproduct}),
where
\begin{enumerate}
\item \( C=\integers \oplus \cbar  \)
\item \( \alpha  \) is a length-\( n \) sequence of nonnegative integers with \( |\alpha |=\sum _{i=1}^{n}\alpha _{i} \). 
\item given \( \alpha =\{\alpha _{1},\dots ,\alpha _{n}\} \), \( \alpha '=\{\alpha _{1}+1,\dots ,\alpha _{n}+1\} \)
\item \( {\cbar }^{\alpha '}=\bigotimes _{i=1}^{n}\cbar \otimes {\cbar }^{\alpha _{i}}={\cbar }^{n}\otimes {\cbar }^{|\alpha |} \)---
we distinguish terms corresponding to distinct length-\( n \) sequences, \( \alpha  \).
\end{enumerate}
Form the induced mapping sequence by applying the functor \( \homz (C,*) \)
to \ref{dia:mapsequencetensor2}). 

We may fill this mapping sequence out to a diagram \[\xymatrix{{}&{\bigoplus_{|\alpha|=i}\mathfrak{R}_{\alpha'}}\ar[d]_{g_i}&{} \\ {\cdots}\ar[r]_-{\{\bar{v}_{i-1,*}\}}&{\homz(C,C^{n+1})}\ar[r]_-{\{\bar{v}_{i,*}\}}&{\cdots}}\] where
\( \bar{v}_{i,*}=\homz (1,v_{i,*}) \) and all vertical maps are structure maps
of \( C \) (as an m-coalgebra). Then, there exist maps (represented by dotted
arrows) that make every square in the diagram \[\xymatrix{{\cdots}\ar@{.>}[r]^-{\{y_{i-1,*}\}}&{\bigoplus_{|\alpha|=i}\mathfrak{R}_{\alpha'}}\ar[d]_{g_1}\ar@{.>}[r]^-{\{y_{i,*}\}}&{\cdots} \\ {\cdots}\ar[r]_-{\{\bar{v}_{i-1,*}\}}&{\homz(C,C^{n+1})}\ar[r]_-{\{\bar{v}_{i,*}\}}&{\cdots}} \label{dia:tensproddef2}\]
commute.
\end{prop}
\begin{proof}
In this case, we translate the formula in \ref{prop:mappingsequencetensor}
into compositions in \( \{\mathfrak{R}_{n}\} \):

\[
y_{i,j}=\phi _{i,j,1}+\phi _{i,j,2}+\cdots +\phi _{i,j,n}\]
with \( n \) terms, where
\[
\phi _{i,0,k}=t_{1,k}\otimes 1+(-1)^{n(k-1)}\sum _{\alpha =1}^{i}(-1)^{\alpha }\Delta _{2}\circ _{n(k-1)+\alpha +1}*\]
for \( j=0 \) and and, for \( j>0 \) 
\begin{eqnarray*}
\phi _{i,j,k} & = & (-1)^{n(k-1)}\sum _{\alpha =1}^{i}(-1)^{\alpha }\Delta _{j+2}\circ _{n(k-1)+\alpha +1}*\\
 &  & +t_{j+1,k}\circ \Delta _{j+1}\circ _{n(k-1)+1}
\end{eqnarray*}
where \( t_{*,k}:\mathfrak{R}_{mn}\rightarrow \mathfrak{R}_{mn+1} \) represents
the inclusion induced by raising all indices

The expression for \( \phi _{i,j,1} \) is clearly a direct translation of the
equation for \( f_{i,j} \) in \ref{def:mapsequencecobar} into compositions
in \( \{\mathfrak{R}_{n}\} \). The term \( \phi _{i,j,k} \) represents a translation
of \( \underbrace{{1\otimes \cdots \otimes f_{i,j}\otimes \cdots \otimes 1}}_{k^{\mathrm{th}}\, \mathrm{term}} \). 
\end{proof}
Although every square commutes in diagrams \ref{dia:cobarzdef1} and \ref{dia:tensproddef2},
we cannot \emph{a priori} conclude that the upper rows of these diagrams are
mapping sequences. The following result implies this:

\begin{prop}
The upper rows of \ref{dia:cobarzdef1} and \ref{dia:tensproddef2} are mapping
sequences in the sense of \ref{dia:mapseq1}.
\end{prop}
\begin{proof}
We know the following two things about the upper rows of \ref{dia:cobarzdef1}
and \ref{dia:tensproddef2}:
\begin{enumerate}
\item They are independent of \( C \)
\item The composites of the maps \( \{z_{i,j}\} \) and \( \{y_{i,j}\} \), respectively
with the corresponding vertical maps satisfy the defining identities (see \ref{def:rightmappingsequencetele})
of a mapping sequence (because every square of the diagrams commutes and because
the lower rows are mapping sequences).
\end{enumerate}
It suffices, therefore, to find a \( C \) for which the vertical maps in diagrams
\ref{dia:cobarzdef1} and \ref{dia:tensproddef2} are injective. But the existence
of such a \( C \) follows from \ref{prop:everyoperad} and Claim~\ref{claim:cobarnegative}
(which implies that we can use the coalgebra given by \ref{prop:everyoperad}
even though it is concentrated in dimensions \( \leq 0 \)). This completes
the proof.

\end{proof}
\begin{cor}
\label{cor:mapttelecobartensor}Given an m-coalgebra \( C \) over an \( \ainfty  \)-operad,
\( \mathfrak{R} \), there exist chain-complex morphisms 
\begin{itemize}
\item \( \iota _{n}:Z_{n}'(\mathfrak{R})=T_{\rightarrow }\{\bigoplus _{|\alpha |=i}\mathfrak{R}_{\alpha },z_{i,j}\}\rightarrow \homz (C,\cobar (C)^{n}) \)
\item \( \kappa _{n}:Y_{n}'(\mathfrak{R})=T_{\rightarrow }\{\bigoplus _{|\alpha |=i}\mathfrak{R}_{\alpha '},y_{i,j}\}\rightarrow \homz (C,(C\otimes _{\ell }\cobar (C))^{n}) \)
\end{itemize}
\end{cor}
\begin{prop}
\label{def:that} The boundary maps of \( Z_{n}'(\mathfrak{R}) \) and \( Y_{n}'(\mathfrak{R}) \)
are given, respectively, by: \begin{equation}\partial_{Z}|\Sigma^{-|\alpha|}\mathfrak{R}_{\alpha}= \sum_{i=1}^{|\alpha|}\sum_{j=1}^{\infty}(-1)^{i+|\alpha|j+ij}\desusp^{|\alpha|+j-1}\circ(f(\mathfrak{D}_{j}) \circ_{i} *)\circ\susp^{|\alpha|}\label{eq:thatdef}\end{equation}
and \begin{equation}\partial_{Y}|\Sigma^{-|\alpha|}\mathfrak{R}_{\alpha'}= \sum_{i=1}^{|\alpha|}\sum_{j=1}^{\infty}(-1)^{i+|\alpha|j+ij}\desusp^{|\alpha|+j-1}(f(\mathfrak{D}_{j}) \circ_{i+n} *)\circ\susp^{|\alpha|}\label{eq:thatdefy}\end{equation}
Given a sequence \( \alpha =\{\alpha _{1},\dots ,\alpha _{n}\} \), we will
say that \( i \) \emph{lies inside the} \( \alpha (i)^{\mathrm{th}} \)-term
of \( \alpha  \) if \( i<\sum _{s=1}^{\alpha (i)}\alpha _{s} \) and \( \alpha (i) \)
is the smallest value with this property. Then the image of 
\[
\desusp ^{|\alpha |+j-1}\circ (f(\mathfrak{D}_{j})\circ _{i}*)\circ \susp ^{|\alpha |}(\Sigma ^{-|\alpha |}\mathfrak{R}_{\alpha })\]
 is defined to lie in 
\[
\Sigma ^{-|\alpha |+j-1}\mathfrak{R}_{\{\alpha _{1},\dots ,\alpha _{\alpha (i)}+j-1,\dots ,\alpha _{n}\}}\]

\end{prop}
\begin{rem}
In the boundary map of \( \homz (C,(\Sigma ^{-1}C)^{N})\subset \homz (C,{\cobar _{m}(C)}^{n}) \),
all factors of \( (\Sigma ^{-1}C)^{N} \) are on an equal footing. We only make
distinctions between the various factors of \( \cobar (C) \) that they represent
when we decide where the image of the boundary map lies.

Note that when \( |\alpha |=0 \), the boundary map of \( Y \) is \emph{identical}
to that of \( \mathfrak{R} \).
\end{rem}
We define:

\begin{defn}
\label{def:zn} Let \( \mathfrak{R}=\{\mathfrak{R}_{n}\} \) be an \( \ainfty  \)-operad
and let \( Z_{n}'(\mathfrak{R}) \) and \( Y'(\mathfrak{R}) \) be as defined
in \ref{cor:mapttelecobartensor}. Define \( Z_{n}(\mathfrak{R}) \) and \( Y_{n}(\mathfrak{R}) \),
respectively, to be the results of truncating \( Z_{n}'(\mathfrak{R}) \) and
\( Y_{n}'(\mathfrak{R}) \) in dimension \( -1 \). In addition, we equip \( Z_{n}(\mathfrak{R}) \)
and \( Y_{n}(\mathfrak{R}) \) with a \( \zs{{n}} \)-action defined as follows:

if \( a\in \Sigma ^{-|\alpha |}\mathfrak{R}_{\alpha } \), \( \sigma \in S_{n} \)
then \( \sigma \cdot a=\mathrm{Parity}(\tlist {\alpha }{n}(\sigma ))\tlist {\alpha }{n}(\sigma )\cdot a \)
where the product on the right-hand side is taken in \( \Sigma ^{-|\alpha |}\mathfrak{R}_{|\alpha |} \)
and the target is regarded as an element of \( \Sigma ^{-|\alpha |}\mathfrak{R}_{\sigma ^{-1}\{\alpha _{1},\dots ,\alpha _{n}\}} \).
Here \( \sigma ^{-1}\{\alpha _{1},\dots ,\alpha _{n}\}=\{\alpha _{\sigma (1)},\dots ,\alpha _{\sigma (n)}\} \)
is the result of permuting the elements of \( \alpha  \) via \( \sigma ^{-1} \).
Note that we are defining the action of \( S_{n} \) to permute the summands
of \( Z_{n}(\mathfrak{R}) \), as well as twisting it by the parity of the permutations.
We define a corresponding action on \( Y_{n}(\mathfrak{R}) \) by replacing
the sequence \( \{\alpha _{1},\dots ,\alpha _{n}\} \) by \( \alpha '=\{\alpha _{1}+1,\dots ,\alpha _{n}+1\} \).
\end{defn}
\begin{prop}
\label{cor:uniqueznm} Let \( \mathfrak{R} \) be an operad with two \( A_{\infty } \)-structures
\[
f_{i}:\mathfrak{A}\rightarrow \mathfrak{R}\]
and with \( n \)-components that are \( \zs{{n}} \)-free for all \( n>1 \).
Then any equivalence between \( f_{1} \) and \( f_{2} \) induces canonical
isomorphisms of chain complexes 
\[
Z_{n}(\mathfrak{R},\partial _{f_{1}})\rightarrow Z_{n}(\mathfrak{R},\partial _{f_{2}})\]
and 
\[
Y_{n}(\mathfrak{R},\partial _{f_{1}})\rightarrow Y_{n}(\mathfrak{R},\partial _{f_{2}})\]

\end{prop}
\begin{proof}
We begin with an observation: Any \( \ainfty  \)-structure on \( \mathfrak{R} \)
induces one on \( \coend (C) \). It follows that we can form \( Z_{n}(\coend (C)) \)
and \( Y_{n}(\coend (C)) \) and that there exist canonical morphisms
\begin{equation}
\label{eq:morphieqa}
u:Z_{n}(\mathfrak{R})\rightarrow Z_{n}(\coend (C))
\end{equation}
and 
\begin{equation}
\label{eq:morphieqb}
v:Y_{n}(\mathfrak{R})\rightarrow Y_{n}(\coend (C))
\end{equation}
It is also not hard to see that these morphisms will be injective if the m-structure
morphism of \( C \) is.

\emph{Claim 1:} The morphisms
\[
\iota _{n}:Z_{n}'(\coend (C))\rightarrow \homz (C,\cobar (C)^{n})\]
 and
\[
\kappa _{n}:Y_{n}'(\coend (C))\rightarrow \homz (C,(C\otimes _{\ell }\cobar (C))^{n})\]
are \emph{isomorphisms}. This follows from the fact that the maps in stage of
the corresponding mapping telescopes are isomorphisms. We conclude that equivalent
\( \ainfty  \)-structures give rise to \emph{isomorphic} \( Z_{n}(\coend (C)) \)
and \( Y_{n}(\coend (C)) \).

\emph{Claim 2:} The \emph{images} of \( u \) and \( v \) in equations \ref{eq:morphieqa}
and \ref{eq:morphieqb} are \emph{independent} of the \( \ainfty  \)-structure
since they are induced by structure maps of \( C \) (which only depend on the
m-structure of \( C \)).
\end{proof}
Now we appeal to \ref{prop:everyoperad} to conclude the existence of a coalgebra,
\( C \), with injective structure morphism. Claim \ref{claim:cobarnegative}
implies that this coalgebra does not have to be an m-coalgebra. In this case,
both \( Z_{n}(\mathfrak{R},\partial _{f_{1}}) \) and \( Z_{n}(\mathfrak{R},\partial _{f_{2}}) \)
will have the same image in \( Z_{n}(\coend (C)) \). The corresponding statement
holds for \( Y_{n}(\mathfrak{R},\partial _{f_{1}}) \) and \( Y_{n}(\mathfrak{R},\partial _{f_{2}}) \).

\begin{cor}
\label{cor:einfinityzunique} Let \( \mathfrak{R} \) be an \( E_{\infty } \)-operad.
Then the chain complexes 
\[
Z_{n}(\mathfrak{R}),\mathfrak{T}(f))\]
 
\[
Y_{n}(\mathfrak{R},\partial _{f})\]
are uniquely determined up to an isomorphism. 
\end{cor}
\begin{proof}
See~\ref{prop:einfinityainfinity} and~\ref{cor:uniqueznm}. 
\end{proof}
\begin{cor}
\label{cor:znainftyalg}If \( \mathfrak{R} \) is an operad with an \( \ainfty  \)-structure
with structure map 
\[
f:\mathfrak{A}\rightarrow \mathfrak{R}\]
and \( \mathfrak{R}_{i} \) is \( \zs{{i}} \)-free for all \( i \), then \( (Z_{n}(\mathfrak{R}),\{t_{n}(f(\Delta _{k}))\}) \)
constitutes an \( \ainfty  \)-algebra for all \( n>1 \) where \( t_{n}:\mathfrak{R}\rightarrow \zend (Z_{n}(\mathfrak{R}) \)
is defined in \ref{def:lowertmaps} in appendix \ref{app:cobarduality}. In
addition, there exists a map
\[
Y_{n}(\mathfrak{R},\partial _{f})\otimes Z_{n}(\mathfrak{R})\otimes \cdots \otimes Z_{n}(\mathfrak{R})\rightarrow Y_{n}(\mathfrak{R},\partial _{f})\]
for all \( n>1 \) making \( Y_{n}(\mathfrak{R},\partial _{f}) \) a right \( \ainfty  \)-module
over \( Z(\mathfrak{R}) \).
\end{cor}
\begin{rem}
Although appendix~\ref{app:cobarduality} gives the details of this \( \ainfty  \)-structure,
we will summarize some salient features:
\begin{itemize}
\item If \( C \) is an m-coalgebra over \( \mathfrak{R} \), its structure map sends
this \( \ainfty  \)-product to the natural \( \ainfty  \)-algebra structure
on \( \homz (C,(\cobar C)^{n}) \) defined by duality (see \ref{lem:operadduality}).
\item On summands of \( Z_{n}(\mathfrak{R}) \) (as defined in \ref{def:zn}) it is
a map 
\[
t_{n}(f(\Delta _{k}):Z_{\beta _{1}}(\mathfrak{R})\otimes \cdots Z_{\beta _{k}}(\mathfrak{R})\rightarrow Z_{\alpha }(\mathfrak{R})\]
 where:

\begin{enumerate}
\item each \( \beta _{i} \) is a length-\( n \) sequence of nonnegative integers
\item \( \alpha  \) is a length-\( n \) sequence of integers equal to the \emph{elementwise}
sum of the \( \beta _{i} \)
\end{enumerate}
\item On summands of \( Z_{n}(\mathfrak{R}) \) (as defined in \ref{def:zn}) it is
a map 
\[
t_{n}(f(\Delta _{k}):Y_{\beta _{1}}(\mathfrak{R})\otimes \cdots Z_{\beta _{k}}(\mathfrak{R})\rightarrow Y_{\alpha }(\mathfrak{R})\]
 where:

\begin{enumerate}
\item each \( \beta _{i} \) is a length-\( n \) sequence of nonnegative integers
\item \( \alpha  \) is a length-\( n \) sequence of integers equal to the \emph{elementwise}
sum of the \( \beta _{i} \)
\end{enumerate}
\end{itemize}
\end{rem}
\begin{proof}
Our first claim is proved in appendix \ref{app:cobarduality}. Our second claim
follows by a variation on the construction in \ref{def:lowertmaps}: we replace
each term of \( \beta _{1} \) by \( \mathrm{itself}+1 \) in all the subscripts
in diagram \ref{dia:tnmapdef} in appendix \ref{app:cobarduality}. The conclusion
follows from the definition of the \( Y_{n}(\mathfrak{R},\partial _{f}) \)
--- particularly the fact that they were defined exactly like the \( Z_{n}(\mathfrak{R},\partial _{f}) \)
with a shift in indices.
\end{proof}
\begin{cor}
\label{cor:bigstarmaps}Let \( C \) be a coalgebra over the \( E_{\infty } \)-operad
\( \mathfrak{R} \). Then there exist morphisms of \( \zs{{n}} \)-chain complexes
\begin{equation}
\label{eq:hstarmap}
h^{*}:\barcs (Z_{n}(\mathfrak{R}))\rightarrow \homz (\cobar (C),\cobar (C)^{n})
\end{equation}
and
\begin{equation}
\label{eq:bigstarmap}
g:Y_{n}(\mathfrak{R},\partial _{f})\bigstar _{\rho _{n}}\barcs (Z_{n}(\mathfrak{R}))\rightarrow \homz \left( C\otimes _{\ell }\cobar (C),\left( C\otimes _{\ell }\cobar (C)\right) ^{n}\right) 
\end{equation}
where \( \rho _{n}:\barcs (Z_{n}(\mathfrak{R}))\rightarrow Z_{n}(\mathfrak{R}) \)
is the canonical twisting cochain that makes 
\[
Z_{n}(\mathfrak{R},\partial _{f})\bigstar _{\rho _{n}}\barcs (Z_{n}(\mathfrak{R}))\]
acyclic and \( Y_{n}(\mathfrak{R},\partial _{f})\bigstar _{\rho _{n}}\barcs (Z_{n}(\mathfrak{R})) \)
is the corresponding twisted tensor product.
\end{cor}
\begin{rem}
We have written \( Y_{n}(\mathfrak{R},\partial _{f})\bigstar _{\rho _{n}}\barcs (Z_{n}(\mathfrak{R})) \)
for the twisted tensor product in this result because it is of the form \( \mathrm{fiber}\times \mathrm{base} \)
rather than the usual \( \mathrm{base}\times \mathrm{fiber} \).

If \( F \) is a DGA-algebra, the transposition map defines a \emph{canonical
isomorphism} of chain-complexes \( T:F\bigstar _{{x,f}}C\cong C'\otimes _{x}F \),
where \( C' \) is the DGA coalgebra with the same underlying chain-complex
as \( C \), but whose coproduct \( \Delta _{C'}=T\circ \Delta _{C} \) and
\( \Delta _{C} \) is the coproduct of \( C \). Here is the twisted tensor
product formed with respect to the coalgebra structure of \( C' \) and the
\( A_{\infty } \)-algebra structure of \( F \). 
\end{rem}
\begin{proof}
The first claim follows from \ref{prop:bardualcobar}. To prove the second,
let
\[
\bar{h}:\barcs (Z_{n}(\mathfrak{R}))\otimes \cobar (C)\rightarrow \cobar (C)^{n}\]
 be the adjoint of \( h^{*} \) in \ref{eq:hstarmap} and let 
\[
\bar{\kappa }_{n}:Y_{n}(\mathfrak{R})\otimes C\rightarrow (C\otimes _{\ell }\cobar (C))^{n}\]
be the adjoint of \( \kappa _{n} \) defined in \ref{cor:mapttelecobartensor}.
We form the tensor product
\[
\bar{\kappa }_{n}\otimes \bar{h}:Y_{n}(\mathfrak{R})\otimes C\otimes \barcs (Z_{n}(\mathfrak{R}))\otimes \cobar (C)\rightarrow (C\otimes _{\ell }\cobar (C))^{n}\otimes \cobar (C)^{n}\]
Now we compose this on the \emph{right} with the map \( \mu _{n}:(C\otimes _{\ell }\cobar (C))^{n}\otimes \cobar (C)^{n}\rightarrow (C\otimes _{\ell }\cobar (C))^{n} \)
defining the action of \( \cobar (C)^{n} \) on \( (C\otimes _{\ell }\cobar (C))^{n} \)
and with \( 1\otimes T\otimes 1 \) on the left to get a map
\[
g^{*}:Y_{n}(\mathfrak{R})\otimes \barcs (Z_{n}(\mathfrak{R}))\otimes C\otimes \cobar (C)\rightarrow (C\otimes _{\ell }\cobar (C))^{n}\]
The tensor products on the left are all untwisted, but a simple computation
shows that the effects of replacing \( Y_{n}(\mathfrak{R})\otimes \barcs (Z_{n}(\mathfrak{R})) \)
by \( Y_{n}(\mathfrak{R})\bigstar _{\rho _{n}}\barcs (Z_{n}(\mathfrak{R})) \)
and \( C\otimes \cobar (C) \) by \( C\otimes _{\ell }\cobar (C) \), respectively,
\emph{exactly cancel out} (when mapped by \( v^{*} \)), giving a map
\[
g^{*}:Y_{n}(\mathfrak{R})\bigstar _{\rho }\barcs (Z_{n}(\mathfrak{R}))\otimes C\otimes _{\ell }\cobar (C)\rightarrow (C\otimes _{\ell }\cobar (C))^{n}\]
that is the \emph{adjoint} of the map, \( v \), that we want. 
\end{proof}
Our final observation is that we can pull back the composition-operations of
\( \homz (\cobar (C),\cobar (C)^{n}) \) to get an operad structure on \( \barcs (Z_{n}(\coend (C))) \)
that makes the map \( h^{*} \), defined above, into a morphism of operads.
For a detailed description of these composition operations, see \S~\ref{app:operadcomps}.

\begin{prop}
Let \( \mathfrak{R} \) be an \( E_{\infty } \)-operad. Then the composition
operations defined in \ref{def:cobaroperadcomps} make \( L'(\mathfrak{R})=\{\barcs (Z_{n}(\mathfrak{R}))\} \)
into an operad-coalgebra.
\end{prop}
\begin{proof}
That \( L'(\mathfrak{R}) \) is an operad follows from the fact that \( L'(\coend (C)) \)
is an operad for any comodule over an \( E_{\infty } \)-operad (by the discussion
preceding \ref{def:cobaroperadcomps}) and \ref{prop:everyoperad} and \ref{claim:cobarnegative}.
That these composition operations commute with the coproduct of \( \barcs (Z_{n}(\mathfrak{R})) \)
follows from: 
\begin{enumerate}
\item the definition of the coproduct
\[
\Delta ([r_{1}|\dots |r_{k}]=\sum _{i=0}^{k}[r_{1}|\dots |r_{i}]\otimes [r_{i+1}|\dots |r_{k}]\]

\item the fact that \( [s_{1}|\dots |s_{t}]\circ _{i}[r_{1}|\dots |r_{k}]=[r_{1}'|\dots |r_{k}'] \)
(by \ref{def:cobaroperadcomps}),
\item for each term of \( \Delta ([r_{1}|\dots |r_{k}] \) only a single term of \( \Delta [s_{1}|\dots |s_{t}] \)
will have a nonzero composition with it. This follows from the fact that \( [s_{1}|\dots |s_{t}]\circ _{i}[r_{1}|\dots |r_{k}]=0 \)
unless \( t=\sum _{z=1}^{k}\beta _{i,z} \) (in \ref{def:cobaroperadcomps}).
\end{enumerate}
\end{proof}
We will conclude this section by deriving a variation on \( L(\mathfrak{R}) \)
that is \( E_{\infty } \)-whenever \( \mathfrak{R} \) is an \( E_{\infty } \)-operad.
We begin by observing 

\begin{prop}
If \( \mathfrak{R} \) is an \( E_{\infty } \)-operad, then there exist twisting
cochains
\[
c_{n}:\rs{{n}}\rightarrow Z_{n}(\mathfrak{R})\]
for all \( n>0 \). These are formed with respect to the \( \ainfty  \)-structure
defined in \ref{def:lowertmaps}.
\end{prop}
\begin{proof}
This is a straightforward consequence of the fact that the \( Z_{n}(\mathfrak{R}) \)
are acyclic in above dimension \( -1 \) if \( \mathfrak{R} \) is \( E_{\infty } \).
We can use the identity that a twisting cochain must satisfy --- \ref{d:ainftyalgtwistedtensor}
--- to compute \( c_{n} \) inductively. We set \( c_{n}=\bigoplus z_{\alpha } \),
where \( \alpha =\{\alpha _{1},\dots ,\alpha _{n}\} \) is a sequence of nonnegative
integers and the direct sum is taken over all such sequences.

We begin the induction by requiring that \( z_{\alpha }:\rs{{1}}\rightarrow \Sigma ^{-1}\mathfrak{R}_{1} \)
sends \( 1\in \rs{{1}} \) to \( \Sigma ^{-1}1 \) whenever \( \alpha  \) is
a sequence with a single nonzero element equal to \( 1 \). The identity for
a twisting cochain allows us to extend this to canonical basis elements of \( \rs{{n}} \)
in higher dimensions and we extend it to all of \( \rs{{n}} \) by setting
\[
z_{\alpha }(\sigma \cdot a)=\mathrm{parity}(\sigma )\tlist {\alpha }{n}(\sigma )z_{\sigma ^{-1}(\alpha )}(a)\]

\end{proof}
In \cite{Smith:1994}, I developed software for computing the \( \{z_{\alpha }\} \).
In the present paper, our main interest in them will be:

\begin{cor}
\label{cor:mainresult}If \( \mathfrak{R} \) is an \( E_{\infty } \)-operad,
then for all \( n>0 \) there exist coalgebra morphisms 
\[
\barcs (c_{n}):\rs{{n}}\rightarrow \barcs (Z_{n}(\mathfrak{R}))\]
inducing morphisms
\[
H_{0}(\rs{{n}})=\integers \rightarrow H_{0}(\barcs (Z_{n}(\mathfrak{R})))\]
The pullback of these morphisms over the projection \( \barcs (Z_{n}(\mathfrak{R}))\rightarrow H_{0}(\barcs (Z_{n}(\mathfrak{R}))) \)
is an \( E_{\infty } \)-operad coalgebra, \( L_{n}(\mathfrak{R}) \), \index{L@$L(*)$ definition}\index{definition!$L(*)$ functor}with
the property that \( \cobar (C) \) is an m-coalgebra over \( L(\mathfrak{R}) \)
if \( C \) is an m-coalgebra over \( \mathfrak{R} \).

The twisting cochain \( \rho :\barcs (Z_{n}(\mathfrak{R}))\rightarrow Y_{n}(\mathfrak{R}) \)
restricts to a twisting cochain \( \rho _{n}:L_{n}(\mathfrak{R})\rightarrow Y_{n}(\mathfrak{R}) \)
such that the twisted tensor product \( Y_{n}(\mathfrak{R})\bigstar _{\rho _{n}}L_{n}(\mathfrak{R}) \)
is acyclic in positive dimensions. 

The composite inclusion
\[
\mathfrak{R}_{n}\subset Y_{n}(\mathfrak{R})\subset Y_{n}(\mathfrak{R})\bigstar _{\rho _{n}}L_{n}(\mathfrak{R})\]
induces a map
\[
H_{0}(\rs{{n}})=\integers \rightarrow H_{0}(Y_{n}(\mathfrak{R})\bigstar _{\rho _{n}}L_{n}(\mathfrak{R}))\]
The pullback of these maps over the projection 
\[
Y_{n}(\mathfrak{R})\bigstar _{\rho _{n}}L_{n}(\mathfrak{R})\rightarrow H_{0}(Y_{n}(\mathfrak{R})\bigstar _{\rho _{n}}L_{n}(\mathfrak{R}))\]
 is an \( E_{\infty } \)-operad, \( J(\mathfrak{R})=\{J_{n}(\mathfrak{R})\} \)
such that \( C\otimes _{\lambda }\cobar C \) is a coalgebra over \index{J@$J(*)$!definition}\index{definition!$J(*)$ functor}it.
\end{cor}
\begin{rem}
We view \( H_{0}(\rs{{n}})=\integers  \) as an operad concentrated in dimension
\( 0 \).

The maps \( \barcs (c_{n}):\rs{{n}}\rightarrow \barcs (Z_{n}(\mathfrak{R})) \)
are useful for performing explicit computations of cohomology operations on
the cobar construction.
\end{rem}
\begin{proof}
To conclude that the pullbacks can be taken, we invoke \ref{cor:operadpullback}.
\end{proof}
\begin{cor}
\label{cor:cobarmhopf}The functor, \( \cobar (*) \), that maps an m-coalgebra
\( C \) over an \( E_{\infty } \)-operad \( \mathfrak{R} \) to \( \cobar (C) \)
over the \( E_{\infty } \)-operad coalgebra, \( L(\mathfrak{R}) \), defines
a functor
\[
\cobar (*):\mathfrak{M}\rightarrow \mathfrak{H}\]

\end{cor}
\begin{proof}
We must show that it maps equivalent objects to equivalent objects. It is easy
to see that \( \cobar (*) \) maps split injections (of \( C \)) into multiplicative
split injections (of \( \cobar (C) \)) and a simple spectral sequence argument
shows that maps induced by homology equivalences are homology equivalences.
\end{proof}
As one might expect:

\begin{prop}
\label{prop:tensorproductprojection}Let \( C \) be an m-coalgebra over an
\( E_{\infty } \)-operad \( \mathfrak{R} \). Then there exists a morphism
of \index{morphism $J(\mathfrak{R}\rightarrow\mathfrak{R}$}operads
\[
J(\mathfrak{R})\rightarrow \mathfrak{R}\]
supporting a morphism in \( \mathfrak{M}_{0} \):
\[
C\otimes _{\ell }\cobar (C)\rightarrow C\]

\end{prop}
\begin{rem}
Note that there is \emph{always} a canonical morphism of operads 
\[
J(\mathfrak{R})\rightarrow L(\mathfrak{R})\]
defined by projection to the base. It is compatible with the inclusion of the
fiber in 
\[
\cobar C\rightarrow C\otimes _{\lambda }\cobar C\]
The map in the present result is projection to (a suboperad of) the \emph{fiber}
of the twisted tensor product, and only exists because of peculiarities in the
twisted differential of \( Y_{n}(\mathfrak{R})\bigstar _{\rho _{n}}L_{n}(\mathfrak{R}) \)
due to the action of \( Z_{n}(\mathfrak{R}) \) on \( Y_{n}(\mathfrak{R}) \).

On the other hand, because \( J(\mathfrak{R}) \) is not really a twisted tensor
product but only a suboperad of one (namely the pullback over an inclusion of
\( \integers  \) in \( H_{0}(Y_{n}(\mathfrak{R})\bigstar _{\rho _{n}}L_{n}(\mathfrak{R})) \)),
one \emph{cannot} conclude that there exists an \emph{inclusion of the fiber}
\[
\mathfrak{R}\rightarrow J(\mathfrak{R})\]
--- this fails in dimension \( 0 \).
\end{rem}
\begin{proof}
Let \( \alpha =\{\alpha _{1},\dots ,\alpha _{n}\} \) be a sequence of nonnegative
integers and let \( U_{\alpha }=\Sigma ^{-|\alpha |}\mathfrak{R}_{n+|\alpha |}\subset Y_{n}(\mathfrak{R}) \)
be a direct summand (see \ref{eq:thatdefy}). Note that \( Y_{n}(\mathfrak{R})\bigstar _{\rho _{n}}L_{n}(\mathfrak{R})\supset J(\mathfrak{R})_{n} \)
contains a summand 
\[
S=U_{\{0,\dots ,0\}}\bigstar _{\rho _{n}}L_{n}(\mathfrak{R})\]
where \( U_{\{0,\dots ,0\}}=\Sigma ^{-0}\mathfrak{R}_{n}=\mathfrak{R}_{n} \).
Proposition~\ref{prop:ycontainsr} implies that the composition operations of
\( U_{\{0,\dots ,0\}} \) coincide with those of \( \mathfrak{R} \). Since
the twisted differential of \( S \) consists of terms
\[
U_{\{\alpha _{1},\dots ,\alpha _{n}\}}\bigstar _{\rho _{n}}L_{n}(\mathfrak{R})\]
with at least one of the \( \alpha _{i}>0 \) (see \ref{cor:znainftyalg} and
\ref{def:lowertmaps}), it follows that we can project \( S \) to \( U_{\{0,\dots ,0\}}=\mathfrak{R}_{n} \)
--- despite the fact that it is in the \emph{fiber} of the twisted tensor product.
An examination of the composition-operations in \( Y_{n}(\mathfrak{R})\bigstar _{\rho _{n}}L_{n}(\mathfrak{R}) \)
shows that this projection is compatible with compositions because, restricted
to \( S\subset Y_{n}(\mathfrak{R})\bigstar _{\rho _{n}}L_{n}(\mathfrak{R}) \),
they are a \emph{perturbation} of the compositions of \( \mathfrak{R} \) ---
and the perturbations \emph{lie} in summands 
\[
U_{\{\alpha _{1},\dots ,\alpha _{n}\}}\bigstar _{\rho _{n}}L_{n}(\mathfrak{R})\]
with at least one of the \( \alpha _{i}>0 \).

This projection induces the map in the conclusion. The projection 
\[
C\otimes _{\lambda }\cobar C\rightarrow C\]
 is clearly compatible with this morphism of \( E_{\infty } \)-operads because
the action of \( L_{n}(\mathfrak{R}) \) gets killed off.
\begin{lem}
\label{lem:liftinglemma}\textbf{(Lifting lemma)} \index{Lifting lemma}Let
\( A \), \( B \), and \( C \) be m-coalgebras over the operad \( \mathfrak{R} \)
and suppose there exists maps in the category \( \mathfrak{M}_{0} \)\[\xymatrix{{A}\ar[r]^f&{B}\ar[r]^g&{C}}\] such
that \( g\circ f=0 \) in positive dimensions. Compose the structure map of
\( A \) with the projection
\[
z:J(\mathfrak{R})\rightarrow \mathfrak{R}\]
so that \( A \) becomes an m-coalgebra over \( J(\mathfrak{R}) \). Then the
map \( a\mapsto f(a)\otimes 1\in B\otimes _{\alpha \circ g}\cobar C_{0}\subset B\otimes _{\alpha \circ g}\cobar C \),
for all \( a\in A \), defines a morphism in \( \mathfrak{M}_{0} \)
\[
f':A\rightarrow B\otimes _{\alpha \circ g}\cobar C=(B\otimes _{\alpha }\cobar B)\otimes _{\cobar B}\cobar C\]

\end{lem}
\begin{rem}
The qualifier ``in positive dimensions'' was necessary for our statement that
\( g\circ f=0 \): All morphism in \( \mathfrak{M}_{0} \) induce isomorphisms
in dimension \( 0 \) (which, for any m-coalgebra is equal to \( \integers  \).
\end{rem}
First, note that the map \( f' \) is a chain-map. This follows from the fact
that the twisting cochain, \( \alpha \circ g:B\rightarrow C \) vanishes on
\( \im f\subset B \). So the main thing to be proved is that it preserves m-structures.

Let the structure maps of \( A \) and \( B\otimes _{\alpha \circ g}\cobar C \)
be
\begin{eqnarray*}
\mathfrak{a}:\mathfrak{R} & \rightarrow  & \coend (A)\\
\mathfrak{d}:J(\mathfrak{R}) & \rightarrow  & \coend (B\otimes _{\alpha \circ g}\cobar C)
\end{eqnarray*}

As in the proof of \ref{prop:tensorproductprojection}, let \( \alpha =\{\alpha _{1},\dots ,\alpha _{n}\} \)
be a sequence of nonnegative integers and let \( U_{\alpha }=\Sigma ^{-|\alpha |}\mathfrak{R}_{n+|\alpha |}\subset Y_{n}(\mathfrak{R}) \)
be a direct summand (see \ref{eq:thatdefy}). Now recall the definition of the
action of \( U_{\alpha }=\Sigma ^{-|\alpha |}\mathfrak{R}_{n+|\alpha |} \)
on \( B\otimes 1\subset B\otimes _{\alpha }\cobar B \) in \ref{prop:ydefmappingsequence}:
the summand \( U_{\{0,\dots ,0\}}=\mathfrak{R}_{n} \) acts on \( B\otimes 1 \)
in exactly the same manner as \( \mathfrak{R}_{n} \) acts on \( B \). Consequently,
the fact that \( f:A\rightarrow B \) is a morphism in \( \mathfrak{R}_{0} \)
implies that the action of \( Y_{0} \) on \( \im f' \) is compatible with
the action of \( \mathfrak{R} \) on \( A \). It follows that it will suffice
to show that the kernel of 
\[
z:J(\mathfrak{R})\rightarrow \mathfrak{R}\]
acts trivially on \( B\otimes _{\alpha \circ g}\cobar C \), at least when restricted
to \( \im f' \) --- or
\[
\mathfrak{d}(\ker z)(\im f')=0\]
In the notation of \ref{prop:ydefmappingsequence}, \( \ker z \) consists of 

\begin{enumerate}
\item \( Y_{n}(\mathfrak{R})\otimes L_{n}(\mathfrak{R})^{+} \), where \( L_{n}(\mathfrak{R})^{+} \)
denotes anything above dimension \( 0 \).

Our conclusion in this case follows simply by the definition of the m-structure
of \( B\otimes _{\alpha }\cobar B \) since the image of \( f' \) consists
of elements of the form \( e\otimes 1\in B\otimes _{\alpha }\cobar B \) and
the action of \( L_{n}(\mathfrak{R})^{+} \) on \( 1\in \cobar B \) is trivial. 

\item the summands 
\[
U_{\alpha }\otimes L(\mathfrak{R})_{0}\subset Y(\mathfrak{R})\otimes _{\rho }L(\mathfrak{R})\]
 with \( |\alpha |>0 \).
\end{enumerate}
In this case, our conclusion follows from the fact that we are \emph{really}
interested in the m-structure of \( (B\otimes _{\alpha }\cobar B)\otimes _{\cobar B}\cobar C \)
--- and the action of \( \im f' \) on \( \cobar C \) gets killed off by the
map \( \cobar B\rightarrow \cobar C \).

This completes the proof.

\end{proof}

\section{Fibrations in \protect\( \mathfrak{M}\protect \)\label{sec:fibrationsinmfrak}}

In this section we will consider an analogue to fibrations in the localized
category \( \mathfrak{M} \). We will be particularly concerned with the homotopy
fiber of a morphism.

We begin by defining the homotopy fiber of a morphism:

\begin{defn}
\label{def:homotopyfibermfrac}Let $m=\xymatrix{{A} \ar[r]^{f}&{Z}& {B} \ar[l]_{\iota} }$
be a morphism in \( \mathfrak{M} \) from \( A \) to \( B \) in canonical
form (see \ref{cor:reducehammock}), where \( f \) is a morphism in \( \mathfrak{M}_{0} \)
and \( \iota  \) is an elementary equivalence in \( \mathfrak{M}_{0} \). Then
we define the \emph{homotopy fiber} of \( m \) to be 
\[
F(m)=A\otimes _{\alpha \circ f}\cobar Z=(A\otimes _{\alpha }\cobar A)\otimes _{\cobar A}\cobar Z\]
This comes equipped with a canonical projection
\[
p:F(m)\rightarrow A\]

\end{defn}
\begin{rem}
Note that we only get a twisted tensor product with an object \emph{equivalent}
to \( B \) in \( \mathfrak{M} \) --- not \( B \) itself. This is due to the
fact that there may not even \emph{exist} a twisting cochain 
\[
A\rightarrow \cobar B\]
The point is that, \emph{a priori,} twisted tensor products are only defined
for morphisms in \( \mathfrak{M}_{0} \) --- which are \emph{actual maps.}
\end{rem}
\begin{prop}
\label{prop:homotopyfiberwelldefined}Under the hypotheses of \ref{def:homotopyfibermfrac},
the homotopy fiber of a morphism in \( \mathfrak{M} \) is well-defined.
\end{prop}
\begin{proof}
There are several cases to consider:
\begin{enumerate}
\item Replace \( A \) by \( A' \) equivalent to it \[\xymatrix{{A'} \ar[r]^{j_1}&{Z'}&{A} \ar[l]_{j_2}}\]
where \( j_{1} \) and \( j_{2} \) are both elementary equivalences. If we
combine this equivalence with the morphism above, we get \[\xymatrix{{A'}\ar[r]^{j_1}&{Z'}&{A} \ar[l]_{j_2}\ar[r]^{f}&{Z}&{B}\ar[l]_{\iota}}\] and
we put this into canonical form (see \ref{cor:reducehammock}) by taking the
push out \[\xymatrix{&&{T}&&\\ {A'}\ar[r]^{j_1}&{Z'}\ar@{.>}[ru]^{v_1}&{A} \ar[l]_{j_2}\ar[r]^{f}&{Z}\ar@{.>}[lu]_{v_2}&{B}\ar[l]_{\iota}}\] where
\( j_{2} \) and \( v_{2} \) are elementary equivalences in \( \mathfrak{M}_{0} \).
The commutative square in the center implies that 
\[
Z'\otimes _{\alpha \circ v_{1}}T\cong A\otimes _{\alpha \circ f}Z\]
 and the fact that \( j_{1} \) is an elementary equivalence implies that 
\[
A'\otimes _{\alpha \circ v_{1}\circ j_{1}}T\cong A\otimes _{\alpha \circ f}Z\]
as objects in \( \mathfrak{M} \).
\item Now we consider the case where \( B \) is replaced by an equivalent \( B' \) \[\xymatrix{{B} \ar[r]^{i_1}&{Z'}&{B'} \ar[l]_{i_2}}\] gives
the combined diagram \[\xymatrix{{A}\ar[r]^{f}&{Z}&{B}\ar[l]_{\iota}\ar[r]^{i_1}&{Z'}&{B'} \ar[l]_{i_2}}\] where
\( i_{1} \) and \( i_{2} \) are elementary equivalences in \( \mathfrak{M}_{0} \).
We canonicalize this diagram by taking the push out \[\xymatrix{ &&{S}&&\\ {A}\ar[r]^{f}&{Z}\ar@{.>}[ru]^{u_1}&{B}\ar[l]_{\iota}\ar[r]^{i_1}&{Z'}\ar@{.>}[lu]_{u_2}&{B'} \ar[l]_{i_2}}\] where
all four maps in the center square are elementary equivalences. It is easy to
see that
\[
A\otimes _{\alpha \circ g}\cobar Z\cong A\otimes _{\alpha \circ u_{1}\circ f}S\]
which proves our claim
\item The final case to be considered is that of replacing our original morphism \[\xymatrix{{A}\ar[r]^{f}&{Z}&{B}\ar[l]_{\iota}}\] one
\emph{equivalent} to it in the localized category \( \mathfrak{M} \). By \ref{cor:elementaryhomotopy},
it suffices to consider an arbitrary \emph{elementary homotopy} (\ref{def:ememtaryhomotopy}).
This will fit into a diagram \[\xymatrix@R=10pt{{}&{Z}\ar[dd]_{v}&{}\\{A}\ar[ur]^{f}\ar[dr]_{f'}&&{B}\ar[ul]_{\iota}\ar[dl]^{\iota'}\\{}&{Z'}&{}}\] where
\( \iota  \) and \( \iota ' \) are both elementary equivalences and the morphism
\( v \) may go up or down (without loss of generality, we assume it points
down). We claim that the morphism \( v \) must be a homology equivalence. This
follows from the fact that \( \iota  \) and \( \iota ' \) are elementary equivalences
(hence homology equivalences).

It is not hard to see that \( v \) induces a morphism in \( \mathfrak{M}_{0} \)
that is a homology equivalence
\[
1\otimes v:A\otimes _{\alpha \circ f}\cobar Z\rightarrow A\otimes _{\alpha \circ f'}\cobar Z'\]
Now corollary~\ref{cor:hequivanequiv} implies that
\[
A\otimes _{\alpha \circ f}\cobar Z\cong A\otimes _{\alpha \circ f'}\cobar Z'\]
as objects in \( \mathfrak{M} \).

\end{enumerate}
\end{proof}
\begin{thm}
\label{th:exactsquaremfrak}Let \begin{equation}\xymatrix{{A} \ar[r]^{f}\ar[d]_{u}& {B}\ar[d]^{v} \\ {C}\ar[r]_{g}&{D}}\label{eq:commsquaremfrak}\end{equation}
be a commutative square of morphisms in \( \mathfrak{M} \). Then there exists
an induced morphism
\[
s:F(u)\rightarrow F(v)\]
that makes \begin{equation}\xymatrix{{F(u)}\ar[r]^{s}\ar[d]&{F(v)}\ar[d] \\ {A}\ar[r]_{f} & {B}}\label{eq:commsquaremfrak2}\end{equation}
commute, where the vertical maps are canonical projections. If the horizontal
morphisms in \ref{eq:commsquaremfrak} are equivalences, then so is \( s \).
\end{thm}
\begin{rem}
Note that the morphisms in \ref{eq:commsquaremfrak} are only morphisms in the
\emph{localized} category, \( \mathfrak{M} \) --- and the square only commutes
in this category. If we replace the morphisms by underlying chain-maps, the
square only commutes up to a \emph{chain-homotopy.}
\end{rem}
\begin{proof}
First replace \ref{eq:commsquaremfrak} by canonical representations in terms
of morphisms in the nonlocalized category \( \mathfrak{M}_{0} \). We get the
somewhat gory diagram \[\xymatrix{ {A} \ar[r]^{\mathfrak{f}}\ar[d]_{\mathfrak{u}} & {Z_1} & {B}\ar[l]_{i_1}\ar[d]^{\mathfrak{v}} \\ {T_1} & & {T_2} \\ {C}\ar[u]^{j_1}\ar[r]_{\mathfrak{g}} & {Z_2} & {D}\ar[l]^{i_2} \ar[u]_{j_2} }\]
where
\begin{enumerate}
\item all maps are morphisms in \( \mathfrak{M}_{0} \).
\item the maps \( i_{1} \), \( i_{2} \), \( j_{1} \), and \( j_{2} \) are elementary
equivalences in \( \mathfrak{M}_{0} \).
\end{enumerate}
In order to put the composites \( v\circ f \) and \( g\circ u \) into canonical
form (to make use of the diagram's commutativity) we add the push outs \[\xymatrix{ {A} \ar[r]^{\mathfrak{f}}\ar[dd]_{\mathfrak{u}} & {Z_1} \ar@{.>}[d]_{d_2}& {B}\ar[l]_{i_1}\ar[d]^{\mathfrak{v}} \\ {} &{R_2} & {T_2}\ar@{.>}[l]_{e_2} \\ {T_1}\ar@{.>}[r]_{d_1}&{R_1}&{}\\ {C}\ar[u]^{j_1}\ar[r]_{\mathfrak{g}} & {Z_2}\ar@{.>}[u]^{e_1} & {D}\ar[l]^{i_2} \ar[uu]_{j_2} }\] where
\( e_{1} \) and \( e_{2} \) are elementary equivalences. 

The homotopy fiber of \( v\circ f \) is \( A\otimes _{\alpha \circ d_{2}\circ \mathfrak{f}}\cobar R_{2} \)
and the homotopy fiber of \( g\circ u \) is \( A\otimes _{\alpha \circ d_{1}\circ \mathfrak{u}}\cobar R_{1} \).
The commutativity of the original diagram \ref{eq:commsquaremfrak} (as morphisms
in \( \mathfrak{M} \)) implies (by \ref{prop:homotopyfiberwelldefined}) that
these homotopy fibers are equal in \( \mathfrak{M} \). Choose a specific equivalence
\[
1\otimes k:A\otimes _{\alpha \circ d_{1}\circ \mathfrak{u}}\cobar R_{1}\rightarrow A\otimes _{\alpha \circ d_{2}\circ \mathfrak{f}}\cobar R_{1}\]
The morphism 
\[
s:F(u)\rightarrow F(v)\]
is the composite \[\xymatrix{{A\otimes_{\alpha\circ\mathfrak{u}}\cobar T_1}\ar[r]^-{1\otimes d_1}&{A\otimes_{\alpha\circ d_1 \circ \mathfrak{u}}\cobar R_1}\ar[r]^{1\otimes k} & {A\otimes_{\alpha\circ d_2\circ \mathfrak{f} }\cobar R_2}\ar[r]^{\mathfrak{f}\otimes 1} & {Z_1\otimes_{\alpha\circ d_2}\cobar R_2} & {}\\ &&&{B\otimes_{\alpha\circ \mathfrak{v}}\cobar T_2}\ar[u]_{i_1\otimes e_2}}\] where

\begin{enumerate}
\item the single vertical map, \( i_{1}\otimes e_{2} \), is \emph{always} an elementary
equivalence in \( \mathfrak{M}_{0} \)
\item all other maps, except for \( k \), are morphisms in \( \mathfrak{M}_{0} \)
\item \( k \) is an equivalence in the localized category \( \mathfrak{M} \)
\end{enumerate}
This clearly defines a morphism in the localized category, \( \mathfrak{M} \)
that makes \ref{eq:commsquaremfrak2}. If the horizontal morphisms in our original
diagram, \ref{eq:commsquaremfrak}, were equivalences, then \( d_{1} \) and
\( \mathfrak{f} \) will be equivalences, and so will our morphism \( s:F(u)\rightarrow F(v) \).

\end{proof}

\section{Geometricity \label{sec:geometricity}}

\fancyhead[RO,LE]{\rightmark} \fancyhead[RE,LO]{Justin R. Smith}

\subsection{Definitions}

Since we are working with simplicial sets, it will be necessary to recall some
basic constructs.

\begin{defn}
\label{d:twistcartesian} Let\index{twisted Cartesian products}\index{products!twisted Cartesian}\index{fibration!in category of simplicial sets}
\( B \) and \( F \) be simplicial sets and let \( A \) be a simplicial group
acting on \( F \). Let \( \xi :B_{n}\rightarrow F_{n-1} \) be a twisting function
satisfying the identities: \( \xi (F_{n}b)\cdot F_{n-1}\xi (b)=\xi (F_{n-1}b) \);
\( \xi (D_{n}b)= \) the unit of \( A_{n} \) if \( b\in B_{n} \). We define
\( B\times _{\xi }F \), as follows: 
\begin{enumerate}
\item as a simplicial set it is the cartesian product \( B\times F \);
\item The face operators are given by \( F_{i}(b,f)=(F_{i}b,F_{i}f) \), \( 0<i \),
where \( b\in B_{n} \), \( f\in F_{n} \);
\item \( F_{0}(b,f)=(F_{0}b,\xi (b)\cdot F_{0}f) \);
\item the degeneracy operators are defined as in the cartesian product. 
\end{enumerate}
\end{defn}
\begin{rem}
Twisted cartesian products are fibrations in the context of simplicial sets
--- see \index{Gugenheim@\textsc{Victor K. A. M. Gugenheim}}\cite{Gugenheim:1959}. 
\end{rem}
Next, we recall a familiar simplicial definition of the loop space

\begin{defn}
Let \( X \) be a pointed,\index{Kan loop-group functor} simply-connected,
2-reduced simplicial set. The \emph{Kan loop-group functor} of \( X \) is a
simplicial group \( \Omega X \) defined as follows:

\( \Omega X_{n} \) is the free group generated by the simplices of \( X_{n+1} \)
subject to the relation \( D^{\Omega X}_{0}x=1 \) for any \( x\in X_{n+1} \)
\begin{enumerate}
\item \( F_{0}^{\Omega X}x=(F_{0}^{X}x)^{-1}(F_{1}^{X}x) \) for \( x\in X_{n+1} \).
\item \( F_{i}^{\Omega X}x=F_{i+1}^{X}x \) for \( x\in X_{n+1} \) and \( i>0 \).
\item \( D_{i}^{\Omega X}x=D_{i+1}^{X}x \) for \( x\in X_{n+1} \) and \( i\geq 0 \).
\end{enumerate}
\end{defn}
In \cite{Kan:1957}, Kan \index{Kan@\textsc{D. Kan}}proved that the geometric
realization of \( \Omega X \) is homotopy equivalent to the loop space of \( X \).

\begin{defn}
Let \( \eta :X\rightarrow \Omega X \) be the twisting function that sends \( x\in X_{n} \)
to \( x\in \Omega X_{n-1} \). It is not hard to see that the twisted Cartesian
product \( X\times _{\eta }\Omega X \) is contractible. This is called the
\emph{canonical acyclic twisted Cartesian product} of \( X \).
\end{defn}

\subsection{Homotopy fibers\label{sec:geometricitymstruct}}

Now we will show that the m-structure on the cobar construction (computed in
\S~\ref{sec:cobar}) is \emph{geometrically valid}, in the that it is equivalent
to the natural m-structure on the loop space. This result simplifies (and corrects)
a proof given in \cite{Smith:1994}. In in more generality, we will compute
m-structures on the homotopy fibers of maps.

Let 
\[
g:X\rightarrow Y\]
 be a map of pointed, simply-connected, 2-reduced simplicial sets. The \emph{homotopy
fiber} of this map is defined to be 
\[
F(g)=X\times _{\eta \circ g}\Omega Y\]
We can form the twisted Cartesian product of \( F(g) \) with \( Y \) --- we
use the canonical contractible twisting function sending \( Y \) to \( \mathrm{pt}\times \Omega Y \)
(where \( \mathrm{pt} \) is the basepoint of \( X \)). It has the property
that there exists a map
\begin{equation}
\label{eq:homotopfiber1}
c:X\rightarrow Y\times _{\eta '}(X\times _{\eta \circ g}\Omega Y)
\end{equation}
with the following properties:

\begin{enumerate}
\item The map \( c \) is a homotopy equivalence
\item The diagram \[\xymatrix{{}&{Y\times_{\eta'}(X\times_{\eta\circ g}\Omega Y)}\ar[d]^{p} \\ {X}\ar[r]_{g}\ar[ur]^{c} & {Y}}\]
commutes, where 
\[
p:Y\times _{\eta '}(X\times _{\eta \circ g}\Omega Y)\rightarrow Y\]
is the projection of the fibration.
\end{enumerate}
\begin{thm}
\label{theorem:mapfiber}Let 
\[
g:X\rightarrow Y\]
 be a map of pointed, simply-connected, 2-reduced simplicial sets and let
\[
f:J(\mathfrak{S})\rightarrow \mathfrak{S}\]
be the operad morphism described in \ref{prop:tensorproductprojection}. Now
regard \( \cf{{X\times _{g\circ \eta }\Omega Y}} \) as an m-coalgebra over
the operad \( J(\mathfrak{S}) \) via composition with the map \( f \). In
addition, regard \( \cobar \cf{{X}} \) as an m-coalgebra over \( J(\mathfrak{S}) \)
via projection to the base
\[
J(\mathfrak{R})\rightarrow L(\mathfrak{R})\]
so that we can set 
\[
\cf{{X}}\otimes _{\alpha \circ \cf{{g}}}\cobar \cf{{Y}}=\cf{{X}}\otimes _{\alpha }\cobar \cf{{X}}\otimes _{\cobar \cf{{X}}}\cobar \cf{{Y}}\]
making the left hand side an m-coalgebra over \( J(\mathfrak{R}) \). Then there
exists an m-coalgebra, \( V(g) \), over \( J(\mathfrak{R}) \) that is functorial
with respect to mappings, and natural elementary equivalences
\[
\cf{{X\times _{g\circ \eta }\Omega Y}}\rightarrow V(g)\leftarrow \cf{{X}}\otimes _{\alpha \circ \cf{{g}}}\cobar \cf{{Y}}\]

Consequently, \( \cf{{X}}\otimes _{\alpha \circ \cf{{g}}}\cobar \cf{{Y}} \)
is equal to the homotopy fiber in \( \mathfrak{M} \).
\end{thm}
\begin{proof}
Let 
\begin{equation}
\label{eq:homotopfiber2}
\iota :X\times _{g\circ \eta }\Omega Y\rightarrow Y\times _{\eta '}(X\times _{g\circ \eta }\Omega Y)
\end{equation}
be inclusion of the fiber and consider the commutative diagram of spaces \[\xymatrix{{X}\ar[r]^-{c}\ar[d]_{g} & {Y\times_{\eta'} (X\times_{\eta\circ g}\Omega Y)}\ar[d]^{p}\\ {Y}\ar@{=}[r] &{Y}}\]
where: 
\begin{enumerate}
\item both rows are homotopy equivalences 
\item \( c \) is the map in \ref{eq:homotopfiber1}
\item \( p:Y\times _{\eta '}(X\times _{g\circ \eta }\Omega Y)\rightarrow Y \) is
the projection of the fibration to its base
\end{enumerate}
This gives a commutative diagram of m-coalgebras, where all maps are defined
in the category \( \mathfrak{M}_{0} \) and the horizontal maps are homology
equivalences: \[\xymatrix{{\cf{X}}\ar[r]^-{\cf{c}} \ar[d]_{\cf{g}} & {\cf{Y\times_{\eta'} (X\times_{\eta\circ g}\Omega Y)}}\ar[d]^{\cf{p}}\\ {\cf{Y}}\ar@{=}[r] &{\cf{Y}}}\] Now
we take canonical acyclic twisted tensor products the copies of \( Y \) in
the lower row by \( \cobar \cf{{Y}} \) and pull them back to the upper row
to get a commutative diagram \[\xymatrix{{\cf{X}\otimes_{\alpha\circ g}\cobar \cf{Y}}\ar[r]^-{\cf{c}\otimes 1} \ar[d]_{\cf{g}\otimes 1} & {\cf{Y\times_{\eta'} (X\times_{\eta\circ g}\Omega Y)}\otimes_{\alpha\circ p} \cobar\cf{Y}}\ar[d]^{\cf{p}\otimes 1}\\ {\cf{Y}\otimes_{\alpha}\cobar\cf{Y}}\ar@{=}[r] &{\cf{Y}\otimes_{\alpha}\cobar\cf{Y}}}\] where
all maps are morphisms in \( \mathfrak{M}_{0} \). As before, both rows are
homology equivalences. We will mainly be interested in the upper row.

Since the twisting cochain \( \alpha \circ p \) vanishes on the image of 
\[
\cf{{\iota }}:\cf{{X\times _{\eta \circ g}\Omega Y}}\rightarrow \cf{{Y\times _{\eta '}(X\times _{\eta \circ g}\Omega Y)}}\]
(inclusion of the fiber), it follows that \( \cf{{\iota }} \) gives rise to
a \emph{chain-map}
\[
z:\cf{{X\times _{\eta \circ g}\Omega Y}}\rightarrow \cf{{Y\times _{\eta '}(X\times _{\eta \circ g}\Omega Y)}}\otimes _{\alpha \circ p}\cobar \cf{{Y}}\]
sending \( x\in \cf{{X\times _{\eta \circ g}\Omega Y}} \) to 
\[
\cf{{\iota }}(x)\otimes 1\in \cf{{Y\times _{\eta '}(X\times _{\eta \circ g}\Omega Y)}}\otimes _{\alpha \circ p}\cobar \cf{{Y}}\]
The lifting lemma (\ref{lem:liftinglemma}) implies that \( z \) is a morphism
in \( \mathfrak{M}_{0} \), i.e., it preserves m-structures exactly.
\begin{claim*}
The map \( z \) is a homology equivalence.
\end{claim*}
There are several ways to see this. Since we have established that \( z \)
preserves m-structures, we can forget about them and convert \( \cf{{Y\times _{\eta '}(X\times _{\eta \circ g}\Omega Y)}} \)
into multiple twisted tensor products, using the twisted Eilenberg-Zilber theorem
in \cite{Gugenheim:1972}. This will give a commutative diagram of chain-maps
(not preserving m-structures) \[\xymatrix{{\cf{X\times_{\eta\circ g}\Omega Y}}\ar[r]^-{z}\ar[rd]_{r\circ z} & {\cf{Y\times_{\eta'} (X\times_{\eta\circ g}\Omega Y)}\otimes_{\alpha\circ p}\cobar \cf{Y}}\ar[d]^{r}\\ {} & {(\cf{Y\otimes_{\eta''} \cf{X\times_{\eta\circ g}\Omega Y}}\otimes_{\alpha\circ p}\cobar \cf{Y}}}\]
and noting that we can ``rearrange'' the twisted tensor products on the bottom
term to get
\[
\cf{{Y}}\otimes _{\alpha }\cobar \cf{{Y}}\otimes _{\eta ''}\cf{{X\times _{\eta \circ g}\Omega Y}}\]
and noting that this is a twisted tensor product over an acyclic base, so that
inclusion of the fiber (which \( r\circ z \) is) induces homology isomorphisms.

\end{proof}
We can immediately conclude that our m-structure on the cobar construction is
geometrically valid:

\begin{cor}
\label{cor:cobarmstructgeo}Let \( X \) be a pointed, simply-connected, 2-reduced
simplicial set and let 
\[
f:J(\mathfrak{S})\rightarrow \mathfrak{S}\]
be the operad morphism described in \ref{prop:tensorproductprojection}. Now
regard \( \cf{{\Omega X}} \) as an m-coalgebra over \( J(\mathfrak{S}) \)
via composition of its structure map with \( f \). In addition, regard \( \cobar \cf{{X}} \)
as an m-coalgebra over \( J(\mathfrak{S}) \) via projection to the base
\[
J(\mathfrak{R})\rightarrow L(\mathfrak{R})\]
Then \( \cf{{\Omega X}} \) and \( \cobar \cf{{X}} \) are equivalent objects
in \( \mathfrak{M} \), via natural elementary equivalences
\[
\cf{{\Omega X}}\rightarrow G(X)\leftarrow \cobar \cf{{X}}\]

\end{cor}
\begin{rem}
The underlying map of this equivalence is very complex but it and \( G(X) \)
are natural with respect to simplicial maps.

Besides being much more lengthy than the present proof, the proof of this result
given in \cite{Smith:1994} is \emph{incorrect}. It relied on forming the bar
construction of \( \cf{{G}} \) where \( G \) is a simplicial group. Unfortunately,
the product operation of \( \cf{{G}} \) does \emph{not} preserve m-structures
(because it incorporates Eilenberg-Zilber maps) so that the bar construction
has no well-defined m-structure.
\end{rem}
\begin{proof}
This follows from \ref{theorem:mapfiber} by making \( X \) equal the basepoint
of \( Y \).
\end{proof}

\section{An equivalence of categories\label{ch:equivcategories}}

\subsection{Cellular m-coalgebras}

\begin{defn}
\label{def:scellular} Let \( C \) be an m-coalgebra with structure map 
\[
f:U\rightarrow \coend (C)\]
 where \( U \) is some \( E_{\infty } \)-operad. Then \( C \) will be called
\textit{strictly cellular} if there exist morphisms of \( E_{\infty } \)-operads:
\[
g_{k}:U\rightarrow \mathfrak{S}\]
 making \( S_{k,n_{k}} \) an m-coalgebra over \( U \), where \( S_{k,n_{k}} \)
is the canonical m-coalgebra of the singular complex of a wedge of spheres (see
4.2 on page 30 of \cite{Smith:1994}) and isomorphisms of m-coalgebras 
\[
f_{k}:S_{k,n_{k}}=\cf{{\Delta (\left( \bigvee _{i=1}^{n_{k}}S^{k-1}\right) )}}\rightarrow C(k-1)\]
 such that 
\[
C(k)=\cf{\Delta (\left( \bigvee _{i=1}^{n_{k}}D^{k}\right) )}\bigcup _{f_{k}}C(k-1)\]
 for all \( k\geq 0 \), where \( \Delta (*) \) is the singular complex functor.
Here, \( C(k) \) denotes the \( k \)-skeleton of \( C \) and the \( D^{k} \)
are disks whose boundaries are the \( S^{k-1} \).

We define the union over \( f \), in this case, to be the push out of the diagram
\[ \xymatrix{
{S_{k,n_{k}}}\ar[r]^{f_k}\ar[d]_{i_k}&{C(k-1)}\\
{\bigvee_{i=1}^{n_{k}}D^{k}}
}
\] 

We will call an m-coalgebra \textit{cellular} if it is equivalent (in \( \mathfrak{M} \))
to a strictly cellular m-coalgebra. 
\end{defn}
\begin{rem}
Note that cellularity requires the m-structure of an m-coalgebra to be an iterated
extension of m-structures of spheres.

Clearly, the canonical m-coalgebra of any CW-complex is cellular. The converse
also turns out to be true --- see \ref{cor:realizability}.
\end{rem}

\subsection{Morphisms}

\label{sec:mainresults}

\begin{lem}
\label{lem:realizationmapsrealizable}In this section, we will prove the main
results involving the topological realization of m-coalgebras and morphisms. 

Let \( X \) and \( Y \) be pointed, simply connected, 2-reduced simplicial
sets, and let
\[
f:|X|\rightarrow |Y|\]
 be a map of topological realizations. Then \( f \) induces a morphism
\[
\cf{{f}}:\cf{{X}}\rightarrow \cf{{Y}}\]
in \( \mathfrak{M} \). Homotopic maps of the topological realizations induce
the same morphism.
\end{lem}
\begin{rem}
Let \( i_{X}:\cf{{X}}\rightarrow \cf{{\ddelta (X)}} \) and \( i_{Y}:\cf{{Y}}\rightarrow \cf{{\ddelta (Y)}} \)
be the elementary equivalences defined in \ref{cor:singularequivalent}. Now
define \[\cf{f}=\xymatrix@C+20pt{{\cf{X}}\ar[r]^-{\cf{\ddelta(f)} \circ i_X}&{\cf{\ddelta(Y)}} & {\cf{Y}}\ar[l]_-{i_Y}}\]
\end{rem}
\begin{proof}
This follows immediately from the equivalence of homotopy theories of pointed,
2-reduced, simplicial sets and pointed simply connected spaces --- by the adjoint
functors \( \ddelta (*) \) (2-reduced singular complex) and \( |*| \) (topological
realization). See \cite{Quillen:1967} for details.
\end{proof}
We begin with a proof that \textit{equivalences} of topologically realizable
m-coalgebras are topologically realizable.

\begin{thm}
\label{th:equiv} Let \( X_{1} \) and \( X_{2} \) be pointed, simply-connected,
2-reduced simplicial sets, with associated canonical m-coalgebras, \( \mathscr{C}(X_{i}) \),
\( i=1,2 \).

In addition, suppose there exists an equivalence of m-coalgebras 
\[
f:\mathscr{C}(X_{1})\rightarrow \mathscr{C}(X_{2})\]
 in the localized category, \( \mathfrak{M} \)

Then there exists a map of topological realizations

\[
\hat{f}:|X_{1}|\rightarrow |X_{2}|\]
 making: \begin{equation}\xymatrix{ {\cf{X_1}}\ar[r]^{f}\ar[d] & {\cf{X_2}} \ar[d]\\ {\cf{\ddelta(X_1)}}\ar[r]_{\hat{f}} & {\cf{\ddelta(X_2)}} }\label{eq:equivconclusion}\end{equation}
commute in \( \mathfrak{M} \), where \( \ddelta (*) \) is the 2-reduced singular
complex and the downward maps are induced by the inclusions. Consequently, any
equivalence in \( \mathfrak{M} \) of topologically realizable m-coalgebras
is topologically realizable. 
\end{thm}
\begin{rem}
We work in the simplicial category because the functors \( \mathscr{C}(*) \)
were originally defined over it.
\end{rem}
\begin{proof}
The hypothesis implies that the chain-complexes are chain-homotopy equivalent,
hence that the \( X_{i} \), \( i=1,2 \), have the same homology. This implies
that the lowest-dimensional nonvanishing homology groups --- say \( M \) in
dimension \( k \) --- are isomorphic. We get a diagram of morphisms in \( \mathfrak{M}_{0} \)
 \begin{equation} \xymatrix{
{\mathscr{C}(X_{1})}\ar[r]^{\mathfrak{f}}
\ar[d]_{\mathscr{C}(c_{1})}
&{Z}&{\mathscr{C}(X_{2})}\ar[l]_{\iota}
\ar[d]^{\mathscr{C}(c_{2})}\\
{\mathscr{C}(K(M,k))}\ar[rr]_{=}&&{\mathscr{C}(K(M,k))}
}
\label{dia:firstmap}
  \end{equation}Where
\begin{enumerate}
\item The maps \( \{\mathscr{C}(c_{i})\} \), \( i=1,2 \), are induced by geometric
classifying maps;
\item horizontal maps are equivalences in \( \mathfrak{M} \);
\end{enumerate}
\begin{claim*}
Without loss of generality, we may assume that this diagram commutes in the
localized category, \( \mathfrak{M} \). Note that, on the chain-level, this
means that the diagram is only homotopy-commutative.
\end{claim*}
This claim follows from:

\begin{enumerate}
\item we can make the diagram commute \emph{exactly} if we replace \( X_{1} \) and
\( X_{2} \) by the associated \emph{, 2-reduced singular} simplicial sets,
\( \ddelta (X_{i}) \). This is because the obstruction to altering the maps
\( \cf{{c_{i}}} \) vanishes above dimension \( k \). 
\item \( \cf{{*}} \) of the singular complexes is canonically equivalent, in \( \mathfrak{M} \),
to \( \cf{*} \) of the originals simplicial sets, by \ref{cor:singularequivalent}.
\end{enumerate}
Theorem~\ref{th:exactsquaremfrak} implies that there exists an equivalence
in \( \mathfrak{M} \)
\[
\hat{f}:{\mathscr{C}(X_{1})}\otimes _{\alpha \circ \mathscr{C}(g_{1})}{\cobar \mathscr{C}(K(M,k))}\rightarrow {\mathscr{C}(X_{2})}\otimes _{\alpha \circ \mathscr{C}(g_{1})}{\cobar \mathscr{C}(K(M,k))}\]
 such that the following diagram commutes:  \begin{equation} \xymatrix{
{\mathscr{C}(X_1)}\otimes_{\alpha\circ \mathscr{C}(g_1)}
{\cobar \mathscr{C}(K(M,k))}\ar[r]^{\hat{f}}\ar[d]_{1\otimes \epsilon}
&{\mathscr{C}(X_2)}\otimes_{\alpha\circ \mathscr{C}(g_1)}
{\cobar \mathscr{C}(K(M,k))}\ar[d]_{1\otimes \epsilon}\\
{\mathscr{C}(X_1)}\ar[r]_{f}
&{\mathscr{C}(X_2)}
}
\label{dia:cobarcommute2}
\end{equation} 

Now, theorem~\ref{theorem:mapfiber} implies equivalences in \( \mathfrak{M} \)
\[
\mathscr{C}(X_{i}\times _{\hat{\alpha }\circ g_{i}}\Omega K(M,k))\rightarrow {\mathscr{C}(X_{i})}\otimes _{\alpha \circ \mathscr{C}(g_{i})}{\cobar \mathscr{C}(K(M,k))}\]
 for \( i=1,2 \)

We conclude that there is an equivalence 
\[
\hat{F}:\mathscr{C}(X_{1}\times _{\hat{\alpha }\circ g_{1}}\Omega K(M,k))\rightarrow \mathscr{C}(X_{2}\times _{\hat{\alpha }\circ g_{2}}\Omega K(M,k))\]
 where \( \Omega (*) \) denotes the loop space functor and \( \hat{\alpha }:K(M,k)\rightarrow \Omega K(M,k) \)
is the canonical twisting function (defining a fibration as twisted Cartesian
product --- see \cite{Gugenheim:1959}).

In addition, the commutativity of \ref{dia:cobarcommute2} implies that 
\[
f^{*}(\mu _{2})=\mu _{1}\in H^{k+1}(X_{1},M)\]
 where \( \mu _{1} \) and \( \mu _{2} \) are the \( k \)-invariants of the
fibrations \( X_{1}\times _{\hat{\alpha }\circ g_{1}}\Omega K(M,k) \) and \( X_{2}\times _{\hat{\alpha }\circ g_{2}}\Omega K(M,k) \),
respectively.

Since the \( X_{i}\times _{\hat{\alpha }\circ g_{i}}\Omega K(M,k) \) are homotopy
fibers of the \( g_{i} \) maps for \( i=1,2 \), respectively, we conclude
that the second stage of the Postnikov towers of \( X_{1} \) and \( X_{2} \)
are equivalent.

A straightforward induction implies that all \textit{finite stages} of the Postnikov
tower of \( X_{1} \) are equivalent to corresponding finite stages of the Postnikov
tower of \( X_{2} \). We ultimately get the following commutative diagram in
\( \mathfrak{M} \):

\begin{equation}\xymatrix{{\cf{P_1}}\ar[r]\ar[d]_{\cf{p_1}}& {\cf{P_2}}\ar[d]^{\cf{p_2}} \\ {\cf{X_1}}\ar@{.>}[r]_{f} & {\cf{X_2}}}\label{eq:equiv1}\end{equation}

Passing to topological realizations gives us a \emph{homotopy} commutative diagram of
\emph{spaces} \begin{equation}\xymatrix{{|P_1|}\ar[r]\ar[d]_{p_1}& {|P_2|}\ar[d]^{p_2} \\ {|X_1|}\ar@{.>}[r]_{f'} & {|X_2|}}\label{eq:equiv2}\end{equation}
where:
\begin{enumerate}
\item the \( P_{i} \) are the Postnikov towers
\item \( p_{i}:P_{i}\rightarrow X_{i} \) are the canonical projections (homotopy
equivalences)
\item \( f':|X_{1}|\rightarrow |X_{2}| \) is a map inserted into the diagram to make
it homotopy-commutative.
\end{enumerate}
By lemma~\ref{lem:realizationmapsrealizable}, diagram~\ref{eq:equiv2} induces
a diagram that commutes in \( \mathfrak{M} \): \begin{equation}\xymatrix{{\cf{\ddelta(P_1)}}\ar[r]\ar[d]_{\cf{\ddelta(p_1)}}& {\cf{\ddelta(P_2)}}\ar[d]^{\cf{\ddelta)(p_2)}} \\ {\cf{\ddelta(X_1)}}\ar[r]_{\cf{\ddelta(f')}} & {\cf{\ddelta(X_2)}}}\label{eq:equiv3}\end{equation}

Combining diagram~\ref{eq:equiv1} with \ref{eq:equiv3} (and invoking \ref{lem:realizationmapsrealizable})
give the following commutative diagram in \( \mathfrak{M} \) \begin{equation}\xymatrix{{\cf{X_1}}\ar[r]^{f}\ar[d]_{i_1}& {\cf{X_2)}}\ar[d]^{i_2} \\ {\cf{\ddelta(X_1)}}\ar[r]_{\cf{\ddelta(f')}} & {\cf{\ddelta(X_2)}}}\end{equation}
where \( i_{j}:\cf{{X_{j}}}\rightarrow \cf{{\ddelta (X_{j})}} \), \( j=1,2 \)
are the elementary equivalences induced by the canonical inclusions. This proves
the result.

\end{proof}
Next, we prove that arbitrary morphisms in \( \mathfrak{M} \) are topologically
realizable.

\begin{thm}
\label{th:morphism} Let \( X_{1} \) and \( X_{2} \) be pointed, simply-connected,
locally-finite, semisimplicial sets complexes, with associated canonical m-coalgebras,
\( \mathscr{C}(X_{i}) \), \( i=1,2 \).

In addition, suppose there exists a morphism in the category \( \mathfrak{M} \)
(see~ \ref{def:mcat}): 
\[
f:\mathscr{C}(X_{1})\rightarrow \mathscr{C}(X_{2})\]

Then there exists a map of topological realizations
\[
\hat{f}:|X_{1}|\rightarrow |X_{2}|\]
 making \begin{equation}\xymatrix{ {\cf{X_1}}\ar[r]^{f}\ar[d] & {\cf{X_2}} \ar[d]\\ {\cf{\ddelta(X_1)}}\ar[r]_{\hat{f}} & {\cf{\ddelta(X_2)}} }\label{eq:equivconclusion2}\end{equation} commute
in \( \mathfrak{M} \), where \( \ddelta (*) \) is the 2-reduced singular complex
and the downward maps are induced by the inclusions.

Consequently, any morphism of m-coalgebras is topologically realizable up to
equivalence in \( \mathfrak{M} \). 
\end{thm}
\begin{proof}
We prove this result by an inductive argument somewhat different from that used
in theorem~\ref{th:equiv}.

We build a sequence of fibrations \[ \xymatrix{
F_i\ar[d]_{p_i}\\
X_2
}
\] over \( X_{2} \) in such a way that 
\begin{enumerate}
\item the morphism \( f:\mathscr{C}(X_{1})\rightarrow \mathscr{C}(X_{2}) \) lifts
to \( \mathscr{C}(F_{i}) \) --- i.e., we have commutative diagrams  \[ \xymatrix{
&{\mathscr{C}(F_i)}\ar[d]^{\mathscr{C}(p_i)}\\
{\mathscr{C}(X_1)}\ar[r]_{f}\ar[ur]^{f_i}&{\mathscr{C}(X_2)}
}
\] 

For all \( i>0 \), \( F_{i} \) will be a fibration over \( F_{i-1} \) with
fiber a suitable Eilenberg-MacLane space.

\item The map \( f_{i} \) is \( i \)-connected in homology. 
\end{enumerate}
If the morphism \( f \) were geometric, we would be building its \textit{Postnikov
tower.}

Assuming that this inductive procedure can be carried out, we note that it forms
a convergent sequence of fibrations (see \cite{Spanier:1966}, chapter 8, \S~3).
This implies that we may pass to the inverse limit and get a commutative diagram
\[ \xymatrix{
&{\mathscr{C}(F_{\infty})}\ar[d]^{\mathscr{C}(p_{\infty})}\\
{\mathscr{C}(X_1)}\ar[r]_{f}\ar[ur]^{f_{\infty}}&{\mathscr{C}(X_2)}
}
\] where \( \tilde{f}_{\infty } \) is a morphism of m-coalgebras that is a
\textit{homology equivalence.} Now Lemma~\ref{lem:homologyequival} implies
that \( f_{\infty } \) is an equivalence of m-coalgebras, and \ref{th:equiv}
implies that it is topologically realizable.

It follows that we get a (geometric) map 
\[
\bar{f}_{\infty }:X_{1}\rightarrow F_{\infty }\]
 and the composite of this with the projection \( p_{\infty }:F_{\infty }\rightarrow X_{2} \)
is a topological realization of the original map \( f:\mathscr{C}(X_{1})\rightarrow \mathscr{C}(X_{2}) \).
The commutativity of diagram~\ref{eq:equivconclusion2} follows from that of
\ref{eq:equivconclusion}.

It only remains to verify the inductive step:

Suppose we are in the \( k^{\mathrm{th}} \) iteration of this inductive procedure.
Then the mapping cone, \( \mathscr{A}(f) \) is acyclic below dimension \( k \).
Suppose that \( H_{k}(\mathscr{A}(f_{k}))=M \). Then we get a long exact sequence
in cohomology: \begin{multline} {\dots}\to H^{k+1}(X_{1},M)\to H^{k}(\mathscr{A}(f_{k}),M)=\homz(M,M)\\
\to H^{k}(F_{k},M)\to H^{k}(X_{1},M)\to 0
\label{fib:exact}
\end{multline}

Let \( \mu \in H^{k}(F_{k},M) \) be the image of \( 1_{M}\in H^{k}(\mathscr{A}(f_{k}),M)=\homz (M,M) \)
and consider the map 
\[
h_{\mu }:X_{2}\rightarrow K(M,k)\]
 classified by \( \mu  \). We pull back the contractible fibration 
\[
K(M,k)\times _{\bar{\alpha }}\Omega K(M,k)\]
 over \( h_{\mu } \) (or form the homotopy fiber of \( h_{\mu } \)) to get
a fibration 
\[
F_{k+1}=F_{k}\times _{\bar{\alpha }\circ h_{\mu }}\Omega K(M,k)\]
 where, as before, \( \Omega (*) \) represents the loop space.

\textbf{Claim:} The morphism \( f_{k} \) lifts to a morphism \( f_{k+1}:\mathscr{C}(X_{1})\rightarrow \mathscr{C}(F_{k+1}) \)
in such a way that the following diagram commutes:  \[ \xymatrix{
{\mathscr{C}(X_1)}\ar[r]^-{f_{k+1}}\ar[rd]_{f_k}
&{\mathscr{C}(F_{k+1})}\ar[d]^{p}\\
&{\mathscr{C}(F_k)}
}
\] where \( p_{k+1}':F_{k+1}\rightarrow F_{k} \) is that fibration's projection
map.

\textbf{Proof of Claim:} We begin by theorem~\ref{theorem:mapfiber} to conclude
the existence of a commutative diagram:  \begin{equation} \xymatrix{
{\mathscr{C}(F_k\times_{\bar{\alpha}\circ h_{\mu}}\Omega K(M,k))}
\ar[r]^{e}\ar[dr]_{p'_{k+1}}
&{\mathscr{C}(F_k)
\otimes_{\bar{\alpha}\circ h_{\mu}}\cobar \mathscr{C}(K(M,k))}
\ar[d]^{1\otimes \epsilon}\\
&{\mathscr{C}(F_k)}
}
\label{dia:timesequiv}
\end{equation}  where \( e \) is an m-coalgebra equivalence.

If we pull back this twisted tensor product over the map \( f_{k} \), we get
a trivial twisted tensor product (i.e., an untwisted tensor product), because
the image of \( f^{*}(\mu )=0\in H^{k}(X_{1},M) \), by the exactness of \ref{fib:exact}.
The Lifting Lemma (\ref{lem:liftinglemma}) implies the existence of a morphism
\begin{multline} {\mathscr{C}(X_{1})}\to {\mathscr{C}(X_{1})}\otimes 1\subset 
{\mathscr{C}(X_{1})}\otimes \cobar {\mathscr{C}(K(M,k))}\\ \to
{\mathscr{C}(F_{k})}\otimes_{\bar{\alpha}\circ h_{\mu}}
\cobar{\mathscr{C}(K(M,k))}
\end{multline} The composition of this map with \( e \) in \ref{dia:timesequiv}
is the required map 
\[
{\mathscr{C}(X_{1})}\rightarrow {\mathscr{C}(F_{k+1})}\]

To see that \( H_{k}(\mathscr{A}(f_{k+1}))=0 \), note that: 
\begin{enumerate}
\item \( \mu \in H^{k}(F_{k},M)=H^{k}({\mathscr{C}(F_{k})},M) \) is the pullback
of the class in \( H^{k}(\mathscr{A}(f_{k}),M) \) inducing a homology isomorphism
\[
\mu :H_{k}(\mathscr{A}(f_{k}))\rightarrow H_{k}(K(M,k))\]
 (by abuse of notation, we identify \( \mu  \) with a cochain) or 
\[
\mu :H_{k}(\mathscr{A}(f_{k}))\rightarrow H_{k}({\mathscr{C}(K(M,k))})\]

\item in the stable range, \( {\mathscr{C}(F_{k})}\otimes _{\bar{\alpha }\circ h_{\mu }}\cobar {\mathscr{C}(K(M,k))} \)
is nothing but the \textit{algebraic mapping cone} of the chain-map, \( \mu  \),
above. But the algebraic mapping cone of \( \mu  \) clearly has vanishing homology
in dimension \( k \) since \( \mu  \) induces homology isomorphisms. 
\end{enumerate}
The conclusion follows by obstruction theory.

\end{proof}
\begin{cor}
\label{cor:realizability} An m-coalgebra is topologically realizable if and
only if it is cellular (see \ref{def:scellular}). 
\end{cor}
\begin{rem}
Clearly, topologically realizable m-coalgebras are cellular.
\end{rem}
Now we combine the results above with \ref{th:cfhomotopmfrak} to conclude:

\begin{thm}
\label{cor:realizable}The functor 
\[
\mathscr{C}(*):\underline{\mathrm{Homotop}}_{0}\rightarrow \mathfrak{M}^{+}\]
 defines an equivalence of categories and homotopy theories (in the sense of
\cite{Quillen:1967}), where \( \underline{\mathrm{Homotop}}_{0} \) is the
homotopy category of pointed, simply-connected CW-complexes and continuous maps
and \( \mathfrak{M}^{+}\subset \mathfrak{M} \) is the subcategory of topologically
realizable m-coalgebras. 
\end{thm}
\begin{rem}
Recall that, in \cite{Quillen:1967}, Quillen defined a \emph{homotopy theory}
to a homotopy \emph{category} augmented with a \emph{loop-space} and \emph{suspension}
functor. These additional structures are necessary to carry out many of the
standard constructions in homotopy theory.
\end{rem}
\begin{proof}
We have already proved most of this. The statement that it defines an equivalence
of \emph{homotopy theories} follows from the fact that it carries suspensions
of CW-complexes to suspensions of m-coalgebras and loop-spaces to the cobar
construction.
\end{proof}
\begin{cor}
\label{cor:lrealizable}The restriction of the canonical functor
\[
\mathfrak{f}:\mathfrak{L}\rightarrow \mathfrak{M}\]
 to topologically realizable objects of \( \mathfrak{L} \) is full and faithful.
Consequently, the functor 
\[
\mathscr{C}(*):\underline{\mathrm{Homotop}}_{0}\rightarrow \mathfrak{L}^{+}\]
is an equivalence of categories, where \( \mathfrak{L}^{+}\subset \mathfrak{L} \)
is the subcategory of topologically realizable m-coalgebras over \( \mathfrak{S} \).
\end{cor}
\begin{rem}
This result implies that the category \( \mathfrak{L} \) encapsulates homotopy
theory to the \emph{same} extent as the category \( \mathfrak{M} \). In other
words, to study homotopy types, it suffices to consider coalgebras over the
operad \( \mathfrak{S} \) --- one \emph{never} needs to go to other \( E_{\infty } \)-operads.
We only considered these more complex operads because the bar and cobar constructions
seemed to require them. 

I conjecture our definitions of the bar and cobar constructions can be reworked
so that one \emph{remains} in the category \( \mathfrak{L} \).
\end{rem}
\begin{proof}
Suppose \( X \) and \( Y \) are semi-simplicial sets such that \( \cf{{X}},\cf{{Y}}\in \mathfrak{L} \)
become equivalent when regarded as objects of \( \mathfrak{M} \). Then there
exists a homotopy equivalence 
\[
|X|\rightarrow |Y|\]
which induces a homology equivalence
\[
\cf{{\ddelta (X)}}\rightarrow \cf{{\ddelta (Y)}}\]
where \( \ddelta (X) \) and \( \ddelta (Y) \) are the singular complexes of
\( X \) and \( Y \), respectively. The conclusion follows from the fact that
the inclusions 
\begin{eqnarray*}
\cf{{X}} & \rightarrow  & \cf{{\ddelta (X)}}\\
\cf{{Y}} & \rightarrow  & \cf{{\ddelta (Y)}}
\end{eqnarray*}
are equivalences in \( \mathfrak{L} \) --- see \ref{cor:singularequivalent}.
Given any morphism in \( \mathfrak{M} \)
\[
f:\cf{{X}}\rightarrow \cf{{Y}}\]
 we know that there exists a geometric map
\[
|X|\rightarrow |Y|\]
 realizing it. But this implies the existence of a corresponding morphism of
singular complexes and a morphism in \( \mathfrak{L} \).
\end{proof}
In the previous results' proof we used intermediaries that are (almost always)
\emph{uncountably generated} to encapsulate homotopy theory. It is interesting
to note that this is not necessary when one restricts oneself to \emph{finite}
semi-simplicial sets.

Recall the definition of \( \mathfrak{F}_{0} \) in \ref{def:L0category} and
the inclusion functor
\[
\mathfrak{F}_{0}\subset \mathfrak{L}_{0}\]
inducing a functor of localizations
\[
\mathfrak{i}:\mathfrak{F}\rightarrow \mathfrak{L}\]

It is not clear that this functor is faithful. Our main result will be that
its restriction to \emph{topologically realizable m-coalgebras} is faithful.

Recall the well-known fact that the category of finite semi-simplicial sets
is equivalent to that of finite simplicial sets --- semi-simplicial sets whose
simplices are all \emph{non-degenerate.} 

We recall the important notion of \emph{barycentric subdivision:}

\begin{defn}
\label{def:barycentric}If \( X \) is a simplicial set and \( p \) is a point
(not contained in \( X \)) let \( \delta (X,p) \) denote the \emph{cone} on
\( X \) with apex \( p \). This is functorial with respect to morphisms of
pairs \( (X,p) \).

Inductively define a functor \( b(*) \) on simplices via:
\begin{enumerate}
\item If \( \sigma  \) is a \( 0 \)-simplex, \( b(\sigma )=\sigma  \).
\item If \( \sigma ^{n} \) is an \( n \)-simplex, \( b(\sigma ^{n})=\delta (b(\partial \sigma ^{n}),\beta _{\sigma }) \),
where \( \beta _{\sigma } \) is a point disjoint from \( X \) and, as the
notation implies, depending on \( \sigma  \).
\end{enumerate}
This extends to a functor \( b(*) \) on simplicial sets, called \emph{barycentric
subdivision}.
\end{defn}
\begin{rem}
It is well-known that there exists a piecewise-linear homeomorphism
\[
|b(X)|\rightarrow |X|\]
sending \( \beta _{\sigma } \) to the (geometric) barycenter of \( \sigma  \)
for all \( \sigma \in X \). This implies that \( \cf{{X}}\equiv \cf{{b(X)}} \)
in \( \mathfrak{L} \), by \ref{cor:lrealizable}. Our main result will be that,
if \( X \) is \emph{finite,} then this equivalence \emph{also} holds in \( \mathfrak{F} \).
\end{rem}
We need the following:

\begin{lem}
\label{lemma:zfunctorconstruction}There exists a functor \( Z(*) \) on the
category of simplicial sets such that there exist functorial inclusions
\begin{eqnarray}
X & \subset  & Z(X)\label{eq:zfunctorinclusions1} \\
b(X) & \subset  & Z(X)\label{eq:zfunctorinclusions2} 
\end{eqnarray}
that are both homology equivalences.
\end{lem}
\begin{proof}
In fact \( |X| \) will be a strong deformation retract of \( Z(X) \) via a
functorial retraction, and the inclusion of \( |b(X)|\subset |Z(X)| \) will
be a homotopy equivalence. We define \( Z(*) \) inductively on simplices via:
\begin{enumerate}
\item If \( \sigma  \) is a \( 0 \)-simplex, \( Z(\sigma )=\delta (\sigma ,p_{\sigma }) \),
where \( p_{\sigma } \) is some point disjoint from \( \sigma  \) and depending
on \( \sigma  \). Note that \( Z(\sigma )=\sigma \times I \).
\item If \( \sigma  \) is an \( n \)-simplex, define \( Z(\sigma )=\delta (\sigma \cup Z(\partial \sigma ),p_{\sigma }) \),
where \( p_{\sigma } \) is some point disjoint from \( \sigma \cup Z(\partial \sigma ) \)
--- and depending on \( \sigma  \).
\end{enumerate}
The following diagram illustrates this construction:

\vspace{0.3cm}
{\par\centering \resizebox*{9cm}{8cm}{\includegraphics{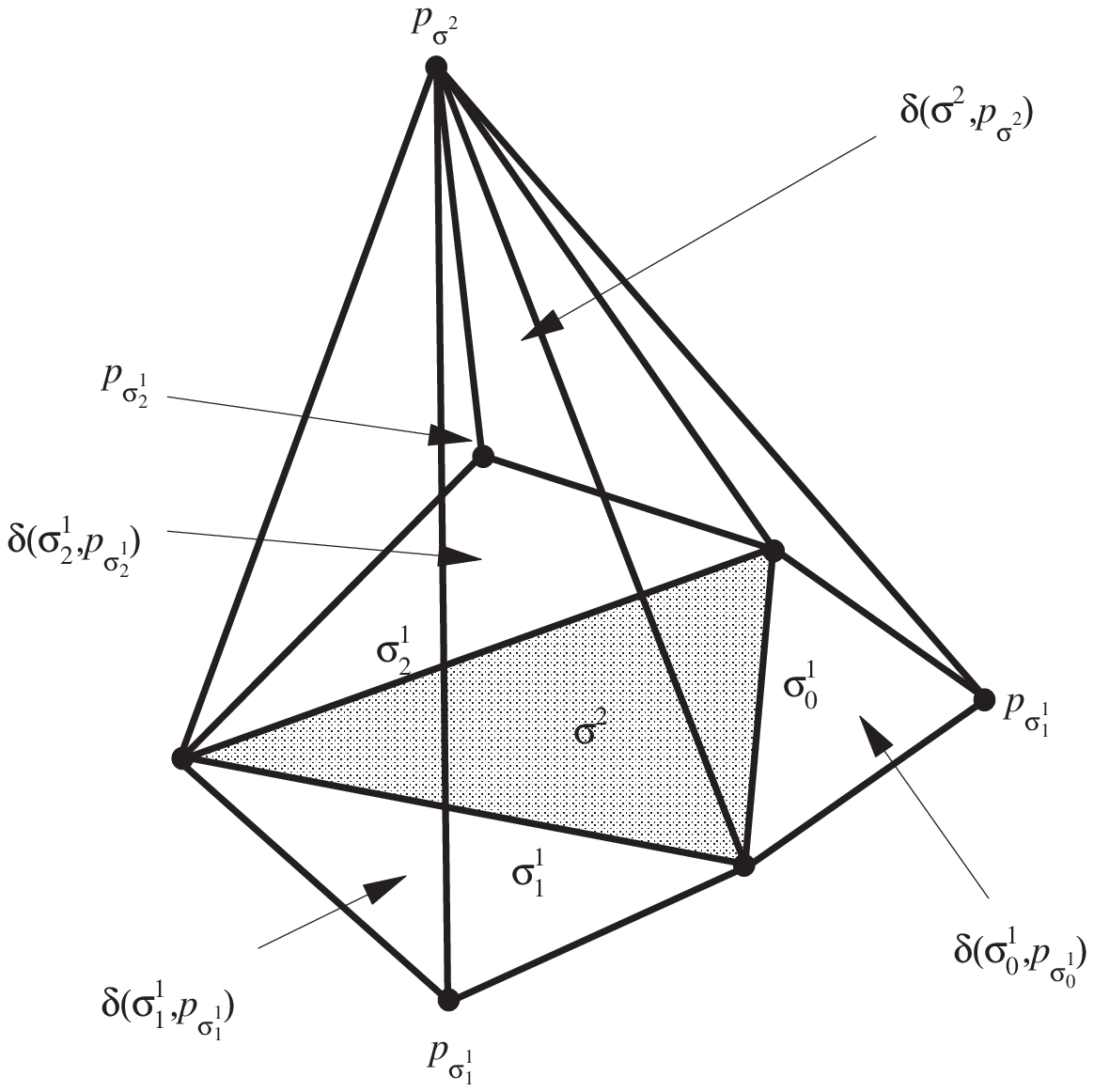}} \par}
\vspace{0.3cm}

where the shaded area (\( \sigma ^{2} \)) was in \( X \).

It is easy to see that this construction is functorial with respect to inclusions
of faces so it extends to a functor on the category of simplicial sets. For
every simplicial set \( X \) there exists a canonical inclusion \( X\subset Z(X) \)
and a set of \( 0 \)-simplices \( p_{*}\subset Z(X) \) indexed by the simplices
of \( X \).

It is also trivial (see \ref{def:barycentric}) to see that the inclusion \ref{eq:zfunctorinclusions2}
exists and is functorial. 

We claim that there exists a functorial retraction
\[
r_{X}:|Z(X)|\rightarrow |X|\]
 making \( |X| \) a strong deformation retract of \( |Z(X)| \). Simply define
\( r_{X} \) to be the unique piecewise-linear map fixing \( |X|\subset |Z(X)| \)
and sending \( p_{\sigma } \) to the barycenter of \( |\sigma | \) for all
\( \sigma \in X \). The restriction of this map to \( |b(X)|\subset |Z(X)| \)
coincides with the well-known homeomorphism \( |b(X)|\rightarrow |X| \). It
follows that the inclusions \ref{eq:zfunctorinclusions1} and \ref{eq:zfunctorinclusions2}
both induce homology isomorphisms.
\end{proof}
\begin{cor}
\label{cor:barycentricequiv}If \( X \) is a finite, pointed, simply-connected
simplicial set, then \( \cf{{X}} \) and \( \cf{{b(X)}} \) define the same
element of \( \mathfrak{F} \).
\end{cor}
Now we can appeal to the well-known \emph{simplicial approximation theorem}
to conclude that:

\begin{prop}
Let \( X \) and \( Y \) be finite, pointed, simply-connected simplicial sets
and let 
\[
f:|X|\rightarrow |Y|\]
 be a map. Then there exists a well-defined induced morphism
\[
f_{*}:\cf{{X}}\rightarrow \cf{{Y}}\]
 in the category \( \mathfrak{F} \) depending on the homotopy class of \( f \).
This morphism is an equivalence if and only if \( f \) was a homotopy equivalence.
\end{prop}
\begin{proof}
The simplicial approximation theorem says that there exist integers \( n,m>0 \)
and a simplicial map
\[
g:b^{n}(X)\rightarrow b^{m}(Y)\]
such that 
\[
|g|:|X|\rightarrow |Y|\]
 is homotopic to \( f \) --- where we have identified \( |X| \) with \( |b^{n}(X)| \)
and \( |Y| \) with \( |b^{m}(Y)| \), respectively, via the well know homeomorphisms
between them. To see that this induced map only depends on the homotopy class
of \( f \), apply this argument to \( X\times I \) and \( Y \).
\end{proof}
\begin{cor}
\label{cor:mainfiniteresult}Let \( \mathrf {F} \) denote the homotopy category
of finite, pointed, simply-connected simplicial sets. Then the functor \( \cf{{*}} \)
defines an equivalence of categories
\[
\cf{{*}}:\mathrf {F}\rightarrow \mathfrak{F}^{+}\]
where \( \mathfrak{F}^{+}\subset \mathfrak{F} \) is the sub-category of topologically
realizable, finitely-generated m-coalgebras over \( \mathfrak{S} \).
\end{cor}
We may state our main result in non-categorical terms:

\begin{cor}
\label{cor:mainnoncat}Let \( X \) and \( Y \) be pointed, simply-connected
semisimplicial sets and let 
\[
f:C(X)\rightarrow C(Y)\]
be a chain-map between canonical chain-complexes. Then \( f \) is topologically
realizable if and only if there exists an m-coalgebra \( C \) over \( \mathfrak{S} \)
and a factorization \( f=f_{\beta }\circ f_{\alpha } \) \[\xymatrix{{\cf{X_1}}\ar[r]^-{f_{\alpha}}&{C}\ar@<.5ex>@{ .>}[r]^-{f_{\beta}}&{\cf{X_2}}\ar@<.5ex>@{ >-}[l]^-{\iota}}\]
where \( f_{\alpha } \) is a morphism of m-coalgebras, \( \iota  \) is an
\emph{elementary equivalence} --- an injection of m-coalgebras with acyclic,
\( \integers  \)-free cokernel --- and \( f_{\beta } \) is a chain map that
is a left inverse to \( \iota  \). If \( X \) and \( Y \) are finite, we
may require \( C \) to be finitely generated.
\end{cor}
\begin{proof}
This follows directly from \ref{cor:mainfiniteresult} and \ref{cor:reducehammock}.
\end{proof}
\appendix

\section{Cobar Duality}

\label{app:cobarduality}

In this appendix, we lay the groundwork for proving \ref{cor:bigstarmaps} by
defining a duality morphism between the bar and cobar constructions somewhat
like that in \ref{lem:operadduality}.

\subsection{Definitions}

\begin{defn}
\label{def:betashuffle}Let \( \beta _{1},\dots ,\beta _{k} \) be \( k \)
length-\( n \) sequences of nonnegative integers, and let \( q=\sum _{i=1}^{k}|\beta _{i}| \).
Define:
\begin{enumerate}
\item \( \mathscr{Z}\{\beta _{1},\dots ,\beta _{k}\}\in \rs{{q}} \) to be the signed
shuffle-permutation that shuffles the sequence \( \{1,\dots ,q\}=\{\beta _{1,1},\dots ,\beta _{1,n},\dots ,\beta _{k,1},\dots ,\beta _{k,n}\} \)
into \( \{\beta _{1,1},\dots ,\beta _{k,1},\dots ,\beta _{1,n},\dots ,\beta _{k,n}\} \).
The sign is just the parity of the permutation.
\item \( V(\beta _{1},\dots ,\beta _{k}) \) to be the corresponding shuffle-map \( C^{q}\rightarrow C^{q} \).
\end{enumerate}
\end{defn}
\begin{rem}
Note that \( V(\beta _{1},\dots ,\beta _{k}) \) is a map of a tensor product
of copies of \( C \) and \( \mathscr{Z}\{\beta _{1},\dots ,\beta _{k}\} \)
is an element of \( \rs{{q}} \). By definition \( \mathscr{Z}\{\alpha \} \)
(where \( \alpha  \) is a \emph{single} sequence of integers) is the identity
map.

The cobar construction is a free DGA-algebra. Its \( n \)-fold tensor product
\( \cobar (C)^{n} \) is also a free DGA-algebra with product operation defined
by
\[
r_{1}\cdots r_{k}=V(\beta _{1},\dots ,\beta _{k})(r_{1}\otimes \cdots \otimes r_{k})\]
where \( r_{i}\in (\Sigma ^{-1}C)^{\beta _{i}} \)
\end{rem}
\begin{prop}
Let \( u(n,k) \) be the \( k\times n \) matrix 
\[
\left( \begin{array}{ccc}
1 & \dots  & n\\
n+1 & \dots  & 2n\\
\vdots  & \ddots  & \vdots \\
n(k-1)+1 & \dots  & kn
\end{array}\right) \]
let \( S(n,k) \) be the sequence that results from transposing \( u(n,k) \)
and reading off its entries from left to right by rows, and let \( p(n,k) \)
be the permutation that maps \( \{1,\dots ,nk\} \) to \( S(n,k) \). Then 
\[
\mathscr{Z}\{\beta _{1},\dots ,\beta _{k}\}=\mathrm{T}_{\{\beta _{1,1},\dots ,\beta _{1,n},\dots ,\beta _{k,1},\dots ,\beta _{k,n}\}}(p(n,k))\]

\end{prop}

\subsection{The map}

We begin by defining the duality map in the case where our operad is \( \coend (C) \)
for a chain-complex \( C \).

\begin{cor}
\label{cor:hmapdef}The composite \begin{multline}h=\homz(1,\bigoplus _{\{\beta_1,\dots,\beta_k\}}V(\beta_1,\dots,\beta_k))\circ h^k_0:\\ \left(\Sigma\homz(C,\cobar(C)^n)\right)^k \rightarrow \homz((\Sigma^{-1}C)^k,\cobar(C)^n)\end{multline}
defines an injective morphism of \( \zs{{n}} \)-chain complexes, where:
\begin{enumerate}
\item \( h^{k}_{0}:\left( \Sigma \homz (C,\cobar (C)^{n})\right) ^{k}\rightarrow \left( \Sigma \homz (\Sigma ^{-1}C,\cobar (C)^{n})\right) ^{k} \)
is the natural isomorphism
\item The direct sum is taken over all sets of \( k \) length-\( n \) sequences
of nonnegative integers, \( \{\beta _{1},\dots ,\beta _{k}\} \).
\item The summand \( V(\beta _{1},\dots ,\beta _{k}) \) is applied to \begin{multline}\homz\left(\Sigma^{-1}C,(\Sigma^{-1}C)^{\beta_{1,1}}\otimes\cdots\otimes(\Sigma^{-1}C)^{\beta_{1,n}}\right) \\\otimes \homz\left(\Sigma^{-1}C,(\Sigma^{-1}C)^{\beta_{2,1}}\otimes\cdots\otimes(\Sigma^{-1}C)^{\beta_{2,n}}\right) \\\vdots\\ \otimes\homz\left(\Sigma^{-1}C,(\Sigma^{-1}C)^{\beta_{k,1}}\otimes\cdots\otimes(\Sigma^{-1}C)^{\beta_{k,n}}\right)\\\subset \Sigma\homz\left(C,\cobar(C)^n\right) \end{multline}
\item \( S_{n} \) acts on both sides by permuting factors of the rightmost argument
of \( \homz (C,*) \).
\end{enumerate}
\end{cor}
\begin{proof}
This is straightforward: the \( V \) maps shuffle the factors of \( \cobar (C)^{n} \)
into the appropriate positions.
\end{proof}
\begin{enumerate}
\item At this point, we can try to pull back the differential of \( \homz ((\Sigma ^{-1}C)^{k},\cobar (C)^{n}) \)
(induced by \( \homz (\cobar (C),\cobar (C)^{n}) \) to \( \left( \Sigma \homz (C,\cobar (C)^{n})\right) ^{k} \).
As before, this differential is a perturbation of the tensor algebra differential,
\( \partial _{\otimes }+p \). 
\item Since \( \left( \Sigma \homz (C,\cobar (C)^{n})\right) ^{k} \) \emph{already}
has the differential of a tensor algebra, we need only perturb it by a term
corresponding to the \emph{dual} of \( p \) (since it appears as the \emph{left}
argument of a \( \homz (*,*) \)-functor). Since \( p \) is defined by a \emph{coproduct}
(and higher coproducts of \( C \)), the \emph{dual} corresponds to a \emph{product}
(with higher product-operations). We must, consequently, define an \( \ainfty  \)-algebra
structure on \( \bigoplus \left( \Sigma \homz (C,\cobar (C)^{n})\right) ^{k} \)
--- or \( Z_{n}(\coend (C)) \).
\end{enumerate}
We must, consequently, consider:

\begin{prop}
\label{prop:mumapdef}It is possible to define a set \( \{\mu _{i}\} \) of
product-operations on \( \homz (C,\cobar (C)^{n}) \) that make the following
diagram commute: \begin{equation}\xymatrix@C+20pt{{\bigotimes_{i=1}^k\Sigma\homz(C,(\Sigma^{-1}C)^{\beta_i})}\ar[r]^-{h_{\beta_1,\dots,\beta_k}}\ar[d]|-{\mu_k}&{\homz((\Sigma^{-1}C)^k,\cobar(C)^n)}\ar[d]|-{\homz(\delta_k,1)}\\{\homz(C,(\Sigma^{-1}C)^{\alpha})}\ar[r]_-{h_{\alpha}}&{\homz(\Sigma^{-1}C,\cobar(C)^n)}}\label{eq:mumapdefs}\end{equation} where
\begin{enumerate}
\item \( h_{\beta _{1},\dots ,\beta _{k}}=\homz (1,V(\beta _{1},\dots ,\beta _{k}))\circ h_{1} \)
\item \( h_{\alpha }=\homz (1,V(\alpha ))\circ h_{1}=h_{1} \)
\item \( \{\delta _{i}\} \) are the higher coproducts from the \( \ainfty  \)-coalgebra
structure of \( C \).
\end{enumerate}
\end{prop}
\begin{rem}
The \( \{h_{\beta _{1},\dots ,\beta _{k}}\} \) are just restrictions of the
\( h \)-map defined in \ref{cor:hmapdef} to various summands of \( \bigotimes \Sigma \homz (C,\cobar (C)^{n}) \).
\end{rem}
\begin{proof}
This follows from the fact that the term \( \homz (\delta _{i},1) \) is expressible
in terms of the composition operations in \( \coend (C) \), so it pulls back
to \( \Sigma \homz (C,(\Sigma ^{-1}C)^{\alpha }) \).
\end{proof}
Now we formulate a definition of the \( \{\mu _{i}\} \) for an \emph{arbitrary}
operad, rather than just \( \coend (C) \). This is a variation on \ref{prop:everyoperadalgebra}:

\begin{defn}
\label{def:lowertmaps} Given any operad \( \mathfrak{R} \), we can define
a morphism of \( \integers  \)-modules 
\[
t_{n}:\mathfrak{R}_{n}\rightarrow \zend (Z_{n}(\mathfrak{R}))\]
 as follows:

The adjoint, \( t_{n}^{*} \), is the composite \begin{equation}\xymatrix{{\mathfrak{R}_k\otimes\Sigma^{-|\beta_{1}|}\mathfrak{R}_{\beta_{1}}\otimes\dots\otimes\Sigma^{-|\beta_{k}|}\mathfrak{R}_{\beta_{k}}}\ar[d]|-{1\otimes\susp^{\beta_1}\otimes\cdots\otimes\susp^{\beta_k}}\\{\mathfrak{R}_k\otimes\mathfrak{R}_{\beta_{1}}\otimes\cdots\otimes\mathfrak{R}_{\beta_{k}}}\ar[d]|-{f\otimes1^k}\\ {\mathfrak{R}_k\otimes\mathfrak{R}_{\beta_{1}}\otimes\cdots\otimes\mathfrak{R}_{\beta_{k}}} \ar[d]|-{\bar{T}}\\ {\mathfrak{R}_{\beta_{1}}\otimes\cdots\otimes\mathfrak{R}_{\beta_{k}}\otimes\mathfrak{R}_k}\ar[d]|-{\gamma_{\mathfrak{R}}}\\{\mathfrak{R}_{\alpha }}\ar[d]|-{\desusp^{|\alpha|}\sshufp{\beta_{1},\dots,\beta_{k}}}\\ \Sigma^{-|\alpha|}\mathfrak{R}_{\alpha }}\label{dia:tnmapdef}\end{equation}
(see \ref{def:betashuffle}) where 
\[
\gamma _{\mathfrak{R}}:\mathfrak{R}_{i_{1}}\otimes \cdots \otimes \mathfrak{R}_{i_{k}}\otimes \mathfrak{R}_{k}\rightarrow \mathfrak{R}_{i_{1}+\cdots +i_{k}}\]
 defines the operad structure of \( \mathfrak{R} \) (see~\ref{def:operad}),
\[
\bar{T}:\mathfrak{R}_{k}\otimes \mathfrak{R}_{i_{1}}\otimes \cdots \otimes \mathfrak{R}_{i_{k}}\rightarrow \mathfrak{R}_{i_{1}}\otimes \cdots \otimes \mathfrak{R}_{i_{k}}\otimes \mathfrak{R}_{k}\]
 is the map that permutes the factors, and \( \alpha =\sum _{i=1}^{k}\beta _{i} \)
with the summation of sequences taken elementwise. Here the \( \{\beta _{i}\} \)
are \( k \) length-\( n \) sequences of nonnegative integers. 
\end{defn}
\begin{prop}
If \( \mathfrak{R}=\coend (C) \), the maps \( t_{k} \) defined in \ref{def:lowertmaps}
and the \( \mu  \)-maps defined in \ref{prop:mumapdef} satisfy the equation
\[
\mu _{k}=\underbrace{\susp \otimes \cdots \otimes \susp }_{k\, \mathrm{times}}\circ t_{k}\circ \desusp :\left( \Sigma \homz (C,\cobar (C)^{n}\right) ^{k}\rightarrow \homz ((\Sigma ^{-1}C)^{k},\cobar (C)^{n})\]

\end{prop}
\begin{proof}
This follows from:
\end{proof}
\begin{enumerate}
\item The \( t \)-maps in \ref{def:lowertmaps} are a direct translation of \ref{prop:mumapdef}
into operad-compositions.
\item In diagram \ref{eq:mumapdefs}, the maps \( \homz (1,V(\beta _{1},\dots ,\beta _{k})) \)
and \( \homz (\delta _{i},1) \) can be \emph{permuted} --- i.e. \( \homz (\delta _{i},1)\circ \homz (1,V(\beta _{1},\dots ,\beta _{k}))=\homz (1,V(\beta _{1},\dots ,\beta _{k}))\circ \homz (\delta _{i},1) \). 

This means that we can apply the map \( \homz (1,V(\beta _{1},\dots ,\beta _{k})) \)
\emph{last} --- and its translation into operad-compositions is \emph{precisely}
\( \mathscr{Z}\{\beta _{1},\dots ,\beta _{k}\} \).

\end{enumerate}

\section{Compositions\label{app:operadcomps}}

Now we are in a position to define composition operations on \( \barcs (Z_{n}(\mathfrak{R})) \)
and \( Y_{n}(\mathfrak{R})\bigstar _{\rho _{n}}\barcs (Z_{n}(\mathfrak{R})) \)
that make the maps in \ref{cor:bigstarmaps} into operad morphisms. 

On the summand \( \left( \Sigma \homz (C,\cobar (C)^{n}\right) ^{\alpha }\otimes \left( \Sigma \homz (C,\cobar (C)^{n}\right) ^{\beta } \),
we define \( \circ _{i} \) via
\begin{eqnarray*}
\lefteqn {\homz (1,V(\gamma _{1},\dots ,\gamma _{n+m-1})^{-1}\circ }\hspace {.5in} &  & \\
 &  & \homz (1,V(\alpha _{1},\dots ,\alpha _{n}))\circ _{i}\homz (1,V(\beta _{1},\dots ,\beta _{m}))
\end{eqnarray*}
where \( \gamma =\{\gamma _{1},\dots ,\gamma _{n+m-1}\} \) results from combining
the sequences \( \alpha  \) and \( \beta  \) suitably. We illustrate this
with diagram \ref{fig:cobarcomps}. Suppose we are trying to form the composite
\[
a\circ _{i}b\]
where 
\[
a\in \bigotimes _{\ell =1}^{j}\Sigma \homz (C,(\Sigma ^{-1}C^{\alpha _{\ell }})=A\]
 \( \alpha =\{\alpha _{1},\dots ,\alpha _{j}\} \) consists of \( j \) length-\( n \)
sequences of nonnegative integers and 
\[
b\in \bigotimes _{t=1}^{k}\Sigma \homz (C,(\Sigma ^{-1}C^{\beta _{t}})=B\]
\( \beta =\{\beta _{1},\dots ,\beta _{k}\} \) is \( k \) length-\( m \) sequences
of integers. Observe that (as in figure \ref{fig:cobarcomps}): 

\begin{enumerate}
\item We may regard \( a \) and \( b \) as maps
\[
a:(\Sigma ^{-1}C)^{j}\rightarrow \bigotimes _{\ell =1}^{j}(\Sigma ^{-1}C)^{\alpha _{\ell }}\]
and
\[
b:(\Sigma ^{-1}C)^{k}\rightarrow \bigotimes _{t=1}^{k}(\Sigma ^{-1}C)^{\beta _{t}}\]
respectively.
\item In the mapping to \( \homz (\cobar (C),(\cobar (C))^{n}), \) the \( i^{\mathrm{th}} \)
\emph{column} of the diagram corresponds to the \( i^{\mathrm{th}} \) copy
of \( \cobar (C) \) in the \emph{right} argument of \( \homz (*,*) \) ---
and that the \( t^{\mathrm{th}} \) \emph{row} of the diagram corresponds to
the \( t^{\mathrm{th}} \) copy of \( \Sigma ^{-1}C\subset \cobar (C) \) in
the \emph{left} argument. 
\end{enumerate}
Consequently, the composite \( a\circ _{i}b \) is \emph{only} nonzero if \( j=\sum \beta _{*,i} \).

\begin{figure} \begin{xy}0;<1cm,0cm>: (0,0)*+{(\Sigma^{-1}C)^k}; (7,0)*+{\begin{array}{c} (\Sigma^{-1}C)^{\beta_{1,1}}\otimes\cdots\otimes (\Sigma^{-1}C)^{\beta_{1,i}}\otimes\cdots\otimes (\Sigma^{-1}C)^{\beta_{1,m}}\\ \otimes\\\vdots\\ \otimes\\(\Sigma^{-1}C)^{\beta_{k,1}}\otimes\cdots\otimes (\Sigma^{-1}C)^{\beta_{k,i}}\otimes\cdots\otimes (\Sigma^{-1}C)^{\beta_{k,m}} \end{array}}**\dir{-} ?<>(.4)*_!/10pt/{b}?>*\dir{>}; (7,3.4)*+\txt{$i^{\mathrm{th}}$ \textit{column} \\ corresponds to the \\$i^{\mathrm{th}}$ copy of $\cobar(C)$\\ in target, $\cobar(C)^m$}; (6.85,1.6)*+{} **\dir{-} ?>*\dir{>};
(6.95,0)*+{\hbox to 1.7cm{\vbox to 2.5cm{}}} *\frm{o-}; (6.95,-1.5)*+{\hbox to 1.7cm{}} *\frm{_\}}; (6.95,-1.8)*+{}; (6.5,-3.6)*+{}**\dir{=}; (7,-5)*+{\left.\begin{array}{c} \Sigma^{-1}C\\ \otimes\\ \vdots\\ \otimes\\ \Sigma^{-1}C\end{array}\right\}j\mathrm{~times}}; (0,-5)*+{\begin{array}{c} (\Sigma^{-1}C)^{\alpha_{1,1}} \otimes \cdots \otimes (\Sigma^{-1}C)^{\alpha_{1,n}}\\ \otimes \\ \vdots \\ \otimes \\ (\Sigma^{-1}C)^{\alpha_{j,1}}\otimes\cdots\otimes (\Sigma^{-1}C)^{\alpha_{j,n}} \end{array}} **\dir{-} ?<>(.4)*_!/10pt/{a}?>*\dir{>};
\end{xy}\caption{Forming the composition $a\circ_i b$ in $\barcs{Z_n(\coend(C))}$}\label{fig:cobarcomps}\end{figure}

Given these considerations, we consider what the composition \emph{does}. We
\emph{insert} the rows of the target of \( a \) into the rows of the target
of \( b \). There is one final subtlety to this, however: we want to insert
\emph{copies of} \( \cobar (C) \) from the target of \( a \) into column \( i \)
of the target of \( b \). This means we must shuffle the copies of \( \Sigma ^{-1}C \)
in the target of \( a \) before inserting them into the target of \( b \)
--- see figure \ref{fig:operadcomps2}.

\begin{figure} $$\xymatrix{{\cdots\otimes{(\Sigma^{-1}C)^{\beta_{t,i-1}}}\otimes{(\Sigma^{-1}C)^{\beta_{t,i}}}\otimes{(\Sigma^{-1}C)^{\beta_{t,i+1}}}\otimes\cdots}\ar@{=}[d]\\ {\cdots\otimes{(\Sigma^{-1}C)^{\beta_{t,i-1}}}\otimes{\left\{\begin{array}{c}\Sigma^{-1}C\\ \otimes\\ \vdots\\ \otimes\\ \Sigma^{-1}C\end{array}\right\}}\otimes{(\Sigma^{-1}C)^{\beta_{t,i+1}}}\otimes\cdots} \ar[d]_{\cdots\otimes1\otimes a\otimes1\otimes\cdots}\\ {\cdots\otimes{\left\{\begin{array}{c}(\Sigma^{-1}C)^{\alpha_{\tilde{\beta}_{t,i},1}} \otimes \cdots \otimes (\Sigma^{-1}C)^{\alpha_{\tilde{\beta}_{t,i},n}}\\ \otimes \\ \vdots \\ \otimes \\ (\Sigma^{-1}C)^{\alpha_{\tilde{\beta}_{t,i}+\beta_{t,i}-1,1}}\otimes\cdots\otimes (\Sigma^{-1}C)^{\alpha_{\tilde{\beta}_{t,i}+\beta_{t,i}-1,n}}\end{array}\right\}\otimes\cdots}}\ar[d]_{\cdots\otimes 1\otimes V'\otimes 1 \otimes\cdots}\\{\cdots\otimes{\left\{\left[\begin{array}{c}(\Sigma^{-1}C)^{\alpha_{\tilde{\beta}_{t,i},1}} \\ \otimes\\ \vdots \\ \otimes \\ (\Sigma^{-1}C)^{\alpha_{\tilde{\beta}_{t,i}+\beta_{t,i}-1,1}}\end{array}\right]\otimes \cdots \otimes \left[\begin{array}{c}(\Sigma^{-1}C)^{\alpha_{\tilde{\beta}_{t,i},n}}\\ \otimes \\ \vdots \\ \otimes \\ (\Sigma^{-1}C)^{\alpha_{\tilde{\beta}_{t,i}+\beta_{t,i}-1,n}}\end{array}\right]\right\}\otimes\cdots}}} $$ \caption{Forming the composition $a\circ_i b$ in $\barcs{Z_n(\coend(C))}$ \textit{(continued)}}\label{fig:operadcomps2}\end{figure}

As before, we formulate this in terms that are well-defined for \emph{any} operad:

\begin{defn}
\label{def:cobaroperadcomps}Let \( \mathfrak{R} \) be an operad and define
composition-operations on direct summands of \( \barcs (Z_{*}(\mathfrak{R}) \)
as follows: 
\begin{enumerate}
\item let \( \alpha =\{\alpha _{1},\dots ,\alpha _{j}\} \) be a sequence of \( j \)
length-\( n \) sequences of nonnegative integers, and let \( a=[r_{1}|\dots |r_{j}]\in \barcs (Z_{n}(\mathfrak{R})) \)
with \( r_{v}\in \Sigma ^{-|\alpha _{v}|}\mathfrak{R}_{\alpha _{v}} \)
\item let \( \beta =\{\beta _{1},\dots ,\beta _{k}\} \) be a sequence of \( k \)
length-\( m \) sequences of integers and let \( b=[s_{1}|\dots |s_{k}]\in \barcs (Z_{m}(\mathfrak{R}) \)
with \( s_{w}\in \Sigma ^{-|\beta _{w}|}\mathfrak{R}_{\beta _{w}} \). 
\end{enumerate}
Then \( a\circ _{i}b \) is nonzero \emph{only if} \( j=\sum _{t=1}^{k}\beta _{t,i} \),
in which case
\[
a\circ _{i}b=S[r_{1}'|\dots |r_{k}']\]
 where
\begin{eqnarray*}
\lefteqn {r_{t}'=V'(\susp ^{|\tilde{\beta }_{t,i}+\beta _{t,i}-1|}r_{\tilde{\beta }_{t,i}+\beta _{t,i}-1})\circ _{\hat{\beta }_{t,i}+\beta _{t,i}-1}} &  & \\
 &  & \cdots (\susp ^{|\tilde{\beta }_{t,i}+1|}r_{\tilde{\beta }_{t,i}+1})\circ _{\hat{\beta }_{t,i}+1}(\susp ^{|\tilde{\beta }|}r_{\tilde{\beta }_{i,t}})\circ _{\hat{\beta }_{t,i}}(\susp ^{|\beta _{w}|}s_{t})
\end{eqnarray*}
 and:
\begin{enumerate}
\item \( S=(-1)^{AB} \), where \( A=\sum _{t=1}^{n}|\alpha _{t}| \) and \( B=\sum _{t=1}^{m}\sum _{z=1}^{i-1}\beta _{t,z} \)
\item \( \hat{\beta }_{t,i}=\sum _{v=1}^{i-1}\beta _{t,v} \) and \( \tilde{\beta }_{t,i}=\sum _{v=1}^{t-1}\beta _{v,i} \)
\item \( V'=\desusp ^{A}\mathscr{Z}\{\lambda _{1},\dots ,\lambda _{\beta _{t,i}}\} \),
where the \( \{\lambda _{\mu }\} \) are \( \beta _{t,i} \) length-\( n \)
sequences of integers defined by \( \lambda _{u,v}=\alpha _{\tilde{\beta }_{t,i}+u-1,v} \) 
\end{enumerate}
\end{defn}
\begin{rem}
These compositions represent a \emph{direct translation} of figure~\ref{fig:operadcomps2}
into formal composition operations. For instance, 
\begin{enumerate}
\item \( \hat{\beta }_{t,i} \) \emph{counts} the number of copies of \( \Sigma ^{-1}C \)
that \emph{precedes} the \( t^{\mathrm{th}} \) column (where the maps defining
\( a \) get inserted).
\item \( \beta _{v,i} \) is the number of copies of \( \Sigma ^{-1}C \) that are
applied to the \( \nu ^{\mathrm{th}} \) row of the map defining \( b \), which
means that
\item \( \tilde{\beta }_{t,i} \) is the number of the \emph{first factor} of \( \Sigma ^{-1}C \)
in the map defining \( a \) that is \emph{available} to be applied to the \( t^{\mathrm{th}} \)
row of the map defining \( b \) (in other words, it totals the number of copies
of \( \Sigma ^{-1}C \) that were ``used up'' by rows \( 1 \) through \( t-1 \)).
\end{enumerate}
As such, these compositions define the operad structure of \( \coend (\cobar C) \)
in terms of the operad structure of \( \coend (C) \).
\end{rem}
\begin{prop}
\label{prop:cobarstructure-appb}Let \( C \) be an m-coalgebra over the \( E_{\infty } \)-operad
\( \mathfrak{R} \). In addition, let 
\[
f_{n}:\barcs (Z_{n}(\mathfrak{R}))\rightarrow \coend (\cobar C)_{n}\]
 for all \( n>1 \) be the chain-maps induced by the structure map of \( C \).
If \( \barcs (Z_{n}(\mathfrak{R})) \) is equipped with the composition-operations
defined in \ref{def:cobaroperadcomps}, then
\begin{enumerate}
\item \( \barcs (Z_{n}(\mathfrak{R})) \) is an operad
\item the maps \( \{f_{n}\} \) constitute a morphism of operads
\end{enumerate}
\end{prop}
\begin{proof}
The second statement is the easiest to prove --- bearing in mind the remark
following \ref{def:cobaroperadcomps}. The first statement (that \( \barcs (Z_{n}(\mathfrak{R})) \),
with these compositions, is an operad) requires us to verify the defining identities
of an operad in \ref{def:operad}). 

We use the hypothesis and \ref{prop:everyoperad} (and \ref{claim:cobarnegative})
to conclude that there exists a coalgebra, \( D \), over \( \mathfrak{R} \)
with an injective structure map
\[
\mathfrak{R}\rightarrow \coend (D)\]
It is not hard to see that the induced map
\[
\barcs (Z_{n}(\mathfrak{R})\rightarrow \coend (\cobar D)\]
 will also be injective. Since (by statement 2 in the conclusions) we know that
the composition-operations of \( \barcs (Z_{n}(\mathfrak{R}) \) map to those
of \( \coend (\cobar D) \) and since the \emph{latter} satisfy the operad identities,
it follows that the \emph{former} do as well.
\end{proof}
Now we will determine corresponding compositions on \( Y_{n}(\mathfrak{R})\bigstar _{\rho _{n}}\barcs (Z_{n}(\mathfrak{R})) \)
to make the map \( g^{*} \) in \ref{cor:bigstarmaps} a morphism of operads.
As before, we are guided by the behavior of compositions in \( \coend (C\otimes _{\lambda }\cobar C) \).
Several features of these compositions come to mind:

Let \( u\in \coend (C\otimes _{\lambda }\cobar C) \) be in the image of 
\[
y\bigstar z\in Y_{n}(\coend (C))\bigstar _{\rho _{n}}\barcs (Z_{n}(\coend (C)))\]
Then \( y=u|C\otimes 1\subset C\otimes _{\lambda }\cobar C \) and \( z=u|1\otimes \cobar C\subset C\otimes _{\lambda }\cobar C \)
--- since \( C\otimes _{\lambda }\cobar C \) is a free \( \cobar C \)-module.

Given \( y_{1},y_{2}\in Y_{n}(\coend (C)) \) and \( z_{1},z_{2}\in \barcs (Z_{n}(\coend (C))) \),
let 
\begin{equation}
\label{eq:appbtensorcomp}
y\bigstar z=(y_{1}\bigstar z_{1})\circ _{i}(y_{2}\bigstar z_{2})
\end{equation}

\begin{claim}
In equation~\ref{eq:appbtensorcomp}, \( z=z_{1}\circ _{i}z_{2} \).
\end{claim}
\begin{proof}
Note that 
\begin{eqnarray*}
(y_{1}\bigstar z_{1})\circ _{i}(y_{2}\bigstar z_{2})|1\otimes \cobar C & = & (y_{1}\bigstar z_{1})\circ _{i}(y_{2}\bigstar z_{2}|1\otimes \cobar C)\\
 & = & (y_{1}\bigstar z_{1})\circ _{i}1\otimes z_{2}\\
 & = & 1\otimes (z_{1}\circ _{i}z_{2})
\end{eqnarray*}
since \( \cobar C\subset C\otimes _{\lambda }\cobar C \) so that the factor
\( \barcs (Z_{n}(\coend (C))) \) in \( Y_{n}(\coend (C))\bigstar _{\rho _{n}}\barcs (Z_{n}(\coend (C))) \)
is \emph{closed} under compositions.
\end{proof}
Consequently, the only remaining problem is to compute
\[
y\bigstar z|C\otimes 1=(y_{1}\bigstar z_{1})\circ _{i}y_{2}\]
since
\[
y\bigstar z=(y_{1}\bigstar z_{1})\circ _{i}y_{2}\cdot (z_{1}\circ _{i}z_{2})\]

\begin{figure} \begin{xy}0;<1cm,0cm>: (0,0)*+{C}; (5.4,0)*+{\left\{\begin{array}{c} C\otimes(\Sigma^{-1}C)^{\gamma_{1}}\\ \otimes\\\vdots\\ \otimes\\ C\otimes(\Sigma^{-1}C)^{\gamma_{i}}\\\otimes\\\vdots\\ \otimes\\C\otimes(\Sigma^{-1}C)^{\gamma_{k}}\end{array}\right.}**\dir{-} ?<>(.4)*_!/10pt/{y_2}?>*\dir{>}; 
(2,3.4)*+\txt{$i^{\mathrm{th}}$ \textit{row} \\ corresponds to the \\$i^{\mathrm{th}}$ copy of $C\otimes_\lambda\cobar(C)$\\ in target, $(C\otimes_\lambda\cobar(C))^k$}; (4.5,.4)*+{} **\dir{-} ?>*\dir{>};
(4.58,0)*+{\hbox to .15cm{\vbox to .4cm{}}} *\frm<8pt>{-}; 
(.1,-3)*+{U\otimes C\otimes_\lambda (\Sigma^{-1}C)^{\alpha_1}\otimes\cdots\otimes C\otimes_\lambda (\Sigma^{-1}C)^{\alpha_j}}*\frm<8pt>{-} **\dir{-} ?<>(.5)*^!/15pt/{1_U\otimes y_1}?>*\dir{>};
(5.9,0)*+{\hbox to 1.4cm{\vbox to .4cm{}}} *\frm<8pt>{-};(6.5,-3)*+{(\Sigma^{-1}C)^{\bar{\beta}_1}\otimes\cdots\otimes (\Sigma^{-1}C)^{\bar{\beta}_j}\otimes V} *\frm<8pt>{-} **\crv ~c{ (8,-3)&(6,-3.5)} ?<>(.3)*_!/20pt/{z_1\otimes 1_V}?>*\dir{>}; (3.73,-3.05)*+{\bigotimes};
(3,-5)*+{U\otimes C\otimes_{\lambda}(\Sigma^{-1 }C)^{\alpha_1+\bar{\beta}_1}\otimes\cdots \otimes C\otimes_{\lambda}(\Sigma^{-1}C)^{\alpha_j+\bar{\beta}_j}\otimes V}**\dir{-} ?<>(.4)*_!/20pt/{\mathrm{shuffle}}?>*\dir{>}; 
(3,-5)*+{}; (3,-7)*+{(C\otimes_{\lambda}\cobar {C})^{j+k-1}}; 
(3,-6.7)*+{}\ar @{^(->}(3,-5.3)*+{} ; 
\end{xy}\caption{Forming the composition $(y_1\bigstar z_1)\circ_i y_2$ }\label{fig:ycomp}\end{figure}

We will compute this on the canonical summands of the factors. Let: 

\begin{enumerate}
\item \( \gamma =\{\gamma _{1},\dots ,\gamma _{k}\} \) be a sequence of nonnegative
integers and let \( y_{2}\in \Sigma ^{-|\gamma |}\coend (C)_{\{\gamma _{1}',\dots ,\gamma _{k}'\}}\subset Y_{k}(\coend (C)) \)
(recalling the notation \( \gamma _{\ell }'=\gamma _{\ell }+1 \) and \( |\gamma |=\sum \gamma _{\ell } \)
--- see \ref{def:zn})
\item \( \beta =\{\beta _{1},\dots ,\beta _{r}\} \) be a sequence of length-\( j \)
sequences of nonnegative integers (i.e., a two-dimensional array) and let 
\[
z_{1}=[v_{1}|\dots |v_{r}]\in \bigotimes _{t=1}^{r}\Sigma \Sigma ^{-|\beta _{t}|}\coend (C)_{\beta _{t}}\subset \barcs Z_{j}(\coend (C))\]

\item \( \alpha =\{\alpha _{1},\dots ,\alpha _{j}\} \) be a sequence of nonnegative
integers and let \( y_{1}\in \Sigma ^{-|\alpha |}\coend (C)_{\{\alpha _{1}',\dots ,\alpha _{j}'\}}\subset Y_{j}(\coend (C)) \)
\end{enumerate}
Then figure~\ref{fig:ycomp} depicts the maps that enter into \( (y_{1}\bigstar z_{1})\circ _{1}y_{2} \),
where: 

\begin{enumerate}
\item \( \betabar _{i}=\sum _{\ell =1}^{r}\beta _{r,i} \) --- i.e., the summation
by columns. 
\item \( U=\bigotimes _{\ell =1}^{i-1}C\otimes _{\lambda }(\Sigma ^{-1}C)^{\gamma _{\ell }} \)
\item \( V=\bigotimes _{\ell =i+1}^{k}C\otimes _{\lambda }(\Sigma ^{-1}C)^{\gamma _{\ell }} \)
\end{enumerate}
Examination of figure~\ref{fig:ycomp} shows that we must have \( \gamma _{i}=r \)
(so that the copies of \( \Sigma ^{-1}C \) match up with the bars in \( [v_{1}|\dots |v_{r}] \)).
We may write the compositions of maps in terms of the composition operations
of \( \coend (C) \) and the shuffle permutation at the bottom of figure~\ref{fig:ycomp}
translates into the action of a suitable element of a symmetric group. We get

\begin{prop}
\label{prop:ycompositions-appb}Let \( \mathfrak{R} \) be an operad, let \( j \),
\( k \), \( i \), and \( r \) be nonnegative integers with \( i\leq k \).
In addition, let
\begin{enumerate}
\item \( \gamma =\{\gamma _{1},\dots ,\gamma _{k}\} \) be a sequence of nonnegative
integers.
\item \( \beta =\{\beta _{1},\dots ,\beta _{r}\} \) be a sequence of length-\( j \)
sequences of nonnegative integers (i.e., a two-dimensional array).
\item \( \alpha =\{\alpha _{1},\dots ,\alpha _{j}\} \) be a sequence of nonnegative
integers.
\end{enumerate}
Define composition-operations
\[
Y_{j}(\mathfrak{R})\bigstar _{\rho _{j}}\barcs (Z_{j}(\mathfrak{R})\circ _{i}Y_{k}(\mathfrak{R})\]
 on canonical direct summands

\begin{enumerate}
\item \( Y_{1}=\Sigma ^{-|\alpha |}\mathfrak{R}_{\{\alpha _{1}',\dots ,\alpha _{j}'\}}\subset Y_{j}(\mathfrak{R}) \)
\item \( Z_{1}=[v_{1}|\dots |v_{r}]\in \bigotimes _{t=1}^{r}\Sigma \Sigma ^{-|\beta _{t}|}\mathfrak{R}_{\beta _{t}}\subset \barcs Z_{j}(\mathfrak{R}) \)
\item \( Y_{2}=\Sigma ^{-|\gamma |}\mathfrak{R}_{\{\gamma _{1}',\dots ,\gamma _{k}'\}}\subset Y_{k}(\mathfrak{R}) \)
\end{enumerate}
via

\begin{itemize}
\item if \( r\neq \gamma _{i} \) then \( (y_{1}\bigstar z_{1})\circ _{i}y_{2}=0 \)
\item otherwise
\begin{eqnarray*}
(y_{1}\bigstar z_{1})\circ _{i}y_{2} & = & \desusp ^{|\alpha |+|\gamma |+\sum \betabar _{*}}\\
 & \circ  & \mathscr{Z}(\susp ^{|\alpha |}y_{1})\circ _{\mu }(\susp ^{\betabar _{1}}r_{1})\circ _{\mu +1}\cdots \circ _{\mu +j-1}(\susp ^{\betabar _{j}})\circ _{\mu +j}\susp ^{|\gamma |}y_{2}
\end{eqnarray*}
where:
\end{itemize}
\begin{enumerate}
\item \( \mu =i+\sum _{\ell =1}^{i-1}\gamma _{\ell } \)
\item \( \mathscr{Z}=\mathscr{Z}\{(1,\alpha _{1},1,\alpha _{2},\dots ,1,\alpha _{j}),(0,\betabar _{1},\dots ,0,\betabar _{j})\} \)
\end{enumerate}
for all \( y_{1}\in Y_{1} \), \( z_{1}\in Z_{1} \), and \( y_{2}\in Y_{2} \).
Extend these compositions to all of 
\[
Y_{j}(\mathfrak{R})\bigstar _{\rho _{j}}\barcs (Z_{j}(\mathfrak{R})\]
by setting 
\begin{equation}
\label{eq:app-b-operadextend}
(y_{1}\bigstar z_{1})\bigstar (y_{2}\bigstar z_{2})=((y_{1}\bigstar z_{1})\circ _{i}y_{2})\cdot (z_{1}\circ _{i}z_{2})
\end{equation}
for all \( y_{1},y_{2}\in Y_{j}(\mathfrak{R}) \) and \( z_{1},z_{2}\in \barcs (Z_{j}(\mathfrak{R}) \)

Then the \( \{\circ _{i}\} \) define the operad structure of 
\[
Y_{j}(\mathfrak{R})\bigstar _{\rho _{j}}\barcs (Z_{j}(\mathfrak{R})\circ _{i}Y_{k}(\mathfrak{R})\]
in the case where \( \mathfrak{R}=\coend (C) \).

\end{prop}
\begin{rem}
This follows directly from figure~\ref{fig:ycomp}. We have simply translated 
\begin{enumerate}
\item compositions of \emph{maps} into suitable \( \{\circ _{i}\} \) operations in
\( \coend (C) \) (counting the number of factors of \( \Sigma ^{-1}C \) to
the right of the factor we are interested in to determine the value of \( i \))
\item the shuffle \emph{map} becomes the action of a suitable shuffle \emph{permutation}
in the symmetric group.
\item The product map in equation~\ref{eq:app-b-operadextend} translates into a \emph{shuffle
map} (therefore the action of a shuffle \emph{permutation}) because \( \cobar C \)
is a free algebra.
\end{enumerate}
Note that we do not claim, given an \emph{arbitrary} operad \( \mathfrak{R} \),
that
\[
Y_{j}(\mathfrak{R})\bigstar _{\rho _{j}}\barcs (Z_{j}(\mathfrak{R})\]
equipped with these compositions, constitutes an operad. We only \emph{know}
that it is an operad when \( \mathfrak{R}=\coend (C) \).
\end{rem}
\begin{thm}
\label{th:ytensorstructure-appb}Let \( C \) be an m-coalgebra over an \( E_{\infty } \)-operad,
\( \mathfrak{R} \), and let 
\[
T=Y_{j}(\mathfrak{R})\bigstar _{\rho _{j}}\barcs (Z_{j}(\mathfrak{R})\]
 be equipped with the compositions defined in \ref{prop:ycompositions-appb}.
Then:
\begin{enumerate}
\item The set of chain-complexes \( T \) constitute an operad
\item the map
\[
Y_{j}(\mathfrak{R})\bigstar _{\rho _{j}}\barcs (Z_{j}(\mathfrak{R})\rightarrow \coend (C\otimes _{\lambda }\cobar C)\]
defined in \ref{cor:bigstarmaps} is a morphism of operads.
\end{enumerate}
\end{thm}
\begin{proof}
As before, the second statement is the easy one to prove --- it follows directly
from \ref{prop:ycompositions-appb}. As in \ref{prop:cobarstructure-appb},
the first statement (that \( T \) is an operad follows from \ref{prop:everyoperad}
--- which implies that there exists a coalgebra, \( D \), over \( \mathfrak{R} \)
with an injective structure map
\[
\mathfrak{R}\rightarrow \coend (D)\]
It is not hard to see that the induced map
\[
Y_{j}(\mathfrak{R})\bigstar _{\rho _{j}}\barcs (Z_{j}(\mathfrak{R})\rightarrow \coend (D\otimes _{\lambda }\cobar D)\]
will also be injective. This requires a quick check that we have not used any
of the specific properties of m-coalgebras distinguishing them from arbitrary
coalgebras over an operad (i.e., the fact that they are concentrated in positive
dimensions). Since (by statement 2 in the conclusions) we know that the composition-operations
of \( Y_{j}(\mathfrak{R})\bigstar _{\rho _{j}}\barcs (Z_{j}(\mathfrak{R}) \)
map to those of \( \coend (D\otimes _{\lambda }\cobar D) \) and since the \emph{latter}
satisfy the operad identities, it follows that the \emph{former} do as well.
\end{proof}
\begin{prop}
\label{prop:ycontainsr}Under the hypotheses of \ref{prop:ycompositions-appb}
and \ref{th:ytensorstructure-appb}, the composition-operations of 
\[
Y_{j}(\mathfrak{R})\bigstar _{\rho _{j}}\barcs (Z_{j}(\mathfrak{R})\]
 restricted to the canonical summand \( Y_{0}=\Sigma ^{-|\alpha |}\mathfrak{R}_{\{\alpha _{1}',\dots ,\alpha _{j}'\}}\subset Y_{j}(\mathfrak{R}) \)
with all \( \alpha _{\ell }=0 \) coincide with the composition-operations of
\( \mathfrak{R}_{j}=Y_{0} \).
\end{prop}
\begin{proof}
Simply set all \( \gamma _{\ell }=0 \) in the formula in \ref{prop:ycompositions-appb}.
\end{proof}

\section{Concordance between the current notation and that of \cite{Smith:1994}\label{app:concordance}}

{\par\centering \vspace{0.3cm} \par}

{\centering \begin{tabular}{|c|c|}
\hline 
Notation in \cite{Smith:1994}&
Current\\
\hline 
formal coalgebras&
 non-\( \Sigma  \) operad\\
\hline 
symmetric formal coalgebra&
 operad\\
\hline 
f-resolution&
 \( E_{\infty } \)-operad\\
\hline 
weakly-coherent m-coalgebra&
 m-coalgebra\\
\hline 
coherent m-coalgebra&
 m-coalgebra over \( \mathfrak{S} \) \\
\hline 
\( \cobar (\mathfrak{R},n) \)&
\( L(\mathfrak{R}) \)\\
\hline 
\end{tabular}\par}

{\par\centering \par{} \vspace{0.3cm}\par}

\providecommand{\bysame}{\leavevmode\hbox to3em{\hrulefill}\thinspace}


    \end{document}